# Enveloping operads and applications

Victor Carmona*

October 2025

**Abstract**

This work addresses the homotopical analysis of enveloping operads in a general cofibrantly generated symmetric monoidal model category. We show the potential of this analysis by obtaining, in a uniform way, several central results regarding the homotopy theory of operadic algebras.

# 1 Introduction

The *universal enveloping algebra* of a Lie-algebra $\mathfrak{g}$ is a unital associative algebra $\big(\mathsf{U}(\mathfrak{g}), \star, 1\big)$ characterized by the fact that its category of left-modules is equivalent to the category of $\mathfrak{g}$-representations. Alternatively, it is the initial such algebra equipped with a map of Lie-algebras $\big(\mathfrak{g}, [\text{-},\text{-}]\big) \to \big(\mathsf{U}(\mathfrak{g}), \star - \star^{\mathsf{op}}\big)$. This construction, which has found important applications in mathematics and mathematical physics, is just the unary part of a much more general device: the *enveloping operad* $\mathsf{Lie}_{\mathfrak{g}}$. This object can be also characterized by a simple fact: $\mathsf{Lie}_{\mathfrak{g}}$-algebras are the same as Lie-algebras under $\mathfrak{g}$ (see [10]). Alternatively, it is the initial operad equipped with a map of operads $\mathsf{Lie} \to \mathsf{Lie}_{\mathfrak{g}}$ and a map of Lie-algebras $\mathfrak{g} \to \mathsf{Lie}_{\mathfrak{g}}(0)$. Notice that $\mathsf{Lie}_{\mathfrak{g}}$ contains much more information than $\mathsf{U}(\mathfrak{g})$. For instance, the *isomorphism problem* for Lie-algebras (see [37]) is trivial if we consider $\mathsf{Lie}_{\mathfrak{g}}$ instead of $\mathsf{U}(\mathfrak{g})$.

For more general operads $\mathsf{Lie} \rightsquigarrow \mathcal{O}$ and their associated algebras $\mathfrak{g} \rightsquigarrow \mathsf{A}$, the analogous *enveloping operad* $\mathcal{O}_{\mathsf{A}}$, the one that controls $\mathcal{O}$-algebras under $\mathsf{A}$, can also be considered. Again, its unary part $\mathsf{U}_{\mathcal{O}}(\mathsf{A})$ yields a unital associate algebra (or category in the colored case) whose category of left-modules is equivalent to the category of operadic modules over $\mathsf{A}$ (see [5]).

In this document, we analyze the homotopical properties of the two-variable assignment $(\mathcal{O}, \mathsf{A}) \mapsto \mathcal{O}_{\mathsf{A}}$ for colored enriched operads, e.g. we study when such a rule preserves weak equivalences in its variables. Due to a new graphical point of view developed in

*The author was partially supported by the grant PID2020-117971GB-C21 funded by MCIN/AEI/10.13039/501100011033.

*Affiliation:* Max-Planck Institut für Mathematik in den Naturwissenschaften, Leipzig, Germany.
*Email address:* `victor.carmona@mis.mpg.de`

§3, inspired by work of Muro, we are able to improve previous results in this respect and to explain in a more geometrical way convoluted constructions in the literature. Our main original results regarding the homotopical analysis of $(\mathcal{O}, \mathsf{A}) \mapsto \mathcal{O}_{\mathsf{A}}$ concern the restriction of this assignment to cofibrant operads. In fact, they are generalizations of several results in [29] (dealing with one-colored non-symmetric operads) and of [10, Theorem 17.4.B] (dealing with one-colored dg-operads)[i] to colored symmetric operads in a general cofibrantly generated symmetric monoidal model category. For a summary of these results, see the end of §3.3.

Despite being abstract computations, the second part of this work shows how fundamental enveloping operads are in the homotopy theory of operadic algebras by applying the homotopical study of $(\mathcal{O}, \mathsf{A}) \mapsto \mathcal{O}_{\mathsf{A}}$ to uniformly obtain a variety of consequences. Actually, we get novel (resp. recover and sometimes fix classical) results regarding:

- **Admissibility and rectification of operads:** (4.2, 4.7, 4.9 and 4.16) existence of (semi)model structures on categories of operadic algebras and the study of when those are Quillen equivalent.
- **Change of homotopy cosmos:** (4.23 and 4.25) comparison of operadic algebras along Quillen equivalences of the underlying symmetric monoidal model categories.
- **Relative left properness of operadic algebras:** (4.29, 4.36 and 4.38) stability of weak equivalences under pushouts along cofibrations in categories of algebras.
- **Moduli of algebra structures:** (4.46) coincidence of the two competing "moduli spaces of derived algebra structures" over an operad.
- **Homotopy invariance of module categories over operadic algebras:** (4.49 and 4.50) study of when assigning modules over an operadic algebra sends weak equivalences to Quillen equivalences.
- **Coincidence of operadic and categorical left Kan extensions:** (4.52)

In most of these applications, either we broadly generalize previous results or we are able to drop hypotheses in known statements. This is extensively justified at the end of each subsection in §4.

We recommend the reader only interested in applications to not follow a linear order when reading this document. Instead, consult directly the subsection containing your result of interest and trace back to the technical body of the paper if needed. On the other hand, this paper is mostly self-contained and it has been written following a clear linear order, thus it is suitable to also read it in this order.

**Relation to other work.** Operads and algebras have been extensively studied in the literature and there are wonderful documents where one can learn about this subject. We want to name a few of them: [10, 20, 22, 23, 38, 39].

[i]Observe that it is claimed a more general statement in [10], but the proof is restricted to the dg-case. The result turns out to be more complicated outside the dg-setting, and actually his statement must be corrected in general.

The study of enveloping operads and their relevance for the above list of topics, related to the homotopy theory of operadic algebras, comes back to the work of various people, for instance: [3, 4, 5, 10, 14, 15, 16, 17, 25, 26, 29, 30, 40, 41]. Of course, this list is far from exhaustive.

Due to the large list of results obtained in this work, more detailed and specific explanations about connections with existing literature are given at the end of each subsection in §4.

**Organization of the paper.** The document begins with a preliminary section, §2, which quickly revisits the basics of operads and operadic algebras. Section §3 contains the technical core of the paper: the homotopical analysis of the assignment $(\mathcal{O}, \mathsf{A}) \mapsto \mathcal{O}_{\mathsf{A}}$. It starts by defining the graphical machinery used to handle enveloping operads and closes with a summary of the obtained results. The next part, §4, is devoted to applications of our previous homotopical analysis. Each subsection can be consulted independently, but their order in the document can be seen as a suggested logical order. Finally, appendix §A gathers notation regarding semimodel categories, cell complexes and technical observations about cell-attachments and $\mathbb{I}$-cofibrant objects.

**Notation 1.1.** In the sequel we adopt the following notation: symbols like $\mathtt{V}, \mathtt{M}\ldots$ will denote ordinary categories, while $\mathcal{V}, \mathcal{M}\ldots$ will be reserved for (semi)model categories.

# Contents

## 2 Preliminaries

Operads can be presented in many equivalent ways, but the idea that one should keep in mind is that they are a generalization of categories (or monoids in the one-colored case) where one allows for multimorphisms, i.e. morphisms with a finite list of inputs and one output. In this work, the term "operad" will stand for (enriched) symmetric colored operad or (enriched) symmetric multicategory.

Fix a symmetric monoidal category (from now on sm-category) $\mathtt{V}$ which will be the base for enrichment. For simplicity, we assume once and for all that $\mathtt{V}$ is closed and bicomplete.

The paradigmatic example of (one-colored) operad is the *endomorphism operad* of an object $X \in \mathtt{V}$. In the same way in which one axiomatizes what a monoid is by looking at the algebraic structure on $\underline{\mathrm{Hom}}(X, X)$ coming from the composition product, one can deduce what an operad is by studying the family $\{\underline{\mathrm{Hom}}(X^{\otimes n}, X)\}_{n \geqslant 0}$ equipped with the various compositions that it admits. The notion of colored operad is what one gets if instead of an object, one looks at a collection of objects in $\mathtt{V}$.

**Definition 2.1.** An *operad* $\mathcal{O}$ in $\mathtt{V}$, or $\mathtt{V}$-operad, consists of:

- a collection of colors or objects $\mathrm{col}(\mathcal{O})$;
- a $\mathtt{V}$-object of multimorphisms $\mathcal{O}\begin{bmatrix}\{c_i\}_i \\ b\end{bmatrix}$ for any finite collection $\{c_i\}_i$ in $\mathrm{col}(\mathcal{O})$ (possibly with repetitions) and any color $b \in \mathrm{col}(\mathcal{O})$;
- an identity element $\mathrm{id}_c \colon \mathbb{1}_{\mathtt{V}} \to \mathcal{O}\begin{bmatrix}\{c\} \\ c\end{bmatrix}$ for any color $c \in \mathrm{col}(\mathcal{O})$;
- a composition product

$$\mu\colon\ \mathcal{O}\begin{bmatrix}\{b_j\}_{j\in J} \\ d\end{bmatrix} \otimes \left(\bigotimes_{j\in J} \mathcal{O}\begin{bmatrix}\{c_i\}_{i\in f^{-1}(j)} \\ b_j\end{bmatrix}\right) \longrightarrow \mathcal{O}\begin{bmatrix}\{c_i\}_{i\in I} \\ d\end{bmatrix}$$

  for any map of finite sets $f\colon I \to J$.

These data are required to satisfy unitality and associativity conditions.

A *morphism of* $\mathtt{V}$*-operads* $\varphi\colon \mathcal{O} \to \mathcal{P}$ is given by a function $\varphi\colon \mathrm{col}(\mathcal{O}) \to \mathrm{col}(\mathcal{P})$ together with maps $\varphi\colon \mathcal{O}\begin{bmatrix}\{c_i\}_i \\ d\end{bmatrix} \to \mathcal{P}\begin{bmatrix}\{\varphi(c_i)\}_i \\ \varphi(d)\end{bmatrix}$ compatible with compositions and units.

We denote by $\mathtt{Opd}(\mathtt{V})$ the category of $\mathtt{V}$-operads and their morphisms.

*Remark* 2.2. Note that Definition 2.1 encodes the action of $\mathrm{Aut}(I)$ on $\mathcal{O}\begin{bmatrix}\{c_i\}_i \\ d\end{bmatrix}$ for any finite sequence $\{c_i\}_{i\in I}$ indexed by $I$ and all the associated compatibilities. Also, observe that we are implicitly assuming that tensor products indexed by a finite set are unordered (see [22] or [24]). We will assume this convention whenever it is required to simplify notation and the exposition.

One of the main motivations to define operads is to study their representations, also called algebras. Roughly speaking, an operad parametrizes operations of arities $n \geqslant 0$ with their relations and an algebra over an operad consists of an object, or a family of objects, which carries those operations.

**Definition 2.3.** An $\mathcal{O}$-*algebra* A is given by the following data:

- a $\mathrm{col}(\mathcal{O})$-family of objects $\{A(o)\}_o$, and
- action maps $\mu\colon \mathcal{O}\left[{}^{\{c_i\}_i}_{\ d}\right] \otimes \bigotimes_i A(c_i) \to A(d)$,

subject to associativity, unitality and equivariance axioms [43, §13]. A morphism of $\mathcal{O}$-algebras is a family of morphisms $A(o) \to B(o)$ preserving this structure.

The category of $\mathcal{O}$-algebras (in $\mathtt{V}$) and their morphisms will be denoted by $\mathtt{Alg}_{\mathcal{O}}(\mathtt{V})$, or simply $\mathtt{Alg}_{\mathcal{O}}$.

*Remark* 2.4. Of course, it is possible to change the underlying base of enrichment for algebras and operads. Given a lax (symmetric) monoidal functor $|\star|\colon \mathtt{V} \to \mathtt{V}'$, applying it locally or objectwise, one gets base-change functors

$$|\star|\colon \mathtt{Opd}(\mathtt{V}) \to \mathtt{Opd}(\mathtt{V}') \quad \text{and} \quad \mathtt{Alg}_{\mathcal{O}}(\mathtt{V}) \to \mathtt{Alg}_{|\mathcal{O}|}(\mathtt{V}').$$

In particular, $\mathtt{V}(\mathbb{I},\star)\colon \mathtt{V} \to \mathsf{Set}$ is lax monoidal and hence induces

$$(\star)_0\colon \mathtt{Opd}(\mathtt{V}) \to \mathtt{Opd}(\mathsf{Set}).$$

Going one step ahead, if $\mathcal{V}$ is a (closed) sm-model category, the localization functor $\mathcal{V} \to \mathrm{Ho}\,\mathcal{V}$ is lax monoidal and so one gets two important base-change functors

$$\pi\colon \mathtt{Opd}(\mathcal{V}) \to \mathtt{Opd}(\mathrm{Ho}\,\mathcal{V}) \quad \text{and} \quad \pi_0\colon \mathtt{Opd}(\mathcal{V}) \to \mathtt{Opd}(\mathrm{Ho}\,\mathcal{V}) \to \mathtt{Opd}(\mathsf{Set}).$$

**Categorical perspectives.** Let us collect some important facts about the categories $\mathtt{Opd}(\mathtt{V})$ and $\mathtt{Alg}_{\mathcal{O}}(\mathtt{V})$. Consult [10, 12] for complete details in the uncolored case.

Given an operad $\mathcal{O}$ in $\mathtt{V}$, one can produce an *analytic* (equivalently *Schur functor*) $\mathcal{O}\circ\star\colon \mathtt{V}^{\mathsf{O}} \to \mathtt{V}^{\mathsf{O}}$, where $\mathsf{O} = \mathrm{col}(\mathcal{O})$. On objects, this functor is given by

$$X = \{X(o)\}_{o\in \mathsf{O}} \longmapsto \mathcal{O}\circ X = \left\{ \coprod_{\{c_i\}_{i\in I}} \mathcal{O}\begin{bmatrix}\{c_i\}_i \\ o\end{bmatrix} \underset{\mathrm{Aut}(I)}{\otimes} \bigotimes_i X(c_i) \right\}_{o\in \mathsf{O}}.$$

By its very definition, this analytic functor preserves sifted colimits and moreover, the operad structure on $\mathcal{O}$ equips $\mathcal{O}\circ\star$ with the structure of a (sifted) monad. Algebras over $\mathcal{O}\circ\star$ are precisely the $\mathcal{O}$-algebras and hence, we have a free-forgetful adjunction

$$F_{\mathcal{O}}\colon \mathtt{V}^{\mathsf{O}} \rightleftarrows \mathtt{Alg}_{\mathcal{O}}(\mathtt{V}) : \mathrm{fgt}$$

with $\mathcal{O}\circ\star = \mathrm{fgt}\cdot F_{\mathcal{O}}$. That is, $\mathtt{Alg}_{\mathsf{O}}(\mathtt{V})$ is monadic over $\mathtt{V}^{\mathsf{O}}$. This adjunction can be exploited to show that $\mathtt{Alg}_{\mathcal{O}}(\mathtt{V})$ is bicomplete.

To any map of operads $\varphi\colon \mathcal{O} \to \mathcal{P}$, we can associate an obvious restriction functor

$$\varphi^*\colon \mathtt{Alg}_{\mathcal{P}}(\mathtt{V}) \longrightarrow \mathtt{Alg}_{\mathcal{O}}(\mathtt{V})$$

given by precomposition with $\varphi$. Either directly (e.g. [10, §3.3]) or via [7, Theorem 4.5.6], one shows that $\varphi^*$ admits a left adjoint $\varphi_\sharp$ called *operadic left Kan extension along* $\varphi$, or *operadic Lan along* $\varphi$.

*Example* 2.5. By neglecting all multimorphisms of arity different from $1$, one finds that any $\mathtt{V}$-operad $\mathcal{O}$ has an underlying $\mathtt{V}$-category $\overline{\mathcal{O}}$, which comes with a canonical inclusion map $\overline{\mathcal{O}} \hookrightarrow \mathcal{O}$. We will denote by

$$\mathsf{ext}\colon \mathtt{Alg}_{\overline{\mathcal{O}}}(\mathtt{V}) \rightleftarrows \mathtt{Alg}_{\mathcal{O}}(\mathtt{V}) : \overline{(\star)}$$

the associated adjunction given by operadic left Kan extension and restriction along $\overline{\mathcal{O}} \hookrightarrow \mathcal{O}$.

Interestingly, the category $\mathtt{Opd}_{\mathsf{O}}(\mathtt{V})$ of O-colored operads is the category of algebras over an operad and so, all the results above apply to it, in particular $\mathtt{Opd}_{\mathsf{O}}(\mathtt{V})$ is bicomplete. Varying the set of colors, one gets a bifibration $\mathsf{col}\colon \mathtt{Opd}(\mathtt{V}) \to \mathsf{Set}$ with fibers the categories $\mathtt{Opd}_{\mathsf{O}}(\mathtt{V})$. Thus, one observes that the whole $\mathtt{Opd}(\mathtt{V})$ is bicomplete. Despite being theoretically enough for our purposes, this paragraph is a bit too much abstract nonsense. For this reason, we recall a different perspective on $\mathtt{Opd}_{\mathsf{O}}(\mathtt{V})$: O-colored operads are monoids for the $\circ$-product.

We will fix a set of colors O for the following discussion.

To encode the structure of a O-colored collection of $\mathtt{V}$-objects with actions of symmetric groups, it is convenient to consider the groupoid $\mathsf{Fin}^{\simeq}_{\mathsf{O}} = \mathsf{Fin}^{\simeq}\downarrow \mathsf{O}$. It encodes finite sequences of colors and their permutations. For ease of notation, we will denote elements of this groupoid by underlined letters $\underline{c} = \{c_i\}_{i\in I}$.

**Definition 2.6.** The *groupoid of* O*-corollas* is $\mathsf{Fin}^{\simeq}_{\mathsf{O}} \times \mathsf{O}$ and we will denote its elements by symbols $\left[\begin{smallmatrix}\underline{c}\\ d\end{smallmatrix}\right]$ where $\underline{c} \in \mathsf{Fin}^{\simeq}_{\mathsf{O}}$ and $d \in \mathsf{O}$.

The *category of* O*-colored* $\Sigma$*-collections* (in $\mathtt{V}$), or $\Sigma_{\mathsf{O}}$*-collections*, is the diagram category $\Sigma\mathtt{Coll}_{\mathsf{O}}(\mathtt{V}) = [(\mathsf{Fin}^{\simeq}_{\mathsf{O}})^{\mathsf{op}} \times \mathsf{O}, \mathtt{V}]$.

*Remark* 2.7. We will freely pass from symmetric collections to symmetric sequences and vice versa [12, 19]. The difference is not substantial since there is an equivalence of categories

$$\Sigma\mathsf{Seq}_{\mathsf{O}}(\mathtt{V}) = [(\Sigma_{\mathsf{O}})^{\mathsf{op}} \times \mathsf{O}, \mathtt{V}] \simeq \left[(\mathsf{Fin}^{\simeq}_{\mathsf{O}})^{\mathsf{op}} \times \mathsf{O}, \mathtt{V}\right] = \Sigma\mathtt{Coll}_{\mathsf{O}}(\mathtt{V})$$

induced by restriction and left Kan extension along $(\Sigma_{\mathsf{O}})^{\mathsf{op}} \times \mathsf{O} \to (\mathsf{Fin}^{\simeq}_{\mathsf{O}})^{\mathsf{op}} \times \mathsf{O}$. We will apply this equivalence without further mention in the sequel.

It is clear that any O-colored operad $\mathcal{O}$ has an underlying $\Sigma_{\mathsf{O}}$-collection and that this assignment defines a forgetful functor $\mathtt{Opd}_{\mathsf{O}}(\mathtt{V}) \to \Sigma\mathtt{Coll}_{\mathsf{O}}(\mathtt{V})$. This functor admits a left adjoint, which can be constructed via operadic Lan, called *free operad* functor and we will denote the resulting adjunction by

$$\mathcal{F}\colon \Sigma\mathtt{Coll}_{\mathsf{O}}(\mathtt{V}) \rightleftarrows \mathtt{Opd}_{\mathsf{O}}(\mathtt{V}) : \mathsf{fgt}_{\Sigma}.$$

To any $\mathcal{X} \in \Sigma\mathtt{Coll}_{\mathsf{O}}(\mathtt{V})$, one can associate an analytic functor on $\mathtt{V}^{\mathsf{O}}$. This construction yields a functor $\Sigma\mathtt{Coll}_{\mathsf{O}}(\mathtt{V}) \to \mathsf{EndFun}(\mathtt{V}^{\mathsf{O}})$. Any category of endofunctors comes with a (non-symmetric) monoidal structure given by composition of functors and a monad is just a monoid for this monoidal structure. The functor $\Sigma\mathtt{Coll}_{\mathsf{O}}(\mathtt{V}) \to \mathsf{EndFun}(\mathtt{V}^{\mathsf{O}})$ can be upgraded to a (strong) monoidal functor by equipping $\Sigma\mathtt{Coll}_{\mathsf{O}}(\mathtt{V})$ with the $\circ$-product (see [10, 21]). Finally, O-colored operads are monoids in $(\Sigma\mathtt{Coll}_{\mathsf{O}}(\mathtt{V}), \circ)$, and that is why the Schur functor $\mathcal{O} \circ \star$ associated to an operad $\mathcal{O}$ gives a monad on $\mathtt{V}^{\mathsf{O}}$.

**Notation 2.8.** When $\varphi\colon \mathcal{O} \to \mathcal{P}$ is the identity on colors, it is common to denote $\varphi_\sharp(A)$ by $\mathcal{P} \circ_{\mathcal{O}} A$ and call it *relative $\circ$-product* by analogy with the theory of modules over rings.

**Notation 2.9.** There are two useful operations in $\mathrm{Fin}^{\cong}_{O}$ that deserve a name:

$$\underline{a} \boxplus \underline{b} = \{a_k\}_{k\in K} \amalg \{b_j\}_{j\in J} \quad \text{and} \quad \underline{b} \circ_t \underline{c} = \{b_j\}_{j\in J\setminus\{t\}} \amalg \{c_i\}_{i\in I}.$$

**Homotopy theory of operadic algebras** The most structured way to study a homotopy theory on $\mathtt{Alg}_{\mathcal{O}}(\mathtt{V})$ is by endowing this category with a (semi)model structure. A natural candidate is to study the left-transferred model structure on $\mathcal{O}$-algebras, if it exists, coming from the adjunction

$$F_{\mathcal{O}}\colon \mathtt{V}^{O} \rightleftarrows \mathtt{Alg}_{\mathcal{O}}(\mathtt{V}) : \mathrm{fgt}.$$

To have a structure to transfer, let us fix once and for all a closed cofibrantly generated sm-model category structure $\mathcal{V}$ on $\mathtt{V}$. One should think of $\mathcal{V}$ as a *homotopy cosmos* where all the homotopy theory occurs.

**Definition 2.10.** When it exists, the left-transferred model structure on $\mathtt{Alg}_{\mathcal{O}}(\mathtt{V})$ is called *projective model structure* and we denote it $\mathcal{A}\mathrm{lg}_{\mathcal{O}}(\mathcal{V})$ or simply $\mathcal{A}\mathrm{lg}_{\mathcal{O}}$.

We say that $\mathcal{O}$ is *admissible* if the projective model structure $\mathcal{A}\mathrm{lg}_{\mathcal{O}}(\mathcal{V})$ exists. The operad $\mathcal{O}$ is *strongly admissible* if it is admissible and $\mathrm{fgt}\colon \mathcal{A}\mathrm{lg}_{\mathcal{O}}(\mathcal{V}) \to \mathcal{V}^{O}$ preserves core (acyclic) cofibrations (Definition A.1).

By their very definition, the adjunctions induced by maps of operads $\varphi\colon \mathcal{O} \to \mathcal{P}$ are Quillen pairs $\varphi_\sharp\colon \mathcal{A}\mathrm{lg}_{\mathcal{O}}(\mathcal{V}) \rightleftarrows \mathcal{A}\mathrm{lg}_{\mathcal{P}}(\mathcal{V}) : \varphi^*$ when $\mathcal{O}$ and $\mathcal{P}$ are admissible. However, it is not always true that an equivalence of operads (see Definition 2.11) induces a Quillen equivalence (*rectification*). In practice, it is convenient to work with operads satisfying this additional condition; see §4.1.

**Homotopy theory of operads.** Let us briefly comment some details about the homotopy theory of operads and set some notation.

Fixing the set of colors, $\mathtt{Opd}_{O}(\mathtt{V})$ is the category of algebras over an operad, and so one might apply the previous discussion to this category. This way, one obtains, when it exists, a model structure on $\mathtt{Opd}_{O}(\mathtt{V})$ denoted $\mathcal{O}\mathrm{pd}_{O}(\mathcal{V})$, or $\mathcal{O}\mathrm{pd}_{O}$. One might have transferred the model structure using the free operad-forgetful adjunction instead. In other words, considering the adjunction

$$\mathcal{F}\colon \Sigma\mathtt{Coll}_{O}(\mathtt{V}) \rightleftarrows \mathtt{Opd}_{O}(\mathtt{V}) : \mathrm{fgt}_{\Sigma},$$

where the left hand side always has the projective model structure as a diagram category. Even when the transferred notions on $\mathtt{Opd}_{O}(\mathtt{V})$ through $\mathcal{F} \dashv \mathrm{fgt}_{\Sigma}$ do not fulfill all the axioms of a model structure, they are quite important and deserve a name.

**Definition 2.11.** Let $\mathcal{O}$ be an O-colored operad and $\mathcal{O} \to \mathcal{P}$ be a map in $\mathtt{Opd}_{O}(\mathtt{V})$.

- $\mathcal{O} \to \mathcal{P}$ is an *equivalence* if $\mathcal{O}\left[\begin{smallmatrix}\underline{c}\\ d\end{smallmatrix}\right] \xrightarrow{\sim} \mathcal{P}\left[\begin{smallmatrix}\underline{c}\\ d\end{smallmatrix}\right]$ is an equivalence for any O-corolla.
- The map $\mathcal{O} \to \mathcal{P}$ is a *fibration* if $\mathcal{O}\left[\begin{smallmatrix}\underline{c}\\ d\end{smallmatrix}\right] \twoheadrightarrow \mathcal{P}\left[\begin{smallmatrix}\underline{c}\\ d\end{smallmatrix}\right]$ is a fibration for any O-corolla.

- The map $\mathcal{O} \to \mathcal{P}$ is a *$\Sigma$-cofibration* if $\mathrm{fgt}_\Sigma \mathcal{O} \to \mathrm{fgt}_\Sigma \mathcal{P}$ is a cofibration in $\Sigma\mathcal{C}\mathrm{oll}_O(\mathcal{V})$, the projective model structure on O-collections.

- The operad $\mathcal{O}$ is *$\Sigma$-cofibrant* if $\mathrm{fgt}_\Sigma \mathcal{O}$ is cofibrant and $\mathcal{O}$ is *well-pointed* if the map from the initial O-operad $\mathcal{I}_O \to \mathcal{O}$ is a $\Sigma$-cofibration.

*Remarks* 2.12.

(i) There is a chain of implications

$$\text{cofibrant operad} \Longrightarrow \text{well-pointed operad} \overset{(*)}{\Longrightarrow} \Sigma\text{-cofibrant operad},$$

where $(*)$ holds when the monoidal unit in $\mathcal{V}$ is cofibrant. The first implication was addressed in, for instance, [3] (see also Corollary 3.37).

(ii) As observed by Hinich in the erratum of [15], it is always possible to replace an operad by a cofibrant one up to equivalence.

(iii) Frequently, the homotopy theory of algebras over cofibrant operads is well-behaved (see §4.1). In fact, the algebraic structure over a cofibrant operad is usually a homotopy invariant notion in the sense that there are *homotopy transfer theorems*: morally speaking, given an equivalence between an algebra and a second object, one can find an algebra structure on the second object lifting the equivalence to an equivalence of algebras (see [39] and [11]). This fact is particularly important for dg-operads. For instance, it implies that any algebraic structure on a complex governed by a cofibrant dg-operad can be transferred to its cohomology giving rise to an equivalent algebra. There is an extensive literature around this problem; in particular, concerned with explicit formulae for the transferred algebraic structure. See [6, 32, 20, 38].

So far, we have been working with morphisms of operads which are the identity on colors. When the set of colors is not fixed, there is also a natural candidate for equivalence in $\mathtt{Opd}(\mathtt{V})$:

**Definition 2.13.** A map $\varphi\colon \mathcal{O} \to \mathcal{P}$ in $\mathtt{Opd}(\mathtt{V})$ is a *Dwyer-Kan equivalence*, or *DK-equivalence*, if the following conditions hold:

- (homotopy essentially surjective) the induced functor $\pi_0(\overline{\mathcal{O}}) \to \pi_0(\overline{\mathcal{P}})$ is essentially surjective (see Remark 2.4), and

- (homotopy fully-faithful) the map $\mathcal{O}\left[\begin{smallmatrix} c \\ d \end{smallmatrix}\right] \xrightarrow{\sim} \mathcal{P}\left[\begin{smallmatrix} \varphi(c) \\ \varphi(d) \end{smallmatrix}\right]$ is an equivalence in $\mathcal{V}$ for any O-corolla.

Under some technical conditions on $\mathcal{V}$ and using the bifibration $\mathrm{col}\colon \mathtt{Opd}(\mathtt{V}) \to \mathtt{Set}$, Caviglia in [9] glued all the model structures $\{\mathcal{O}\mathrm{pd}_O(\mathcal{V})\}_O$ into a model category $\mathcal{O}\mathrm{pd}(\mathcal{V})$ which presents the homotopy theory of operads with DK-equivalences. This homotopy theory is quite important for theoretical purposes and one of the first questions one can ask is when does the assignment $\mathcal{O} \mapsto \mathcal{A}\mathrm{lg}_\mathcal{O}(\mathcal{V})$ send DK-equivalences to Quillen equivalences, generalizing the usual *rectification problem* of operadic algebras (see §4.1).

# 3 Enveloping operads

## 3.1 Decorating trees and enveloping operads

In the sequel, we will make extensive use of rooted trees and isomorphisms between them, e.g. in order to define enveloping operads and to perform their homotopical analysis.

Recall that there are many different definitions of (rooted) trees in the literature, most of them equivalent. In this work, we will consider rooted trees in the sense of [22, §1.5], or equivalently [12, §A]. Thus, a rooted tree $\tau$ is given by a finite set of vertices $\mathsf{V}(\tau)$, a finite set of edges $\mathsf{E}(\tau)$, with adjacency relations to the vertices, and a distinguished external edge $r$ called the root of $\tau$. An isomorphism of rooted trees $\varphi\colon \tau \xrightarrow{\cong} \tau'$ consists of a pair of bijections $\varphi_{\mathsf{V}}\colon \mathsf{V}(\tau) \xrightarrow{\cong} \mathsf{V}(\tau')$, $\varphi_{\mathsf{E}}\colon \mathsf{E}(\tau) \xrightarrow{\cong} \mathsf{E}(\tau')$ which preserve the adjacency relations and the root.

If we are given a set $\mathsf{O}$, we can define $\mathsf{O}$-colored rooted trees, or rooted $\mathsf{O}$-trees, by simply considering rooted trees $\tau$ equipped with a labeling function $\mathsf{f}(\tau)\colon \mathsf{E}(\tau) \to \mathsf{O}$. With this modification, an isomorphism of rooted $\mathsf{O}$-trees becomes an isomorphism of rooted trees $\varphi\colon \tau \to \tau'$ such that the following triangle commutes

$$\begin{array}{ccc} \mathsf{E}(\tau) & \xrightarrow{\varphi_{\mathsf{E}}} & \mathsf{E}(\tau') \\ & {\scriptstyle \mathsf{f}(\tau)}\searrow \quad \swarrow{\scriptstyle \mathsf{f}(\tau')} & \\ & \mathsf{O} & \end{array} .$$

In other words, it is an isomorphism of rooted trees that preserves the $\mathsf{O}$-labels.

Given a vertex $\nu$ of a rooted $\mathsf{O}$-tree, we will denote by $\left[\begin{smallmatrix} s(\nu) \\ t(\nu) \end{smallmatrix}\right]$ the associated $\mathsf{O}$-corolla; i.e. the finite set of adjacent edges to $\nu$ with their $\mathsf{O}$-labels and the canonical partition into inputs and output.

**Decorating trees.** For any $\Sigma_{\mathsf{O}}$-sequence $\mathcal{X}$, we can consider a functor decorating trees $(\mathcal{X})\colon (\mathsf{Tree}^{\cong}_{\mathsf{O}})^{\mathrm{op}} \to \mathsf{V}$ by setting $(\mathcal{X})(\tau) = \bigotimes_{\nu \in \mathsf{V}(\tau)} \mathcal{X}\left[\begin{smallmatrix} s(\nu) \\ t(\nu) \end{smallmatrix}\right]$, where $\mathsf{Tree}^{\cong}_{\mathsf{O}}$ is the groupoid of rooted $\mathsf{O}$-trees. This functor is easier to understand with a picture:

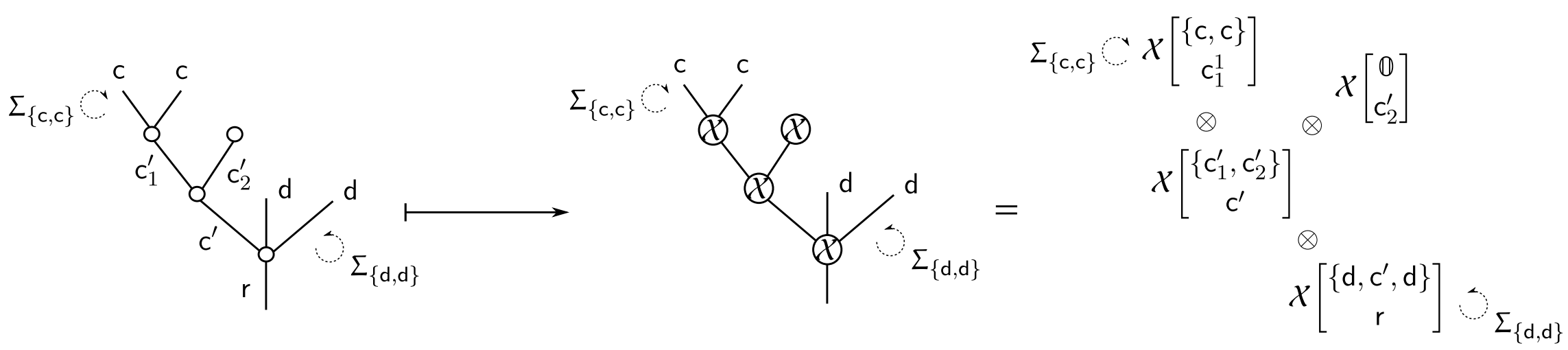

We will need more refined versions of this decoration process in order to deal with enveloping operads. Thus, we introduce two additional groupoids of trees and tagging functors for them.

**Definition 3.1.** The groupoid $\mathsf{Tree}^{\cong}_{\mathsf{O},\lightning}$ of *rooted* $\mathsf{O}$*-trees with two kinds of leaves* is defined by:

**Obj:** Its objects are rooted O-trees $\tau$ whose set of leaves $L(\tau)$ decomposes as a disjoint union, snaky and straight leaves, $L(\tau) = L_{\lightning}(\tau) \amalg L_{\downarrow}(\tau)$. Equivalently, $\tau$ comes equipped with a function $L(\tau) \to \{\lightning, \downarrow\}$.

**Mor:** Its morphisms are isomorphisms of rooted O-trees that preserve types of leaves, or equivalently, isomorphisms whose restriction to leaves commute with the functions into $\{\lightning, \downarrow\}$.

Trees in $\mathsf{Tree}^{\simeq}_{O,\lightning}$ should be decorated by a pair of objects: a $\Sigma_O$-sequence for vertices and a $\Sigma_O$-sequence concentrated in arity 0 for straight leaves (which should be seen as corks).

**Definition 3.2.** Let $\mathcal{X}, \mathcal{C} \in \Sigma\mathtt{Seq}_O(\mathtt{V})$ be symmetric sequences, where $\mathcal{C}$ is concentrated in arity 0. Define the tagging functor $(\mathcal{X};\mathcal{C})\colon (\mathsf{Tree}^{\simeq}_{O,\lightning})^{\mathrm{op}} \to \mathtt{V}$ by:

- on objects, it sends a tree $\tau \in \mathsf{Tree}^{\simeq}_{O,\lightning}$ to

$$(\mathcal{X};\mathcal{C})(\tau) = \bigotimes_{v\in V(\tau)} \mathcal{X}\begin{bmatrix} s(v) \\ t(v) \end{bmatrix} \otimes \bigotimes_{\ell \in L_{\downarrow}(\tau)} \mathcal{C}\begin{bmatrix} \mathbb{0} \\ \ell \end{bmatrix}.$$

- on morphisms, it sends an isomorphism of trees to the corresponding symmetric action on the tensor factors; e.g.

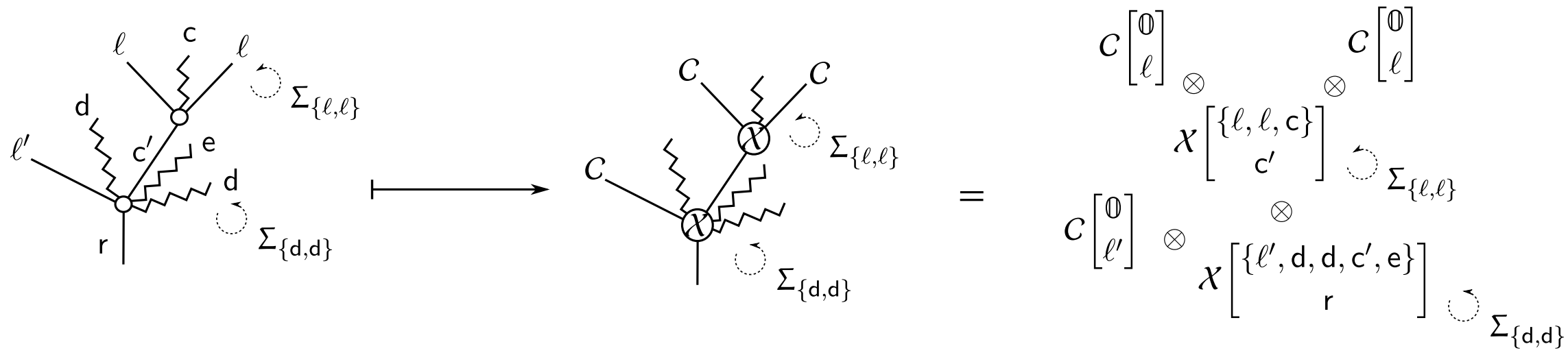

**Definition 3.3.** The groupoid $\mathsf{Tree}^{\mathrm{lv},\simeq}_{O,\lightning}$ of *leveled rooted* O*-trees with two kinds of leaves* is defined by:

**Obj:** Its objects are rooted O-trees $\tau$ with two kinds of leaves equipped with a level structure, that is, their set of vertices decompose as a finite disjoint union according to a height number $V(\tau) = V_0(\tau) \amalg \cdots \amalg V_{h(\tau)}(\tau)$ and the components are related by a chain of functions

$$V_{h(\tau)}(\tau) \to V_{h(\tau)-1}(\tau) \to \cdots \to V_1(\tau) \to V_0(\tau)$$

that assign to every vertex the unique vertex sitting below. The *height* of $\tau$ is the number $h(\tau)$.

**Mor:** Its morphisms are isomorphisms of rooted O-trees that preserve types of leaves and levels.

*Remark* 3.4. One can also associate heights to edges of a tree $\tau \in \mathsf{Tree}^{\mathrm{lv},\simeq}_{O,\lightning}$. The root is always considered to have height 0, internal edges have height equal to the height of their top vertex and leaves have height equal to the height of the unique vertex having it as source plus 1; e.g.

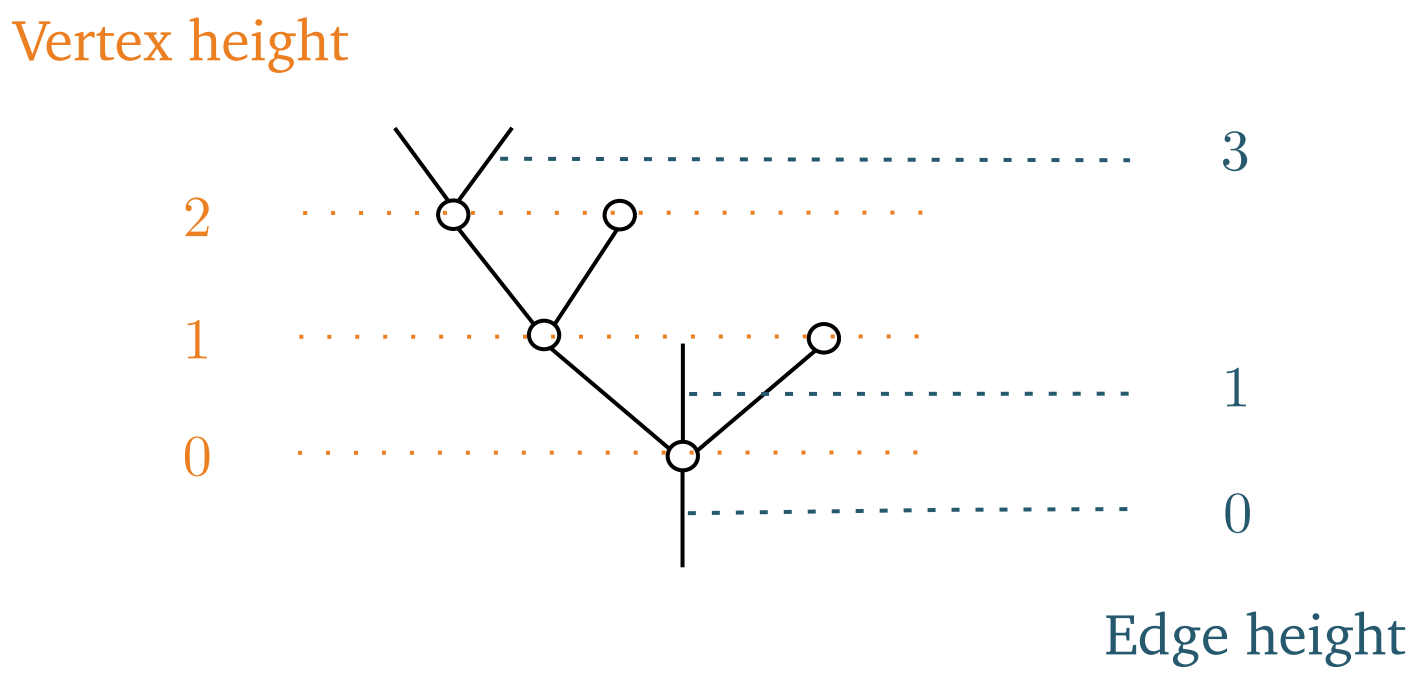

Note that the former notion of leveled tree does not coincide with that of May tree [21, §5]. May-trees have all their leaves in the highest height. For instance, if no nullary operations are considered, the chain of functions for May-trees

$$V_{h(\tau)}(\tau) \to V_{h(\tau)-1}(\tau) \to \cdots \to V_1(\tau) \to V_0(\tau)$$

is comprised of surjections.

Now, we introduce the most general tagging functor that we will need in this chapter.

**Definition 3.5.** Let $\mathcal{E}, \mathcal{O}, \mathcal{C} \in \Sigma\mathrm{Seq}_{\mathrm{O}}(\mathtt{V})$ be symmetric sequences, where $\mathcal{C}$ is concentrated in arity 0. Define the functor $(\mathcal{E}; \mathcal{O}; \mathcal{C})\colon (\mathsf{Tree}^{\mathrm{lv},\simeq}_{\mathrm{O},\notin})^{\mathrm{op}} \to \mathtt{V}$ by:

- on objects, it sends a tree $\tau \in \mathsf{Tree}^{\mathrm{lv},\simeq}_{\mathrm{O},\notin}$ to the object

$$(\mathcal{E}; \mathcal{O}; \mathcal{C})(\tau) = \bigotimes_{v \in V_{\mathrm{even}}(\tau)} \mathcal{E}\begin{bmatrix} s(v) \\ t(v) \end{bmatrix} \otimes \bigotimes_{u \in V_{\mathrm{odd}}(\tau)} \mathcal{O}\begin{bmatrix} s(u) \\ t(u) \end{bmatrix} \otimes \bigotimes_{\ell \in L_{\downarrow}(\tau)} \mathcal{C}\begin{bmatrix} \mathbb{0} \\ \ell \end{bmatrix},$$

  where $V_{\mathrm{even}}(\tau) = \coprod_{i \text{ even}} V_i(\tau)$ and $V_{\mathrm{odd}}(\tau) = \coprod_{i \text{ odd}} V_i(\tau)$.

- on morphisms, it sends an isomorphism of trees to the corresponding symmetric action on the tensor factors (see Definition 3.2).

It is also admissible to replace one symmetric collection in $(\mathcal{X})$, $(\mathcal{X}; \mathcal{C})$ or $(\mathcal{E}; \mathcal{O}; \mathcal{C})$ by a morphism of symmetric collections, resulting in a tagging functor with values in the arrow category $\mathtt{V}^2$. Let us specify one of these variations.

**Definition 3.6.** Let $f\colon \mathcal{O} \to \mathcal{O}'$ be a morphism of symmetric collections and $\mathcal{E}, \mathcal{C}$ two additional symmetric collections, with $\mathcal{C}$ concentrated in arity 0. The functor $(\mathcal{E}; f; \mathcal{C})\colon (\mathsf{Tree}^{\mathrm{lv},\simeq}_{\mathrm{O},\notin})^{\mathrm{op}} \to \mathtt{V}^2$ is defined by:

- on objects, it sends a tree $\tau \in \mathsf{Tree}^{\mathrm{lv},\simeq}_{\mathrm{O},\notin}$ to the morphism in $\mathtt{V}$

$$(\mathcal{E}; f; \mathcal{C})(\tau) = \bigotimes_{v \in V_{\mathrm{even}}(\tau)} \mathcal{E}\begin{bmatrix} s(v) \\ t(v) \end{bmatrix} \otimes \mathop{\Box}_{u \in V_{\mathrm{odd}}(\tau)} f\begin{bmatrix} s(u) \\ t(u) \end{bmatrix} \otimes \bigotimes_{\ell \in L_{\downarrow}(\tau)} \mathcal{C}\begin{bmatrix} \mathbb{0} \\ \ell \end{bmatrix}.$$

- on morphisms, its sends an isomorphism of trees to the natural transformation between morphisms in $\mathtt{V}$ whose components are the evaluations of $(\mathcal{E}; \mathcal{O}; \mathcal{C})$ and $(\mathcal{E}; \mathcal{O}'; \mathcal{C})$ on that isomorphism.

*Remark* 3.7. Checking that the functor $(\mathcal{E}; \mathsf{f}; \mathcal{C})$ is well defined amounts to prove certain compatibilities between permutations and f. They hold since f is a morphism of symmetric collections.

Last but not least, it is important to note that assigning to a tree its set of (snaky) leaves and its root defines a functor $\left[\begin{smallmatrix} s \\ t \end{smallmatrix}\right] \colon \mathsf{Tree}^{(\mathsf{lv}),\simeq}_{\mathsf{O},(\lightning)} \to \mathsf{Fin}^{\simeq}_{\mathsf{O}} \times \mathsf{O}$, where the target can be seen as the category of O-corollas. For a corolla $\left[\begin{smallmatrix} c \\ d \end{smallmatrix}\right] \in \mathsf{Fin}^{\simeq}_{\mathsf{O}} \times \mathsf{O}$, we denote by $\mathsf{Tree}^{(\mathsf{lv}),\simeq}_{\mathsf{O},(\lightning)}\left[\begin{smallmatrix} c \\ d \end{smallmatrix}\right]$ the slice category $\left[\begin{smallmatrix} s \\ t \end{smallmatrix}\right] \downarrow \left[\begin{smallmatrix} c \\ d \end{smallmatrix}\right]$, which consists on (decorated) rooted O-trees whose (snaky) leaves and root have selected labels $\left[\begin{smallmatrix} c \\ d \end{smallmatrix}\right]$. The canonical morphism $\mathsf{Tree}^{(\mathsf{lv}),\simeq}_{\mathsf{O},(\lightning)}\left[\begin{smallmatrix} c \\ d \end{smallmatrix}\right] \to \mathsf{Tree}^{(\mathsf{lv}),\simeq}_{\mathsf{O},(\lightning)}$ allows one to restrict the tagging functors defined above.

**Definition of enveloping operad.** Given an operad $\mathcal{O}$ and an $\mathcal{O}$-algebra A, one can define an operad $\mathcal{O}_{\mathsf{A}}$ whose category of algebras is isomorphic to the category of $\mathcal{O}$-algebras under A, i.e. $\mathtt{Alg}_{\mathcal{O}_\mathsf{A}} \simeq \mathsf{A} \downarrow \mathtt{Alg}_{\mathcal{O}}$. It is called *enveloping operad* and its construction can be informally described as modifying $\mathcal{O}$ to have A in arity 0. Our first goal is to recall the precise definition of this gadget using trees. See [10, 29, 30, 40].

To define the underlying symmetric sequence of $\mathcal{O}_{\mathsf{A}}$ (and later its operad structure), we consider the reflexive coequalizer in $\Sigma\mathtt{Seq}_{\mathsf{O}}(\mathtt{V}) \simeq \Sigma\mathtt{Coll}_{\mathsf{O}}(\mathtt{V})$

$$\mathcal{O}^1_{\mathsf{A}} \underset{d_1}{\overset{d_0}{\rightrightarrows}} \mathcal{O}^0_{\mathsf{A}} \xrightarrow{\ \mathsf{colim}\ } \mathcal{O}_{\mathsf{A}}, \qquad s_0 \colon \mathcal{O}^0_{\mathsf{A}} \dashrightarrow \mathcal{O}^1_{\mathsf{A}}$$

where:

- $\mathcal{O}^0_{\mathsf{A}}$ evaluated on $\left[\begin{smallmatrix} c \\ d \end{smallmatrix}\right]$ is given by the colimit of the tagging functor $(\mathcal{O}; \mathcal{O}; \mathsf{A})$ over the full subgroupoid of $\mathsf{Tree}^{\mathsf{lv},\simeq}_{\mathsf{O},\lightning}\left[\begin{smallmatrix} c \\ d \end{smallmatrix}\right]$ spanned by corollas (trees with height $\leqslant 0$)

$$\mathcal{O}^0_{\mathsf{A}}\begin{bmatrix} c \\ d \end{bmatrix} = \underset{\tau}{\mathsf{colim}}\ (\mathcal{O}; \mathcal{O}; \mathsf{A})(\tau) \equiv \underset{\tau}{\mathsf{colim}}\ (\mathcal{O}; \mathsf{A})(\tau).$$

  $\Sigma_{\mathsf{O}}$ acts on the snaky leaves (and hence on factors of type $\mathcal{O}$); e.g.

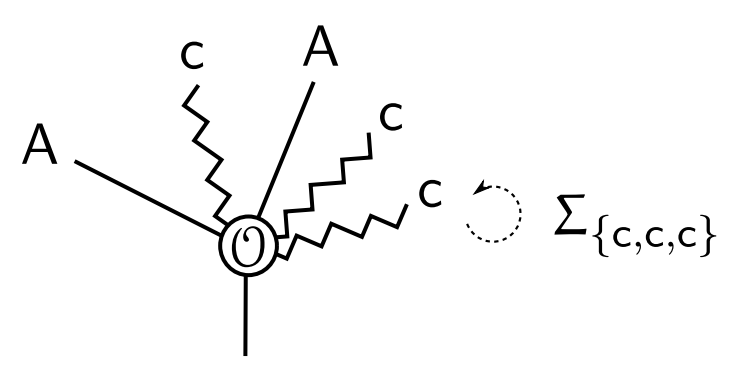

- $\mathcal{O}^1_{\mathsf{A}}$ evaluated on $\left[\begin{smallmatrix} c \\ d \end{smallmatrix}\right]$ is given by the colimit of the tagging functor $(\mathcal{O}; \mathcal{O}; \mathsf{A})$ over the full subgroupoid of $\mathsf{Tree}^{\mathsf{lv},\simeq}_{\mathsf{O},\lightning}\left[\begin{smallmatrix} c \\ d \end{smallmatrix}\right]$ spanned by trees with height $\leqslant 1$ whose level 1 leaves are snaky and whose level 2 leaves (if any) are straight

$$\mathcal{O}^1_{\mathsf{A}}\begin{bmatrix} c \\ d \end{bmatrix} = \underset{\Psi}{\mathsf{colim}}\ (\mathcal{O}; \mathcal{O}; \mathsf{A})(\Psi).$$

  $\Sigma_{\mathsf{O}}$ acts on the indexing category of trees of height $\leqslant 1$ (and hence on factors of type $\mathcal{O}$).

- $s_0\colon \mathcal{O}^0_A \to \mathcal{O}^1_A$ is induced by subdivision of straight leaves, e.g.

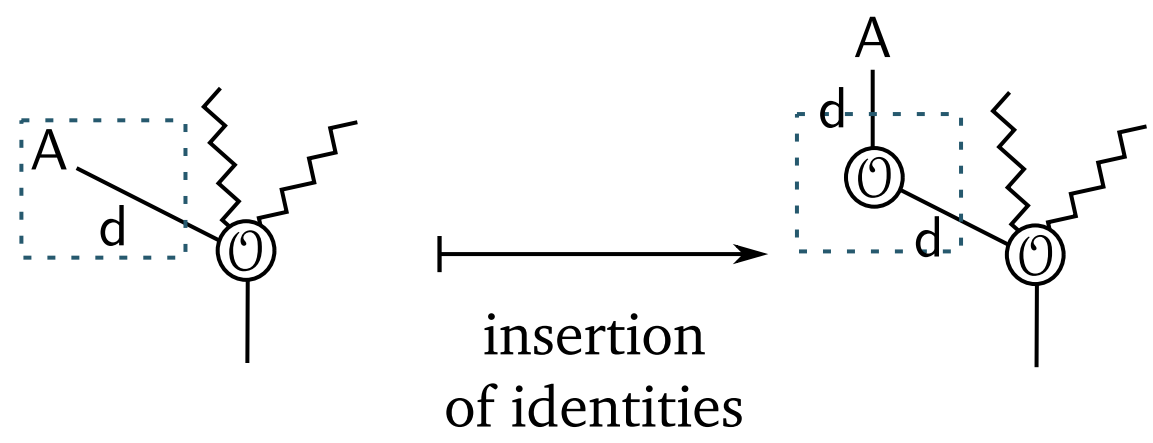

- $d_0\colon \mathcal{O}^1_A \to \mathcal{O}^0_A$ is induced by contracting corollas whose vertex lives in level 1 using the $\mathcal{O}$-algebra structure on A (they may have straight leaves or no leaves), e.g.

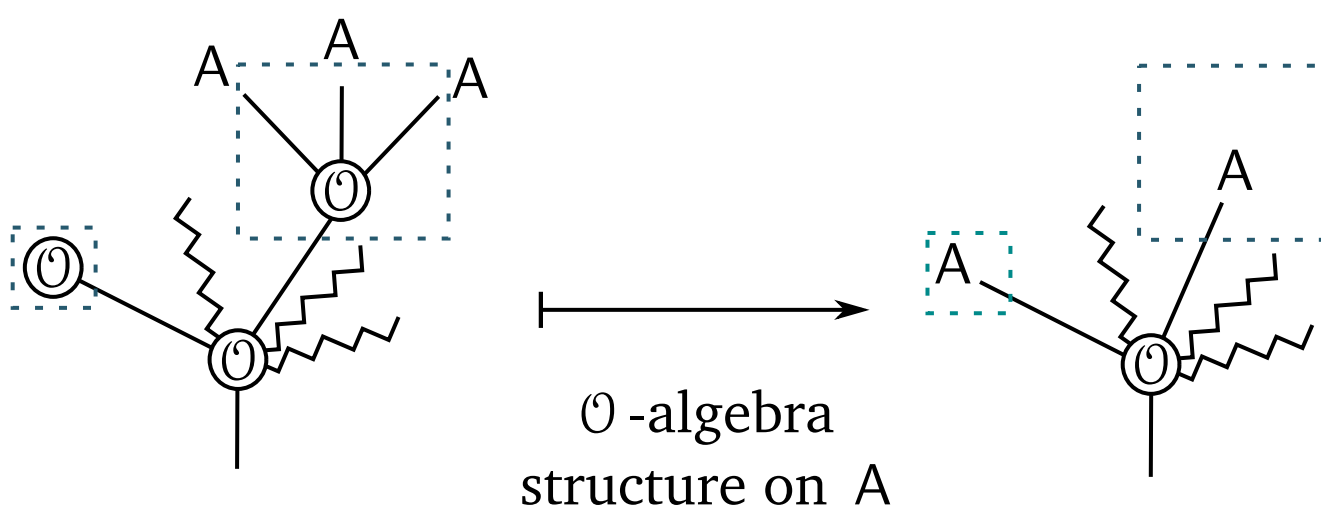

- $d_1\colon \mathcal{O}^1_A \to \mathcal{O}^0_A$ is induced by contracting internal edges in level 1 using the operad structure of $\mathcal{O}$, e.g.

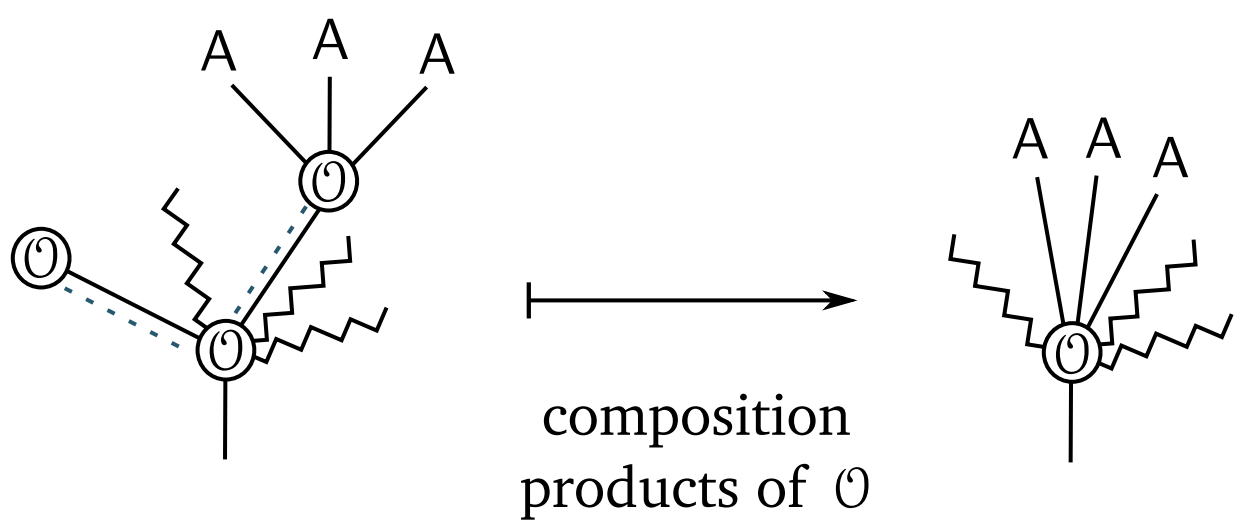

**Definition 3.8.** Let $\mathcal{O}$ be an operad and A an $\mathcal{O}$-algebra. The *enveloping operad* $\mathcal{O}_A$ associated to the pair $(\mathcal{O}, A)$ is determined by:

- its underlying collection is the tip of the colimit diagram defined above,

$$\mathcal{O}^1_A \underset{d_1}{\overset{d_0}{\underset{\xleftarrow{\;s_0\;}}{\rightrightarrows}}} \mathcal{O}^0_A \xrightarrow{\text{colim}} \mathcal{O}_A \;;$$

- its identity at c is $\mathrm{id}_c\colon \mathbb{I} \to \mathcal{O}^0_A \left[\begin{smallmatrix} c \\ c \end{smallmatrix}\right] \to \mathcal{O}_A \left[\begin{smallmatrix} c \\ c \end{smallmatrix}\right]$ where the map on the left is induced by $\mathrm{id}_c\colon \mathbb{I} \to \mathcal{O} \left[\begin{smallmatrix} c \\ c \end{smallmatrix}\right]$ (it picks $\mathrm{id}_c$ in the corolla with one snaky leaf);
- its composition product $\circ_j\colon \mathcal{O}_A \left[\begin{smallmatrix} c \\ d \end{smallmatrix}\right] \otimes \mathcal{O}_A \left[\begin{smallmatrix} a \\ c_j \end{smallmatrix}\right] \to \mathcal{O}_A \left[\begin{smallmatrix} c \circ_j a \\ d \end{smallmatrix}\right]$ is induced by grafting trees and contracting the newly created inner edge via the operad structure of $\mathcal{O}$.

**Proposition 3.9.** *There is an equivalence of categories* $\mathtt{Alg}_{\mathcal{O}_\mathsf{A}}(\mathtt{V}) \simeq \mathsf{A}{\downarrow}\mathtt{Alg}_{\mathcal{O}}(\mathtt{V})$.

*Remark* 3.10. There is an alternative way to characterize the enveloping operad; see the universal property written down in [10, §4.1.1].

*Remark* 3.11. The enveloping operad $\mathcal{O}_\mathsf{A}$ comes with a canonical operad morphism $\mathcal{O} \to \mathcal{O}_\mathsf{A}$ that induces the usual adjunction $\mathtt{Alg}_{\mathcal{O}}(\mathtt{V}) \rightleftarrows \mathsf{A} \downarrow \mathtt{Alg}_{\mathcal{O}}(\mathtt{V})$. Using the description given above, $\mathcal{O}\left[\begin{smallmatrix}\mathsf{c}\\ \mathsf{d}\end{smallmatrix}\right] \to \mathcal{O}_\mathsf{A}\left[\begin{smallmatrix}\mathsf{c}\\ \mathsf{d}\end{smallmatrix}\right]$ is induced by decorating corollas with no straight leaves.

There is an alternative presentation of enveloping operads by means of the $\circ$-product of symmetric sequences given in [40]. Recall that the $\circ$-product of $\mathcal{X}, \mathcal{Y} \in \Sigma\mathtt{Seq}_{\mathsf{O}}(\mathtt{V})$ can be defined by

$$(\mathcal{X} \circ \mathcal{Y})\begin{bmatrix}\mathsf{c}\\ \mathsf{d}\end{bmatrix} = \operatorname*{colim}_{\Psi}\,(\mathcal{X};\mathcal{Y};\mathbb{I})(\Psi),$$

where we take colimit of the tagging functor evaluated on the full subgroupoid of $\mathsf{Tree}^{\mathsf{lv},\simeq}_{\mathsf{O},\not{\downarrow}}\left[\begin{smallmatrix}\mathsf{c}\\ \mathsf{d}\end{smallmatrix}\right]$ spanned by trees of height $= 1$ whose leaves (if any) are snaky and sit in level 2 (i.e. we can forget about having two kinds of leaves).

We end up with:

**Proposition 3.12.** *The underlying symmetric sequence of the enveloping operad* $\mathcal{O}_\mathsf{A}$ *sits into a reflexive coequalizer diagram*

$$\operatorname*{colim}_{\tau}\,(\mathcal{O};\mathcal{O}\circ\mathsf{A})(\tau) \;\underset{d_1}{\overset{d_0}{\rightrightarrows}}\!\!\overset{s_0}{\longleftarrow}\; \operatorname*{colim}_{\tau}\,(\mathcal{O};\mathsf{A})(\tau) \longrightarrow \mathcal{O}_\mathsf{A}\left[\begin{smallmatrix}\mathsf{c}\\ \mathsf{d}\end{smallmatrix}\right],$$

*where* $(\mathrm{i})$ *both colimits are indexed by the full subgroupoid of* $(\mathsf{Tree}^{\simeq}_{\mathsf{O},\not{\downarrow}}\left[\begin{smallmatrix}\mathsf{c}\\ \mathsf{d}\end{smallmatrix}\right])^{\mathsf{op}}$ *spanned by corollas,* $(\mathrm{ii})$ $d_0$ *applies the* $\mathcal{O}$*-algebra structure of* $\mathsf{A}$*,* $(\mathrm{iii})$ *operadic composition maps of* $\mathcal{O}$ *determine* $d_1$ *and* $(\mathrm{iv})$ *insertion of identities in* $\mathcal{O}$ *induces* $s_0$*.*

*Proof.* It follows from a direct comparison between this reflexive coequalizer and the one chosen as a definition for $\mathcal{O}_\mathsf{A}$. □

It is clear that the assignment $(\mathcal{O},\mathsf{A}) \mapsto \mathcal{O}_\mathsf{A}$ extends to a functor

$$\mathtt{OpdAlg}(\mathtt{V}) = \left\{\begin{matrix}\text{operadic algebras as pairs } (\mathcal{O},\mathsf{A})\\ \text{with } \mathcal{O}\in\mathtt{Opd}(\mathtt{V}) \text{ and } \mathsf{A}\in\mathtt{Alg}_{\mathcal{O}}(\mathtt{V})\end{matrix}\right\} \longrightarrow \mathtt{Opd}(\mathtt{V}),$$

which is moreover left adjoint to $\mathcal{O} \mapsto (\mathcal{O}, \mathcal{O}\left[\begin{smallmatrix}\mathbb{0}\\ \star\end{smallmatrix}\right])$, the functor that sends an operad $\mathcal{O}$ to the initial $\mathcal{O}$-algebra. Note that $\mathtt{OpdAlg}(\mathtt{V})$ is just the Grothendieck construction for $\mathcal{O} \mapsto \mathtt{Alg}_{\mathcal{O}}(\mathtt{V})$.

## 3.2 Cell attachments and filtrations

We are interested in studying how the functor $(\mathcal{O},\mathsf{A}) \mapsto \mathcal{O}_\mathsf{A}$ behaves with respect to cellular cofibrations because of the possible homotopical applications. The hardest part to understand such behavior corresponds to the analysis of pushouts in both variables. We will describe filtrations that build such pushouts in a controlled way.

**Filtrations associated to pushouts of free maps.** Let us start with the not so-well studied behavior of enveloping operads with respect to its operad variable, i.e. varying $\mathcal{O}$.

**Proposition 3.13** (Free on operad)**.** *Let $\mathcal{F}(\mathcal{X})$ be the free operad on a symmetric sequence $\mathcal{X} \in \Sigma\mathsf{Seq}_{\mathsf{O}}(\mathtt{V})$ and $\mathsf{A}$ be an $\mathcal{F}(\mathcal{X})$-algebra. Then, there is a natural isomorphism*

$$\big(\mathcal{F}(\mathcal{X})\big)_{\mathsf{A}}\begin{bmatrix}\mathsf{c}\\ \mathsf{d}\end{bmatrix} \cong \underset{\Upsilon}{\mathsf{colim}}\,(\mathcal{X};\mathsf{A})(\Upsilon),$$

*where the colimit is taken over the full subgroupoid of $(\mathsf{Tree}^{\cong}_{\mathsf{O},\natural}\left[\begin{smallmatrix}\mathsf{c}\\ \mathsf{d}\end{smallmatrix}\right])^{\mathsf{op}}$ spanned by trees with no corollas of straight leaves.*

*Proof.* We will find the natural isomorphism in the statement by analyzing the description of $\big(\mathcal{F}(\mathcal{X})\big)_{\mathsf{A}}\left[\begin{smallmatrix}\mathsf{c}\\ \mathsf{d}\end{smallmatrix}\right]$ as a reflexive coequalizer, i.e.

$$\big(\mathcal{F}(\mathcal{X})\big)^{1}_{\mathsf{A}}\left[\begin{smallmatrix}\mathsf{c}\\ \mathsf{d}\end{smallmatrix}\right] \rightrightarrows \big(\mathcal{F}(\mathcal{X})\big)^{0}_{\mathsf{A}}\left[\begin{smallmatrix}\mathsf{c}\\ \mathsf{d}\end{smallmatrix}\right] \xrightarrow{\mathsf{colim}} \big(\mathcal{F}(\mathcal{X})\big)_{\mathsf{A}}\left[\begin{smallmatrix}\mathsf{c}\\ \mathsf{d}\end{smallmatrix}\right] \,.$$

The idea is to apply the description in terms of trees of the free operad $\mathcal{F}(\mathcal{X})$ in this coequalizer. Being more precise, $\mathcal{F}(\mathcal{X})\left[\begin{smallmatrix}\mathsf{a}\\ \mathsf{b}\end{smallmatrix}\right] \cong \mathsf{colim}_{\Psi}(\mathcal{X})(\Psi)$, where the colimit is indexed by the groupoid $(\mathsf{Tree}^{\cong}_{\mathsf{O}}\left[\begin{smallmatrix}\mathsf{a}\\ \mathsf{b}\end{smallmatrix}\right])^{\mathsf{op}}$. Taking this into account, we find an isomorphism of reflexive coequalizers

$$\begin{array}{ccc}
\big(\mathcal{F}(\mathcal{X})\big)^{1}_{\mathsf{A}}\left[\begin{smallmatrix}\mathsf{c}\\ \mathsf{d}\end{smallmatrix}\right] & \xrightarrow{\cong} & \underset{\Psi'}{\mathsf{colim}}\,(\mathcal{X};\mathsf{A})(\Psi') \\
\downarrow\uparrow\downarrow & & \downarrow\uparrow\downarrow \\
\big(\mathcal{F}(\mathcal{X})\big)^{0}_{\mathsf{A}}\left[\begin{smallmatrix}\mathsf{c}\\ \mathsf{d}\end{smallmatrix}\right] & \xrightarrow{\cong} & \underset{\Psi}{\mathsf{colim}}\,(\mathcal{X};\mathsf{A})(\Psi)
\end{array},$$

where $\Psi$ runs over $\mathsf{Tree}^{\cong}_{\mathsf{O},\natural}$ and $\Psi'$ runs over the groupoid of decomposable trees (with isomorphisms preserving decompositions) $\Psi' = \mathsf{grafting}(\Upsilon;\gamma_1,\dots,\gamma_t)$ where $(\mathrm{i})$ $\Upsilon,\gamma_j \in \mathsf{Tree}^{\cong}_{\mathsf{O},\natural}$ for all $\mathrm{j}$, $(\mathrm{ii})$ $\Upsilon$ have exactly $\mathrm{t}$ straight leaves, where $\gamma_j$'s are inserted, although $\mathrm{t}$ may vary for different trees and $(\mathrm{iii})$ all the leaves (if any) of all $\gamma_j$'s are straight. The vertical maps on the right can be readily guessed from the definition of those on the left, e.g. one of the downpointing arrows takes an indexing tree $\Psi' = \mathsf{grafting}(\Upsilon;\gamma_1,\dots,\gamma_t)$ and contracts $\gamma_j$'s using the $\mathcal{F}(\mathcal{X})$-algebra structure on $\mathsf{A}$.

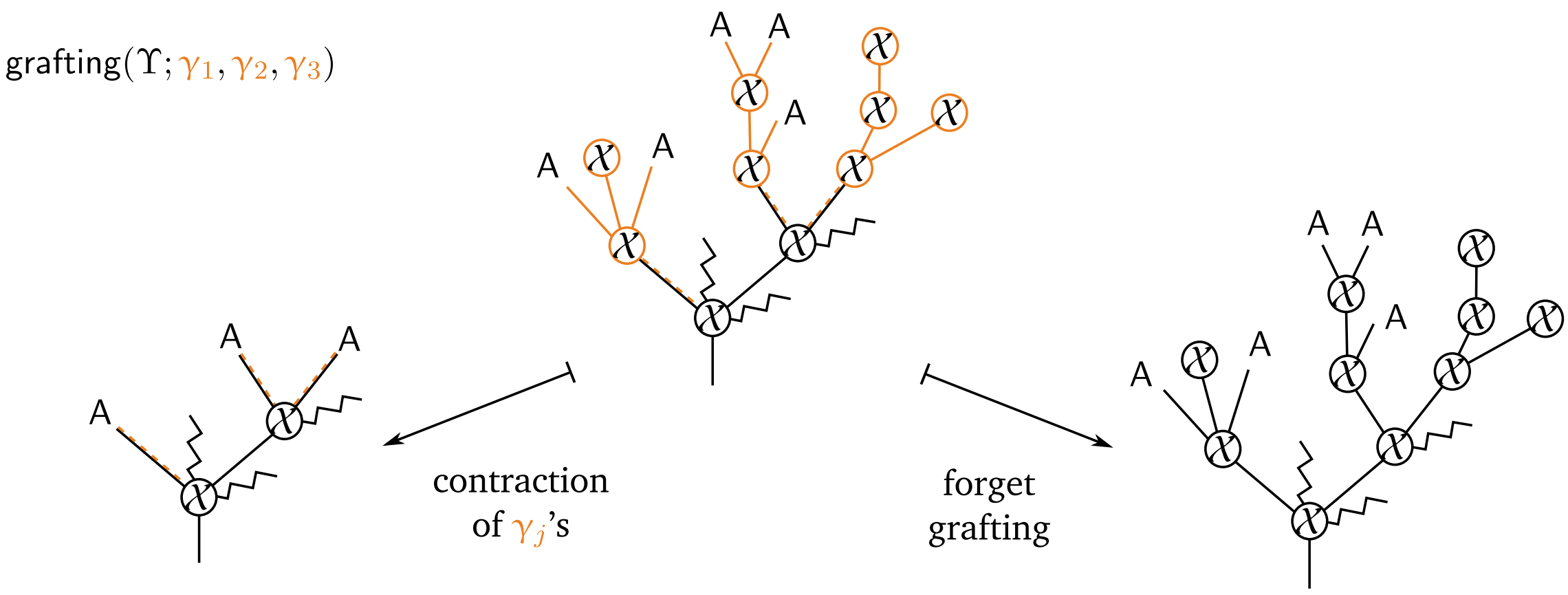

The colimit of the vertical coequalizer on the right can be found by extending it to a split coequalizer

$$\operatorname*{colim}_{\Psi'}(\mathcal{X};\mathsf{A})(\Psi') \rightrightarrows \operatorname*{colim}_{\Psi}(\mathcal{X};\mathsf{A})(\Psi) \longrightarrow \operatorname*{colim}_{\Upsilon}(\mathcal{X};\mathsf{A})(\Upsilon),$$

where $\Upsilon$ runs over the full subgroupoid of $(\mathsf{Tree}^{\cong}_{\mathcal{O},\natural}\left[\begin{smallmatrix}\mathsf{c}\\ \mathsf{d}\end{smallmatrix}\right])^{\mathsf{op}}$ spanned by trees with no corollas of straight leaves, i.e. the groupoid given in the statement. The left curved dashed map

$$\operatorname*{colim}_{\Psi'}(\mathcal{X};\mathsf{A})(\Psi') \dashleftarrow \operatorname*{colim}_{\Psi}(\mathcal{X};\mathsf{A})(\Psi)$$

decomposes $\Psi$ as grafting$(\Upsilon;\gamma_1,\dots,\gamma_t)$, where $\gamma_j$'s are the maximal subtrees of $\Psi$ which have only straight leaves. The right curved dashed map

$$\operatorname*{colim}_{\Psi}(\mathcal{X};\mathsf{A})(\Psi) \dashleftarrow \operatorname*{colim}_{\Upsilon}(\mathcal{X};\mathsf{A})(\Upsilon)$$

is just a colimit factor inclusion. Lastly, the remaining arrow

$$\operatorname*{colim}_{\Psi}(\mathcal{X};\mathsf{A})(\Psi) \longrightarrow \operatorname*{colim}_{\Upsilon}(\mathcal{X};\mathsf{A})(\Upsilon)$$

collapses any corolla of straight leaves, again via the $\mathcal{F}(\mathcal{X})$-algebra structure on $\mathsf{A}$.

This shows that there is a natural isomorphism $\big(\mathcal{F}(\mathcal{X})\big)_{\mathsf{A}}\left[\begin{smallmatrix}\mathsf{c}\\ \mathsf{d}\end{smallmatrix}\right] \cong \operatorname{colim}_{\Upsilon}(\mathcal{X};\mathsf{A})(\Upsilon)$ since both objects share an universal property.

□

**Proposition 3.14** (Pushout on operads)**.** *Let*

$$\begin{array}{ccc} \mathcal{F}(\mathcal{X}) & \longrightarrow & \mathcal{O} \\ {\scriptstyle \mathcal{F}\mathsf{j}}\downarrow & \ulcorner & \downarrow{\scriptstyle \mathsf{g}} \\ \mathcal{F}(\mathcal{Y}) & \longrightarrow & \mathcal{O}[\mathsf{j}] \end{array}$$

*be a pushout square of operads, where* $\mathsf{j}\colon \mathcal{X} \to \mathcal{Y}$ *is a map of symmetric collections, and let* $\mathsf{A}$ *be an* $\mathcal{O}[\mathsf{j}]$*-algebra. Then, the map* $\mathsf{g}_\diamond\colon \mathcal{O}_{\mathsf{A}}\left[\begin{smallmatrix}\mathsf{c}\\ \mathsf{d}\end{smallmatrix}\right] \to \mathcal{O}[\mathsf{j}]_{\mathsf{A}}\left[\begin{smallmatrix}\mathsf{c}\\ \mathsf{d}\end{smallmatrix}\right]$ *is the transfinite composite of a sequence*

$$\mathcal{O}_{\mathsf{A}}\begin{bmatrix}\mathsf{c}\\ \mathsf{d}\end{bmatrix} = \mathcal{O}[\mathsf{j}]_{\mathsf{A},0}\begin{bmatrix}\mathsf{c}\\ \mathsf{d}\end{bmatrix} \longrightarrow \cdots \longrightarrow \mathcal{O}[\mathsf{j}]_{\mathsf{A},\mathsf{t}-1}\begin{bmatrix}\mathsf{c}\\ \mathsf{d}\end{bmatrix} \xrightarrow{\mathsf{g}_\mathsf{t}} \mathcal{O}[\mathsf{j}]_{\mathsf{A},\mathsf{t}}\begin{bmatrix}\mathsf{c}\\ \mathsf{d}\end{bmatrix} \longrightarrow \cdots,$$

*where each factor* $\mathsf{g}_\mathsf{t}$ *fits into a pushout in* $\mathcal{V}$

$$\begin{array}{ccc} \bullet & \longrightarrow & \mathcal{O}[\mathsf{j}]_{\mathsf{A},\mathsf{t}-1}\left[\begin{smallmatrix}\mathsf{c}\\ \mathsf{d}\end{smallmatrix}\right] \\ {\scriptstyle \operatorname*{colim}_{\Lambda}(\mathcal{O}_{\mathsf{A}};\mathsf{j};\mathsf{A})(\Lambda)}\downarrow & \ulcorner & \downarrow{\scriptstyle \mathsf{g}_\mathsf{t}} \\ \bullet & \longrightarrow & \mathcal{O}[\mathsf{j}]_{\mathsf{A},\mathsf{t}}\left[\begin{smallmatrix}\mathsf{c}\\ \mathsf{d}\end{smallmatrix}\right] \end{array},$$

*where the colimit in* $\mathtt{V}^2$ *is taken over the full subcategory of* $(\mathsf{Tree}^{\mathsf{lv},\simeq}_{\mathsf{O},\not\downarrow}\left[\begin{smallmatrix}\mathsf{c}\\ \mathsf{d}\end{smallmatrix}\right])^{\mathsf{op}}$ *spanned by trees whose external vertices are all even and whose leaves are all snaky, there are* $\mathtt{t}$ *odd vertices (all of them without leaves), and such that they contain no forbidden subtrees of the form*

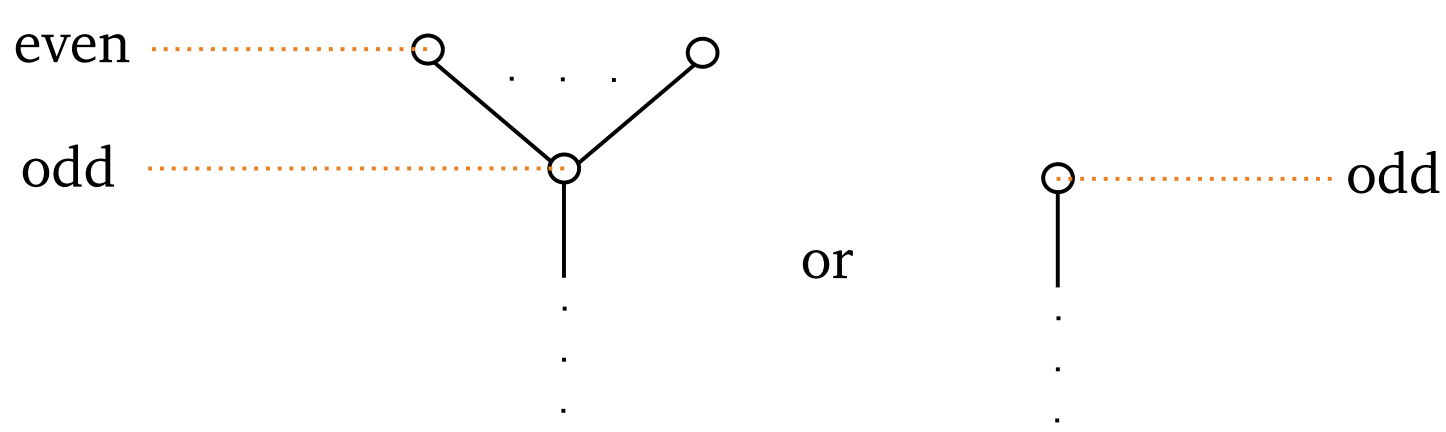

*Proof.* The strategy is a refinement of the one used in Proposition 3.13. Again, we will exploit the reflexive coequalizer in $\Sigma\mathtt{Seq}_{\mathsf{O}}(\mathtt{V})$

$$\mathcal{O}^1_{\mathsf{A}} \rightrightarrows \mathcal{O}^0_{\mathsf{A}} \xrightarrow{\mathsf{colim}} \mathcal{O}_{\mathsf{A}} \ .$$

In this case, instead of a colimit description of $\mathcal{O}[\mathtt{j}]$ as in Proposition 3.13, what we know about the underlying symmetric collection of $\mathcal{O}[\mathtt{j}]$ is that the map $\mathcal{O} \to \mathcal{O}[\mathtt{j}]$ can be seen as the transfinite composite of a sequence in $\Sigma\mathtt{Seq}_{\mathsf{O}}(\mathtt{V})$

$$\mathcal{O} = \mathcal{O}[\mathtt{j}]_0 \longrightarrow \cdots \longrightarrow \mathcal{O}[\mathtt{j}]_{r-1} \xrightarrow{g_r} \mathcal{O}[\mathtt{j}]_r \longrightarrow \cdots ,$$

where each $g_t$ is a cobase change of a map defined in terms of trees (see [3]). Since the abstract non-sense in the following proof is involved, we will split the proof in two parts: the first one deals with the categorical argument from which we deduce the claim and the second one fills the categorical argument with the specific input from the statement.

*First part: categorical argument.*

Using the decomposition of $\mathcal{O} \to \mathcal{O}[\mathtt{j}]$ as a transfinite colimit, one can also write the reflexive coequalizer defining $\mathcal{O}[\mathtt{j}]_{\mathsf{A}}$ as a transfinite colimit of reflexive coequalizers (in $\Sigma\mathtt{Seq}_{\mathsf{O}}(\mathtt{V})$); more concretely, there is a sequence of reflexive coequalizers

$$\begin{array}{ccccccccc}
\mathcal{R}^1_0 & \longrightarrow & \cdots & \longrightarrow & \mathcal{R}^1_{t-1} & \xrightarrow{g^1_t} & \mathcal{R}^1_t & \longrightarrow & \cdots \\
\downdownarrows & & & & \downdownarrows & & \downdownarrows & & \\
\mathcal{R}^0_0 & \longrightarrow & \cdots & \longrightarrow & \mathcal{R}^0_{t-1} & \xrightarrow[g^0_t]{} & \mathcal{R}^0_t & \longrightarrow & \cdots
\end{array}$$

whose transfinite colimit, which we denote by $\mathcal{R}^\bullet_\omega$, is the reflexive coequalizer defining $\mathcal{O}[\mathtt{j}]_{\mathsf{A}}$. Thus, by a commutation of colimits, $\mathcal{O}[\mathtt{j}]_{\mathsf{A}}$ is the transfinite colimit of the sequence of maps obtained above by taking vertical coequalizers, i.e. the sequence

$$\mathsf{coeq}\ \mathcal{R}^\bullet_0 \longrightarrow \cdots \longrightarrow \mathsf{coeq}\ \mathcal{R}^\bullet_{t-1} \xrightarrow{g_t} \mathsf{coeq}\ \mathcal{R}^\bullet_t \longrightarrow \cdots .$$

We have to show that this is the transfinite composite alluded in the statement, so let us use that notation for simplicity, $\mathsf{coeq}\ \mathcal{R}^\bullet_t = \mathcal{O}[\mathtt{j}]_{\mathsf{A},t}$. We should prove that the initial term in the sequence $\mathcal{O}[\mathtt{j}]_{\mathsf{A},0}$ coincides with $\mathcal{O}_{\mathsf{A}}$ and that each sucessive map $g_t$ fits into the pushout that we claimed in the statement (since colimits in diagram categories are computed objectwise).

Note that $\mathcal{O}[j] = \mathrm{colim}_t\ \mathcal{O}[j]_t$ and so

$$\mathcal{R}^0_\omega \begin{bmatrix} c \\ d \end{bmatrix} = \underset{\tau}{\mathrm{colim}} \left(\underset{t}{\mathrm{colim}}\, \mathcal{O}[j]_t; A\right)(\tau) \cong \underset{t}{\mathrm{colim}}\, \underset{\tau}{\mathrm{colim}} \left( \mathcal{O}[j]_t; A\right)(\tau)$$

by a commutation of colimits (since $\tau$ is just a corolla, the tagging functor only takes into account $\mathcal{O}[j]$ tagging the 0-level vertex) and for this reason we set

$$\mathcal{R}^0_t \begin{bmatrix} c \\ d \end{bmatrix} = \mathrm{colim}_\tau\, (\mathcal{O}[j]_t; A)(\tau).$$

Similarly, $\mathcal{R}^1_\omega \left[\begin{smallmatrix} c \\ d \end{smallmatrix}\right]$ is the evaluation of $\mathrm{colim}_\Psi\,(-;-;A)(\Psi)$ at $(\mathrm{colim}_p\, \mathcal{O}[j]_p; \mathrm{colim}_q\, \mathcal{O}[j]_q)$. Hence, $\mathcal{R}^1_t \left[\begin{smallmatrix} c \\ d \end{smallmatrix}\right]$ splits as a coproduct over isomorphism classes of trees $\Psi$ with height 1 and where the component associated to $[\Psi]$ (which has $s$ odd vertices) is

$$\underset{p+q_1+\cdots+q_s \leqslant t}{\mathrm{colim}}\ \underset{\mathrm{Aut}(\Psi)}{\mathrm{colim}}\, (\mathcal{O}[j]_p; (\mathcal{O}[j]_{q_1}|\cdots|\mathcal{O}[j]_{q_s}); A)(\Psi).$$

Here, $(\mathcal{O}[j]_{q_1}|\cdots|\mathcal{O}[j]_{q_s})$ means that we decorate each odd vertex with one $\mathcal{O}[j]_{q_j}$; e.g.

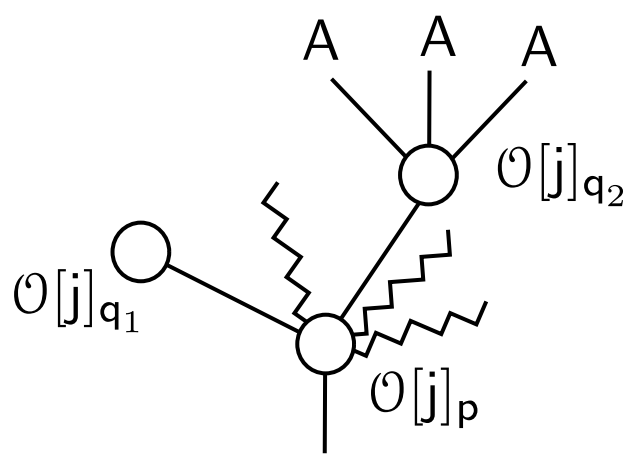

It is clear that $\mathcal{R}^\bullet_0 \left[\begin{smallmatrix} c \\ d \end{smallmatrix}\right]$ coincides with the reflexive coequalizer defining $\mathcal{O}_A \left[\begin{smallmatrix} c \\ d \end{smallmatrix}\right]$.

Let us take for granted, using Lemma A.7 (we will explain this point in the second part of the proof), that $\mathcal{R}^\bullet_{t-1} \to \mathcal{R}^\bullet_t$ fits into a pushout of reflexive coequalizer diagrams

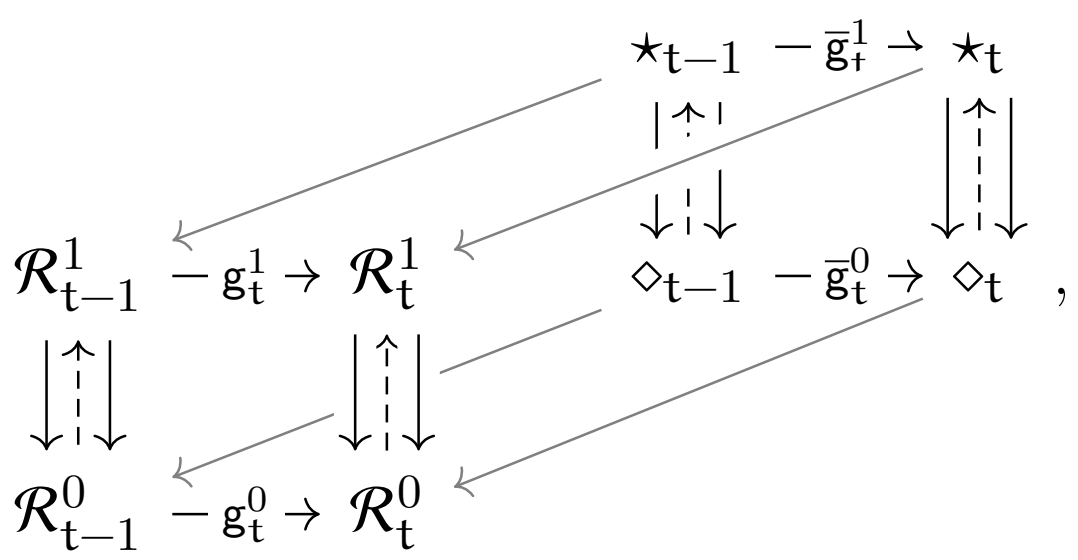

(3.1)

i.e. horizontal faces of the cube are pushouts. Note that taking vertical coequalizers in the previous cube one obtains a diagram whose horizontal squares are pushouts. Therefore, identifying $\overline{g}_t\colon \bullet_{t-1} \to \bullet_t$ with the map in the statement, one concludes the proof. Observe that we have not specified what $\overline{g}^0_t$ and $\overline{g}^1_t$ are (we will do that in the second part of the proof). Once we do so, we will be able to identify $\overline{g}_t = \mathrm{coeq}(\overline{g}^1_t \rightrightarrows \overline{g}^0_t)$ by computing the coequalizer in the back face of the cube.

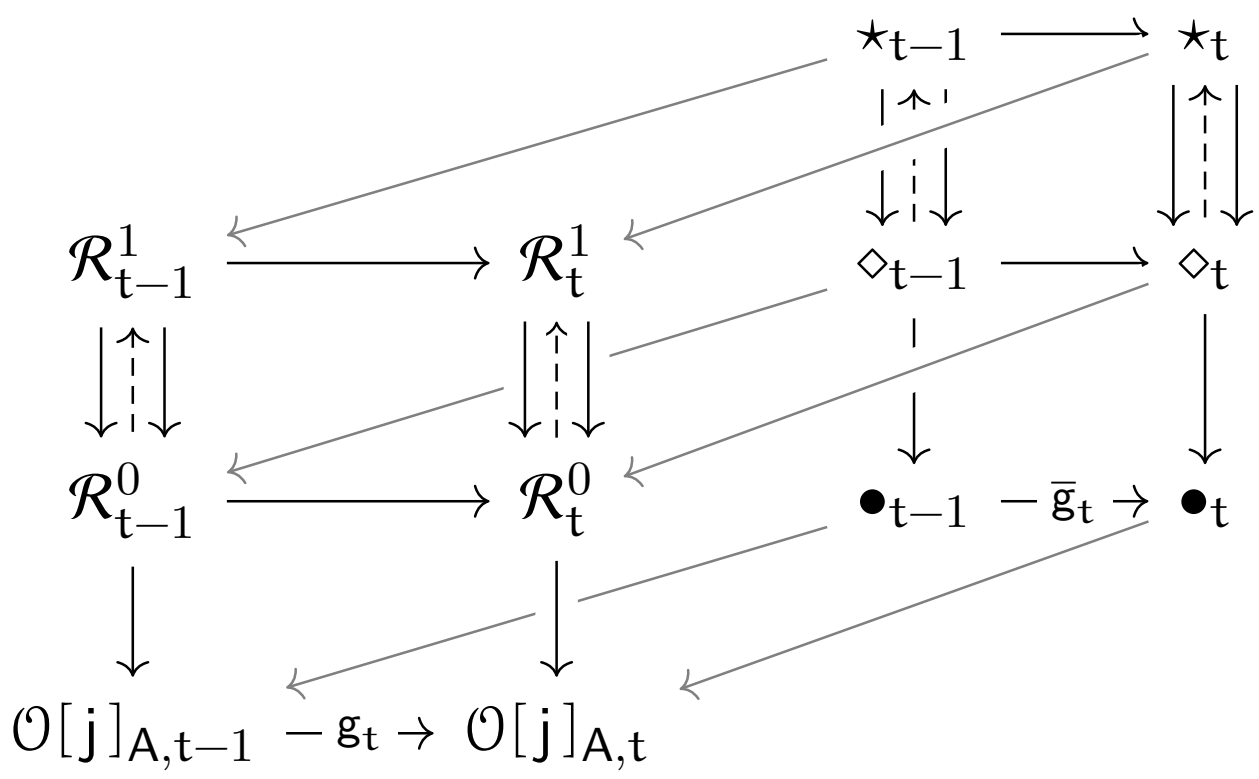

*Second part: specific descriptions.*

We have already defined $\mathcal{R}^1_t$ and $\mathcal{R}^0_t$. Also, the reflexive coequalizer $\mathcal{R}^1_t \rightrightarrows \mathcal{R}^0_t$ can be readily defined adapting the reflexive coequalizer describing enveloping operads; for example, one may visualize the parallel maps as

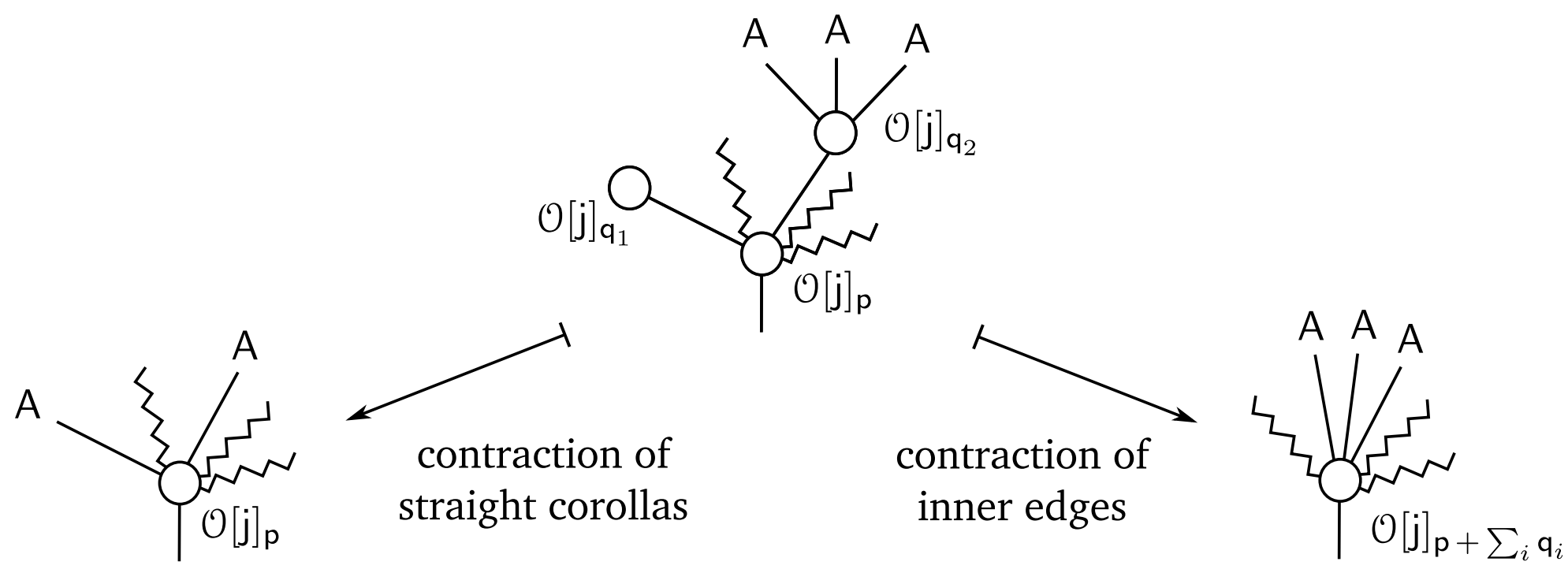

From the filtration of $\mathcal{O} \to \mathcal{O}[j]$, we deduce that $g^0_t\colon \mathcal{R}^0_{t-1}\left[\begin{smallmatrix}c\\d\end{smallmatrix}\right] \to \mathcal{R}^0_t\left[\begin{smallmatrix}c\\d\end{smallmatrix}\right]$ is the cobase change of $\varphi^0_t = \mathrm{colim}_\zeta\,(\mathcal{O};j;A)(\zeta)$ where $\zeta$ runs over trees in $\mathsf{Tree}^{\mathsf{lv},\simeq}_{\mathsf{O},\notin}\left[\begin{smallmatrix}c\\d\end{smallmatrix}\right]$ with t-odd inner vertices and no leaves in even level (see Lemma A.7); e.g.

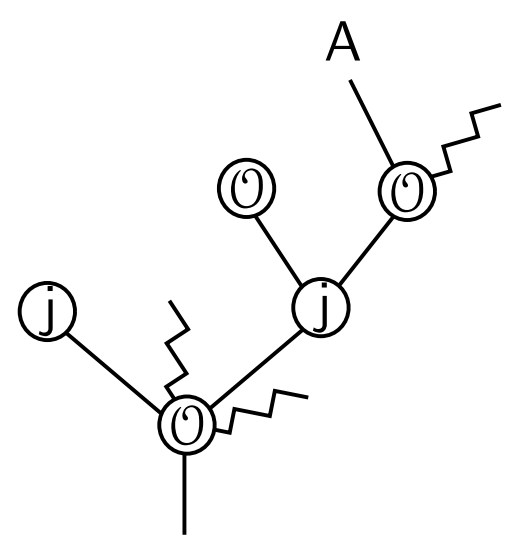

Analogously, one deduces that $\mathcal{R}^1_{t-1}\left[\begin{smallmatrix}c\\d\end{smallmatrix}\right] \to \mathcal{R}^1_t\left[\begin{smallmatrix}c\\d\end{smallmatrix}\right]$ is the cobase change of $\varphi^1_t$, the arrow obtained as the colimit of

$$\widetilde{\Psi} = \mathsf{grafting}(\Upsilon; \gamma_1, \dots, \gamma_m) \longmapsto (\mathcal{O}; j; A)(\Upsilon) \otimes \bigotimes_{1\leqslant r\leqslant m} (\mathcal{O}; j; A)(\gamma_r).$$

This diagram is indexed by trees $\widetilde{\Psi}$ constructed by grafting $\gamma_1, \dots, \gamma_m$ over $\Upsilon$, where: $(1)$ $\Upsilon \in \mathsf{Tree}^{\mathsf{lv},\simeq}_{\mathsf{O},\not{\iota}}\left[\begin{smallmatrix}c\\d\end{smallmatrix}\right]$ has exactly $m$ straight leaves, $(2)$ $\gamma_1, \dots, \gamma_m \in \mathsf{Tree}^{\mathsf{lv},\simeq}_{\mathsf{O},\not{\iota}}$ with none of them having snaky leaves, $(3)$ each $\gamma_r$ is grafted in one of the straight leaves of $\Upsilon$ and $(4)$ $\Upsilon, \gamma_1, \dots, \gamma_m$ have jointly $t$-odd inner vertices; e.g.

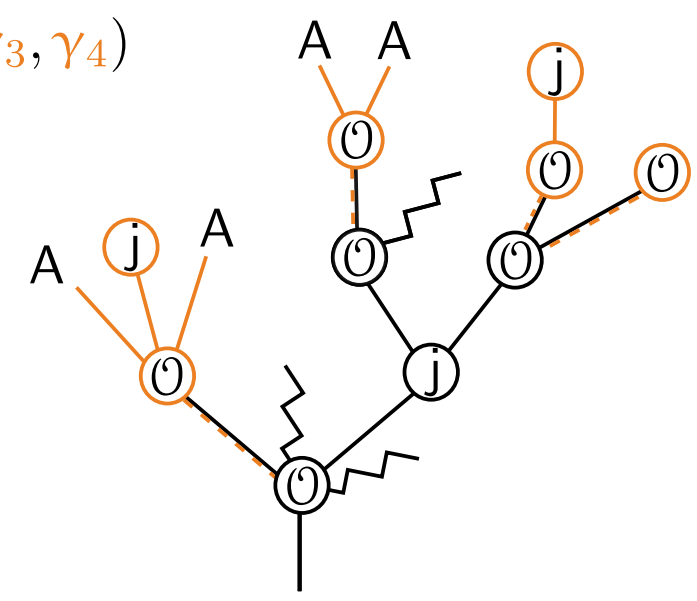

Note that the newly created inner edges in $\widetilde{\Psi}$ by grafting have both surrounding vertices decorated by $\mathcal{O}$. Since we will refer to those edges later, let us call them *grafting edges*.

It could seem that $\varphi^0_t$ and $\varphi^1_t$ fill the cube (3.1) when evaluating at $\left[\begin{smallmatrix}c\\d\end{smallmatrix}\right]$, although this is not the case: one of the vertical maps in the back face of the cube contracts the trees $\gamma_r$'s using the $\mathcal{O}$-algebra structure on $A$ as in

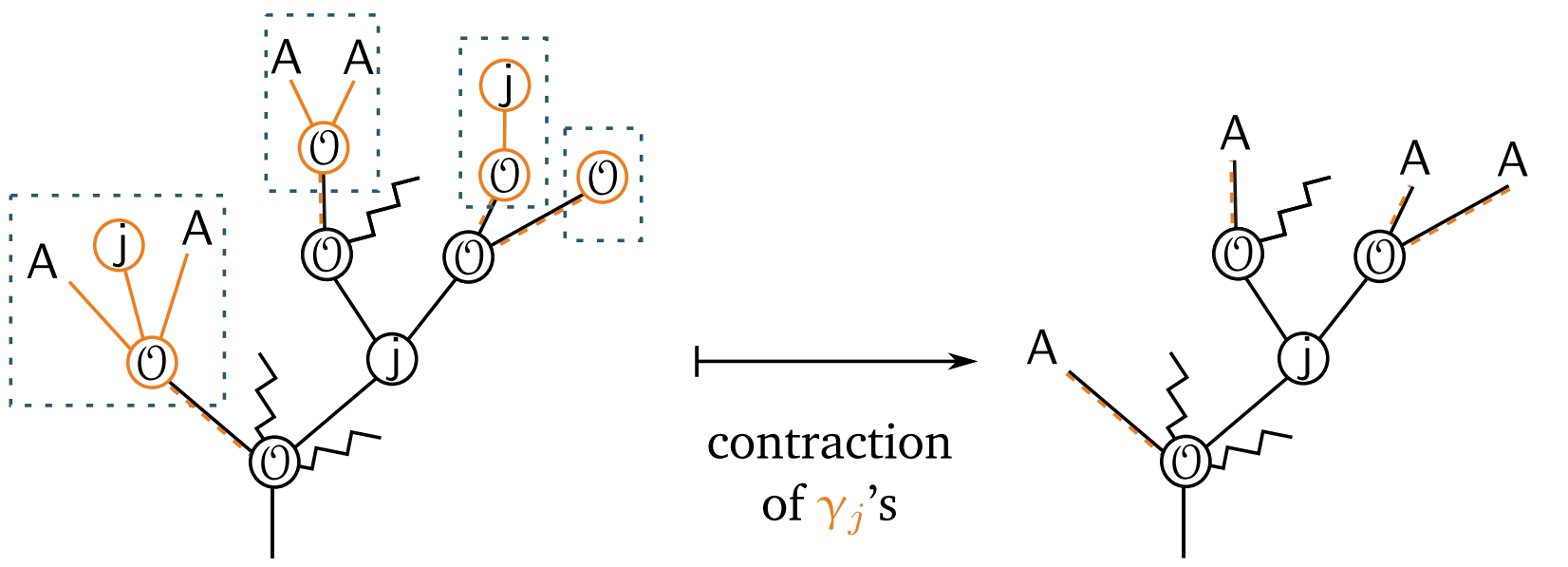

It is clear that the number of vertices decorated by j have decreased, so the map does not land in $\varphi^0_t$ but in a previous step in the filtration.

Hence, we add some wiggle room on the cells without modifying $g^k_t$ in order to obtain the cube (3.1). Simply take $\overline{g}^k_t\left[\begin{smallmatrix}c\\d\end{smallmatrix}\right] = \varphi^k_t \amalg \mathsf{id}_{\mathcal{R}^k_{t-1}\left[\begin{smallmatrix}c\\d\end{smallmatrix}\right]}$ for $k \in \{0,1\}$ and the obvious attaching map (twisted cell complex; see Remark A.8). With such a modification, we had a reflexive coequalizer $\varphi^1_t \amalg \mathsf{id}_{\mathcal{R}^1_{t-1}\left[\begin{smallmatrix}c\\d\end{smallmatrix}\right]} \rightrightarrows \varphi^0_t \amalg \mathsf{id}_{\mathcal{R}^0_{t-1}\left[\begin{smallmatrix}c\\d\end{smallmatrix}\right]}$ in $\mathsf{V}^2$ where the parallel maps $d_0$, $d_1$ are given by: (i) contraction of grafting edges in $\widetilde{\Psi}$ induces $d_0$ and (ii) full contraction of $\gamma_r$'s in $\widetilde{\Psi}$ yields $d_1$. Note that a tree corresponding to $\varphi^1_t$ goes via $d_1$ to the first component of $\varphi^0_t \amalg \mathsf{id}_{\mathcal{R}^0_{t-1}\left[\begin{smallmatrix}c\\d\end{smallmatrix}\right]}$ iff all $\gamma_r$'s are corollas. As mentioned before, see *categorical argument*, it remains to compute

$$\overline{g}_t = \mathsf{coeq}\left(\varphi^1_t \amalg \mathsf{id}_{\mathcal{R}^1_{t-1}\left[\begin{smallmatrix}c\\d\end{smallmatrix}\right]} \rightrightarrows \varphi^0_t \amalg \mathsf{id}_{\mathcal{R}^0_{t-1}\left[\begin{smallmatrix}c\\d\end{smallmatrix}\right]}\right)$$

to conclude the proof. To perform such calculation, observe that the groupoids of trees defining $\varphi_t^0$ and $\varphi_t^1$ can be decomposed into a disjoint union yielding $\varphi_t^k = \varphi_t^{k,c} \amalg \varphi_t^{k,nc}$ (c stands for contributing and nc for non-contributing). More precisely,

- For $k = 0$, $\varphi_t^{0,c}$ is the colimit indexed by trees $\zeta$ containing no subtrees of the following form (called bad subtrees)

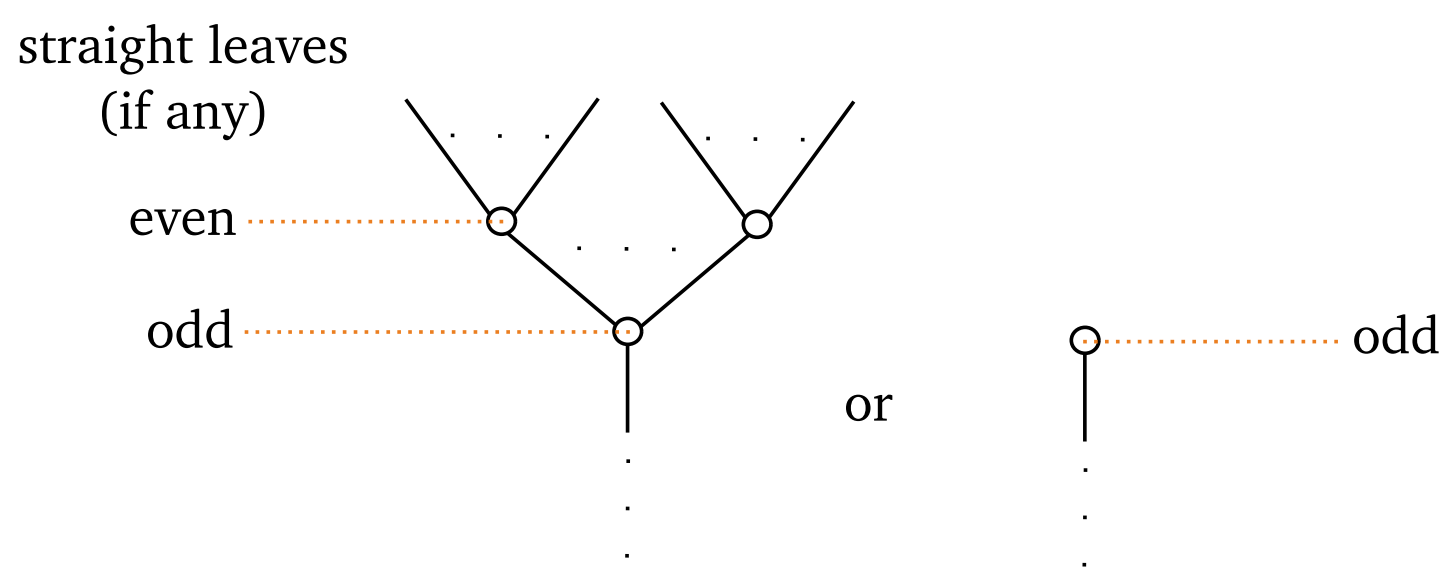

- For $k = 1$, $\varphi_t^{1,c}$ is the colimit indexed by trees $\widetilde{\Psi} = \text{grafting}(\Upsilon; \gamma_1, \ldots, \gamma_m)$ such that $\Upsilon$ contains no bad subtrees and all $\gamma_r$'s are corollas.

Clearly, $\varphi_t^1 \amalg \text{id}_{\mathcal{R}^1_{t-1}\left[\substack{c\\d}\right]} \rightrightarrows \varphi_t^0 \amalg \text{id}_{\mathcal{R}^0_{t-1}\left[\substack{c\\d}\right]}$ decomposes as the coproduct of

$$\varphi_t^{1,c} \rightrightarrows \varphi_t^{0,c} \qquad \text{and} \qquad \varphi_t^{1,nc} \amalg \text{id}_{\mathcal{R}^1_{t-1}\left[\substack{c\\d}\right]} \rightrightarrows \varphi_t^{0,nc} \amalg \text{id}_{\mathcal{R}^0_{t-1}\left[\substack{c\\d}\right]} .$$

By definition, coeq $\left(\varphi_t^{1,c} \rightrightarrows \varphi_t^{0,c}\right)$ is the map $\text{colim}_\Lambda\, (\mathcal{O}_A; j; A)(\Lambda)$ given in the statement and the cocone map just absorbs straight leaves through $\mathcal{O}_A^0 \to \mathcal{O}_A$, e.g.

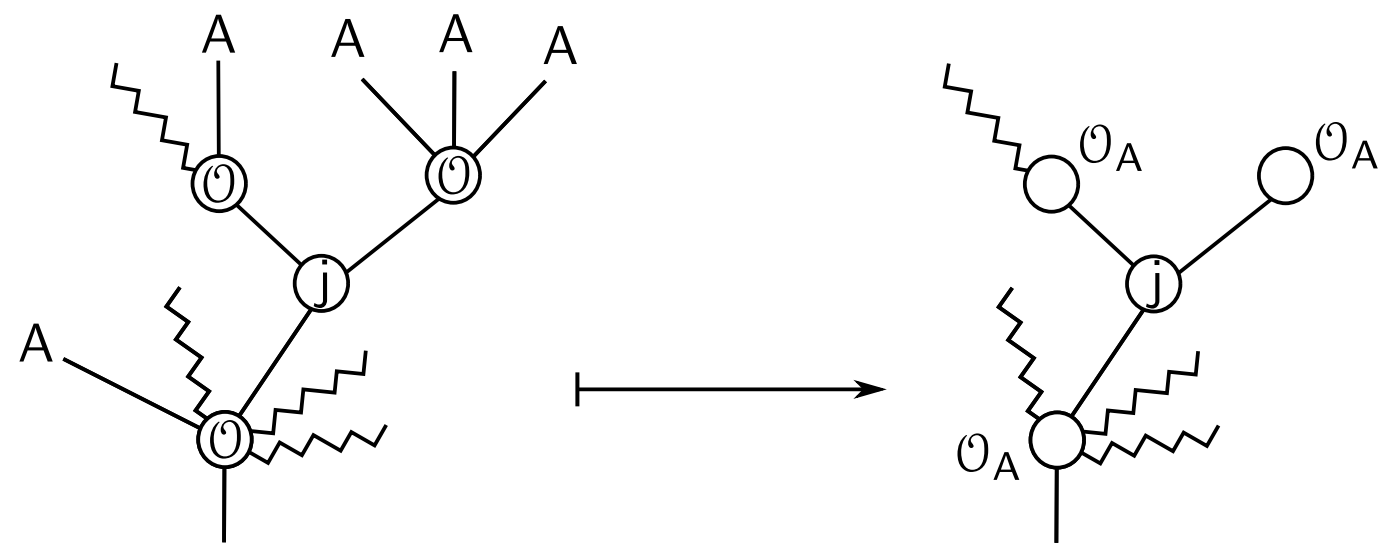

Note that the tree on the right hand side is obtained from the one on the left by removing the straight leaves. Thus, the fact that indexing trees in $\varphi_t^{0,c}$ do not contain bad subtrees yields that the indexing trees in the statement do not contain forbidden trees.

Lastly, we should check that

$$\text{id}_{\mathcal{O}[j]_{A,t-1}\left[\substack{c\\d}\right]} \cong \text{coeq}\left(\varphi_t^{1,nc} \amalg \text{id}_{\mathcal{R}^1_{t-1}\left[\substack{c\\d}\right]} \rightrightarrows \varphi_t^{0,nc} \amalg \text{id}_{\mathcal{R}^0_{t-1}\left[\substack{c\\d}\right]}\right)$$

to finish the computation of $\bar{g}_t$, since $g_t \colon \mathcal{O}[j]_{A,t-1}\left[\substack{c\\d}\right] \to \mathcal{O}[j]_{A,t}\left[\substack{c\\d}\right]$ is the cobase change of $\text{colim}_\Lambda\, (\mathcal{O}_A; j; A)(\Lambda)$ iff it is the pushout of $\text{colim}_\Lambda\, (\mathcal{O}_A; j; A)(\Lambda) \amalg \text{id}_{\mathcal{O}[j]_{A,t-1}\left[\substack{c\\d}\right]}$ with the obvious attaching map (twisted cell attachment).

Just note that we can apply Lemma 3.15 to the previous coequalizer and the following maps in $V^2$:

- $\xi\colon \varphi_t^{0,\mathrm{nc}} \to \varphi_t^{1,\mathrm{nc}}$ subdivides the bottom edge of a chosen bad subtree; e.g.

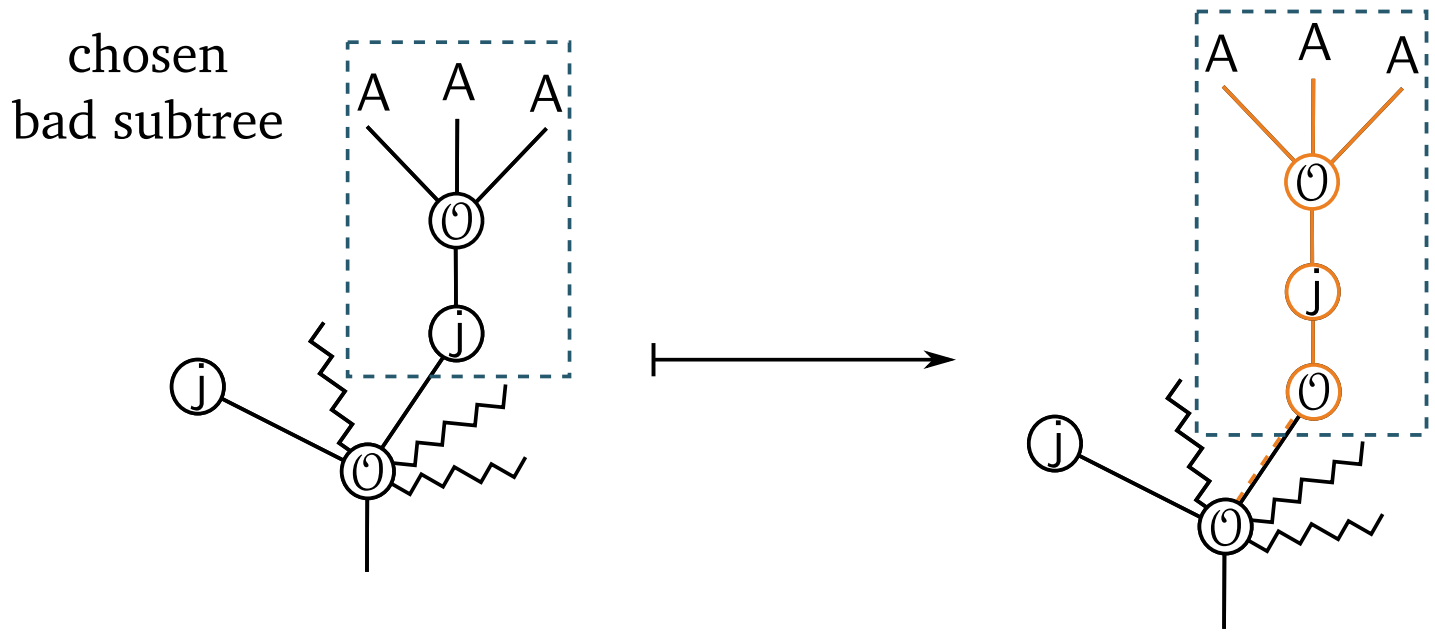

- $\epsilon\colon \varphi_t^{0,\mathrm{nc}} \to \mathrm{id}_{\mathcal{R}^0_{t-1}}\left[\begin{smallmatrix}c\\d\end{smallmatrix}\right]$ is the unique map filling the commutative diagram

$$\begin{array}{ccc} \varphi_t^{0,\mathrm{nc}} & \xrightarrow{\xi} & \varphi_t^{1,\mathrm{nc}} \\ {\scriptstyle\epsilon}\downarrow & & \searrow {\scriptstyle \text{contract } \gamma_r\text{'s}} \\ \mathrm{id}_{\mathcal{R}^0_{t-1}}\left[\begin{smallmatrix}c\\d\end{smallmatrix}\right] & \xrightarrow[\text{inclusion}]{} & \varphi_t^{0,\mathrm{nc}} \amalg \mathrm{id}_{\mathcal{R}^0_{t-1}}\left[\begin{smallmatrix}c\\d\end{smallmatrix}\right] \end{array} .$$

In other words, $\epsilon$ subdivides a chosen bad subtree, fully contracts the resulting subtree $\gamma$ and sends the result via the characteristic map to $\mathrm{id}_{\mathcal{R}^0_{t-1}}\left[\begin{smallmatrix}c\\d\end{smallmatrix}\right]$.

□

**Lemma 3.15.** *Suppose that it is given a coequalizer in a cocomplete category*

$$R^1 \underset{r_R}{\overset{l_R}{\rightrightarrows}} R^0 \xrightarrow{\pi} R$$

*and maps* $l_Q\colon Q^1 \to Q^0$, $r_Q\colon Q^1 \to Q^0 \amalg R^0$, $\epsilon\colon Q^0 \to R^0$. *Then, the fork*

$$Q^1 \amalg R^1 \underset{(r_Q, r_R)}{\overset{l_Q \amalg l_R}{\rightrightarrows}} Q^0 \amalg R^0 \xrightarrow{(\pi\epsilon, \pi)} R$$

*is also a coequalizer provided there is a section* $\xi\colon Q^0 \to Q^1$ *of* $l_Q$, *i.e.* $l_Q\xi = \mathrm{id}$, *making the diagram*

$$\begin{array}{ccccccc} Q^0 & \xrightarrow{\xi} & Q^1 & \xrightarrow{l_Q} & Q^0 & \xrightarrow{\epsilon} & R^0 \\ {\scriptstyle\epsilon}\downarrow & & & \searrow{\scriptstyle r_Q} & & & \downarrow{\scriptstyle\pi} \\ R^0 & & \xrightarrow[\text{inclusion}]{} & & Q^0 \amalg R^0 & \xrightarrow[(\pi\epsilon,\pi)]{} & R \end{array}$$

*fully commute.*

*Proof.* This is just an easy categorical computation. □

We now turn to the behavior of the enveloping operad with respect to its algebra variable, i.e. let us fix the operad $\mathcal{O}$ and study $A \mapsto \mathcal{O}_A$. This part has been studied in the literature before and so, we just collect the results for completeness.

**Proposition 3.16** (Free on algebra). *Let $X \in \mathtt{V}^{\mathsf{O}}$ be a symmetric sequence concentrated in arity 0. Then, there is a natural isomorphism*

$$\mathcal{O}_{\mathcal{O}\circ X}\begin{bmatrix} c \\ d \end{bmatrix} \cong \operatorname*{colim}_{\tau}\,(\mathcal{O};X)(\tau),$$

*where the colimit is taken over the full subgroupoid of $(\mathsf{Tree}^{\simeq}_{\mathsf{O},\natural}\left[\begin{smallmatrix} c \\ d \end{smallmatrix}\right])^{op}$ spanned by corollas.*

*Proof.* This is essentially a rephrasing of [40, Corollary 5.2.6]. □

**Proposition 3.17** (Pushout on algebras). *Consider a pushout square in $\mathtt{Alg}_{\mathcal{O}}(\mathtt{V})$*

$$\begin{array}{ccc} \mathcal{O}\circ X & \longrightarrow & A \\ {\scriptstyle \mathcal{O}\circ j}\downarrow & \ulcorner & \downarrow{\scriptstyle g} \\ \mathcal{O}\circ Y & \longrightarrow & A[j] \end{array} .$$

*Then, the induced map $g_*\colon \mathcal{O}_A\left[\begin{smallmatrix} c \\ d \end{smallmatrix}\right] \to \mathcal{O}_{A[j]}\left[\begin{smallmatrix} c \\ d \end{smallmatrix}\right]$ is the transfinite composite of a sequence*

$$\mathcal{O}_A\begin{bmatrix} c \\ d \end{bmatrix} = \mathcal{O}_{A[j],0}\begin{bmatrix} c \\ d \end{bmatrix} \longrightarrow \cdots \longrightarrow \mathcal{O}_{A[j],t-1}\begin{bmatrix} c \\ d \end{bmatrix} \xrightarrow{g_t} \mathcal{O}_{A[j],t}\begin{bmatrix} c \\ d \end{bmatrix} \longrightarrow \cdots,$$

*where each factor fits into a pushout square in $\mathtt{V}$*

$$\begin{array}{ccc} \bullet & \longrightarrow & \mathcal{O}_{A[j],t-1}\left[\begin{smallmatrix} c \\ d \end{smallmatrix}\right] \\ {\scriptstyle \operatorname*{colim}_{\tau}(\mathcal{O}_A;j)(\tau)}\downarrow & \ulcorner & \downarrow{\scriptstyle g_t} \\ \bullet & \longrightarrow & \mathcal{O}_{A[j],t}\left[\begin{smallmatrix} c \\ d \end{smallmatrix}\right] \end{array}$$

*such that the colimit in $\mathtt{V}^2$ is taken over the full subgroupoid of $(\mathsf{Tree}^{\simeq}_{\mathsf{O},\natural}\left[\begin{smallmatrix} c \\ d \end{smallmatrix}\right])^{op}$ spanned by corollas with $t$ straight leaves.*

*Proof.* Again, this is just a reformulation of [40, Proposition 5.3.2].[i]

Let us just comment that one may follow a similar strategy to that in Proposition 3.14, where we do apply the filtration of $A \to A[j]$ in contrast to [40, Remark 5.3.6]. We reproduce briefly the categorical argument that should be applied, which is a generalization of [29, Lemma 1.2].

Arguing with the filtration of $A \to A[j]$ in [40, Proposition 4.3.16] (generalized to a general $j$ in $\mathtt{V}^{\mathsf{O}}$), one obtains a cubical diagram whose vertical columns are coequalizers and horizontal squares (i.e. squares parallel to horizontal plane) are pushouts. Hence, we are reduced to compute $\overline{g}_t = \text{coeq}(\overline{g}_t^1 \rightrightarrows \overline{g}_t^0)$. In this case, by the specific input given in the statement, one constructs $\overline{g}_t^0$ and $\overline{g}_t^1$ via coequalizers, namely $\overline{g}_t^k = \text{coeq}(\overline{g}_t^{k,1} \rightrightarrows \overline{g}_t^{k,0})$ for $k \in \{0,1\}$. This fact yields a diagram whose diagonal components are coequalizers and its frontal face is a reflexive coequalizer of maps. A commutation of colimits implies

$$\begin{aligned}\overline{g}_t &= \text{coeq}\left(\text{coeq}(\overline{g}_t^{1\bullet}) \rightrightarrows \text{coeq}(\overline{g}_t^{0\bullet})\right)\\ &\cong \text{coeq}\left(\text{coeq}(\overline{g}_t^{\bullet 1}) \rightrightarrows \text{coeq}(\overline{g}_t^{\bullet 0})\right),\end{aligned}$$

that is, we can compute $\overline{g}_t$ by taking vertical coequalizers first and diagonal coequalizers later. The upshot is that this last nested coequalizer can be computed because the vertical coequalizers can be extended to split coequalizers.

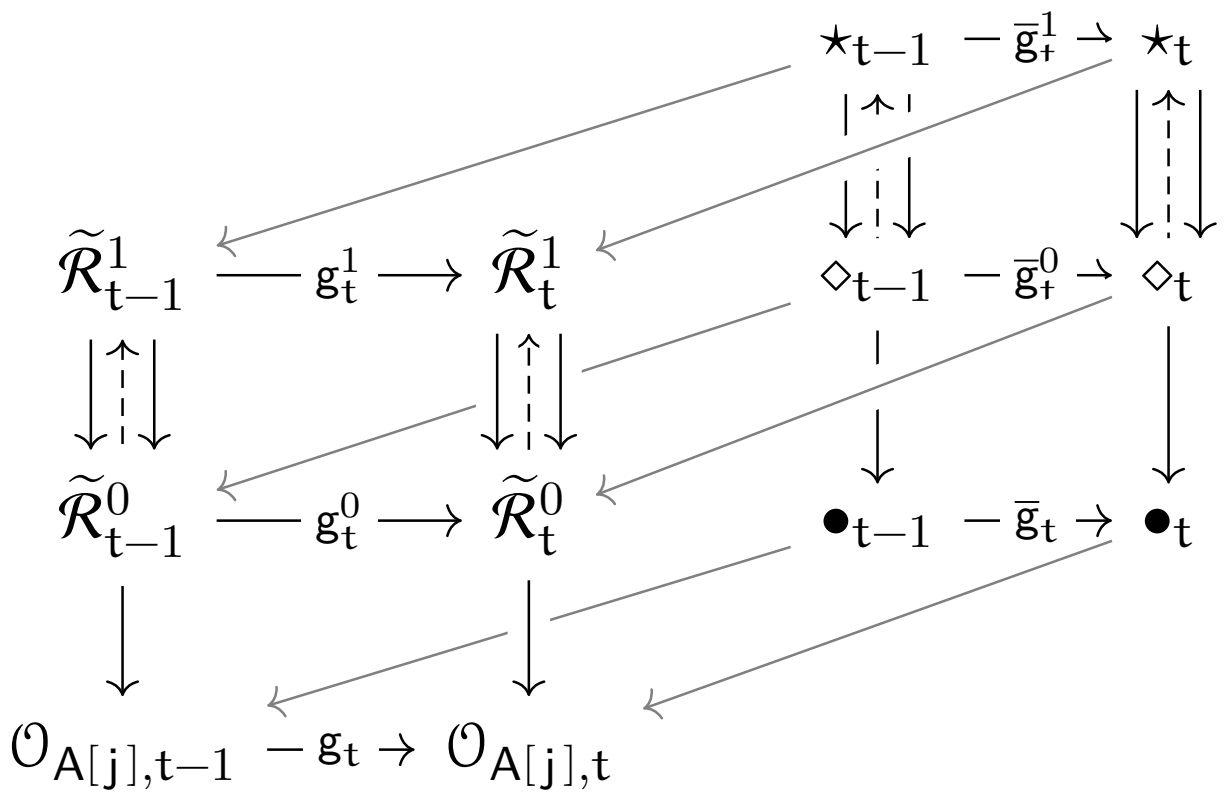

[i][30, Proposition 5.7] also provides a filtration of this form, but it is technically different. The issue concerns units in the operad.

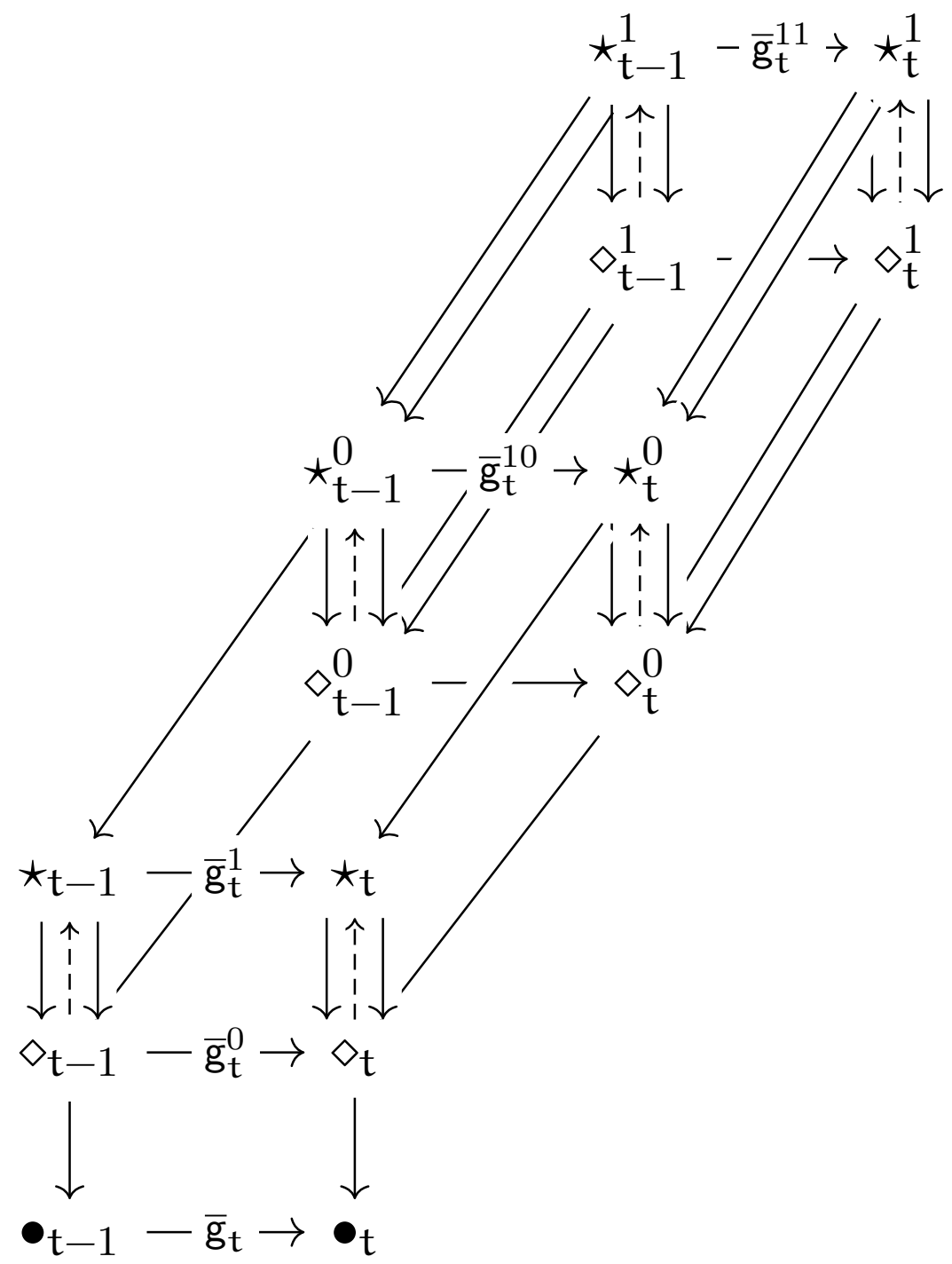

□

*Remark* 3.18. Since $\mathcal{O}_A\begin{bmatrix}0\\ \star\end{bmatrix} \cong A$, Proposition 3.17 gives back the filtration of the map $A \to A[j]$ in $\mathcal{V}^O$. Note that Pavlov-Scholbach in [30] and White-Yau in [40] deduce the filtration independently from that of $A \to A[j]$, while our sketch of proof implies that one can obtain a filtration of $\mathcal{O}_A \to \mathcal{O}_{A[j]}$ from one of $A \to A[j]$. Actually, the mentioned diagram with vertical and diagonal coequalizers that can be applied to compute $\overline{g}_t$ comes from the fact that in the filtration of $A \to A[j]$ the enveloping operad $\mathcal{O}_A$ plays a role.

*Remark* 3.19. In Propositions 3.14 and 3.17, we were not concerned with the operad structure of the resulting symmetric sequences. In all cases, the operadic structure is readily identified; in fact, composition products are induced by a combination of grafting and contractions of the tagged trees. It is also possible to describe explicitly the attaching maps in the filtrations.

**Slightly improving the filtrations.** Now we will discuss a generalization of Proposition 3.17 which will be important for applications, e.g. when comparing categorical and operadic left Kan extensions.

**Lemma 3.20.** *Let*

$$\begin{array}{ccc} \mathcal{O}\circ X & \longrightarrow & A \\ {\scriptstyle \mathrm{id}\circ j}\downarrow & \ulcorner & \downarrow{\scriptstyle g} \\ \mathcal{O}\circ Y & \longrightarrow & A[j] \end{array}$$

*be a pushout in* $\mathtt{Alg}_{\mathcal{O}}(\mathtt{V})$*. Then, the induced filtration on* $g_*\colon \mathcal{O}_A \to \mathcal{O}_{A[j]}$ *(seen as a map in* $\Sigma\mathtt{Coll}_{\mathsf{O}}(\mathtt{V})$*) obtained in Proposition 3.17 lifts to a filtration in* $\left[\mathsf{Fin}_{\mathsf{O}}^{\simeq,\mathrm{op}} \times \overline{\mathcal{O}}, \mathtt{V}\right]$*. In other words, for any* $\underline{a} \in \mathsf{Fin}_{\mathsf{O}}^{\simeq}$*, the morphism*

$$g_*\begin{bmatrix}\underline{a}\\ \star\end{bmatrix}\colon \mathcal{O}_A\begin{bmatrix}\underline{a}\\ \star\end{bmatrix} \to \mathcal{O}_{A[j]}\begin{bmatrix}\underline{a}\\ \star\end{bmatrix}$$

*(seen as a map in* $\mathtt{Alg}_{\overline{\mathcal{O}}}(\mathtt{V})$*) is the transfinite composite of a sequence* $(g_{*,t}\left[\begin{smallmatrix}\underline{a}\\ \star\end{smallmatrix}\right])_{t\in\mathbb{N}}$ *defined by pushouts in* $\mathtt{Alg}_{\overline{\mathcal{O}}}(\mathtt{V})$

$$\begin{array}{ccc} \bullet & \longrightarrow & \mathcal{O}_{A[j],t-1}\left[\begin{smallmatrix}\underline{a}\\ \star\end{smallmatrix}\right] \\ {\scriptstyle \widetilde{g}_{*,t}\left[\begin{smallmatrix}\underline{a}\\ \star\end{smallmatrix}\right]}\downarrow & \ulcorner & \downarrow{\scriptstyle g_{*,t}\left[\begin{smallmatrix}\underline{a}\\ \star\end{smallmatrix}\right]} \\ \bullet & \longrightarrow & \mathcal{O}_{A[j],t}\left[\begin{smallmatrix}\underline{a}\\ \star\end{smallmatrix}\right]. \end{array},$$

*where* $\widetilde{g}_{*,t}\left[\begin{smallmatrix}\underline{a}\\ \underline{b}\end{smallmatrix}\right] = \mathrm{colim}_\tau\,(\mathcal{O}_A; j)(\tau)$ *with* $\tau$ *running over the subcategory of* $(\mathsf{Tree}^{\simeq}_{\mathsf{O},\not\downarrow}\left[\begin{smallmatrix}\underline{a}\\ \underline{b}\end{smallmatrix}\right])^{\mathrm{op}}$ *defined in Proposition 3.17.*

*Remark* 3.21. Consider the situation described in Lemma 3.20. When $j\colon X \to Y$ is concentrated in one color, $c \in \mathsf{O}$, the morphism $\widetilde{g}_{*,t}$ of $\left[\mathsf{Fin}_{\mathsf{O}}^{\simeq,\mathrm{op}} \times \overline{\mathcal{O}}, \mathtt{V}\right]$ admits a workable explicit formula. Evaluated on $\underline{a} \in \mathsf{Fin}_{\mathsf{O}}^{\simeq}$, it becomes

$$\widetilde{g}_{*,t}\left[\begin{smallmatrix}\underline{a}\\ \star\end{smallmatrix}\right]\colon \mathcal{O}_A\left[\begin{smallmatrix}\underline{a}\boxplus c^{\boxplus t}\\ \star\end{smallmatrix}\right] \underset{\Sigma_t}{\otimes} s(j^{\square t}) \xrightarrow{\ \mathrm{id}\otimes_{\Sigma_t} j^{\square t}\ } \mathcal{O}_A\left[\begin{smallmatrix}\underline{a}\boxplus c^{\boxplus t}\\ \star\end{smallmatrix}\right] \underset{\Sigma_t}{\otimes} Y^{\otimes t}.$$

*Proof of Lemma 3.20.* The result holds in $\Sigma\mathtt{Coll}_{\mathsf{O}}(\mathtt{V}) = \left[\mathsf{Fin}_{\mathsf{O}}^{\simeq,\mathrm{op}} \times \mathsf{O}, \mathtt{V}\right]$, i.e. without functoriality over $\overline{\mathcal{O}}$, by Proposition 3.17. Hence, the contribution of this lemma is to show this functoriality. Due to the pointwise-computation of colimits in functor categories, the proof consists of showing that the pushouts in the statement (which by the cited proposition belong to $\mathtt{V}^{\mathsf{O}}$) live in $\mathtt{Alg}_{\overline{\mathcal{O}}}(\mathtt{V})$. For ease of notation, we are going to assume that $j$ is concentrated in the color $c \in \mathsf{O}$. The general case is analogous.

Notice that the vertices are $\overline{\mathcal{O}}$-algebras/functors by induction and since $\mathcal{O}_A$ is a $\mathcal{O}$-module through the canonical map of operads $\mathcal{O} \to \mathcal{O}_A$. So, it remains to show that the edges of the squares live in $\mathtt{Alg}_{\overline{\mathcal{O}}}(\mathtt{V})$, i.e.

$$\begin{array}{ccc} \mathcal{O}_A \begin{bmatrix} a \boxplus c^{\boxplus t} \\ \star \end{bmatrix} \underset{\Sigma_t}{\otimes} s(j^{\square t}) & & \mathcal{O}_A \begin{bmatrix} a \boxplus c^{\boxplus t} \\ \star \end{bmatrix} \underset{\Sigma_t}{\otimes} s(j^{\square t}) \\ \Big\downarrow \text{attach}_t & \text{and} & \Big\downarrow \text{id} \otimes_{\Sigma_t} j^{\square t} \\ \mathcal{O}_{A[j],t-1} \begin{bmatrix} a \\ \star \end{bmatrix} & & \mathcal{O}_A \begin{bmatrix} a \boxplus c^{\boxplus t} \\ \star \end{bmatrix} \underset{\Sigma_t}{\otimes} Y^{\otimes t} \end{array}$$

are compatible with the $\overline{\mathcal{O}}$-actions.

- $\text{id} \otimes_{\Sigma_t} j^{\square t}$ is compatible with the $\overline{\mathcal{O}}$-action because the structure comes from the first tensor factor and the map on this factor is the identity.
- $\text{attach}_n$ is compatible with $\overline{\mathcal{O}}$-action: the explicit description of the attaching map given in [14, Proposition 7.12] can be generalized to the colored-case yielding a decomposition of $\text{attaching}_n$. First note that $s(j^{\square t})$ is the colimit of a punctured $t$-cube diagram constructed out of j, whose vertices are tensor products of $p$ copies of X and $q$ copies of Y such that $p + q = t$ and $p > 0$. With this in mind, it is easy to see that $\text{attaching}_t$ factors as

$$\begin{array}{c} \mathcal{O}_A \begin{bmatrix} a \boxplus c^{\boxplus (p+q)} \\ \star \end{bmatrix} \underset{\Sigma_p \times \Sigma_q}{\otimes} X^{\otimes p} \otimes Y^{\otimes q} \\ \Big\downarrow \text{adjoint map to } \mathcal{O} \circ X \to A \\ \mathcal{O}_A \begin{bmatrix} a \boxplus c^{\boxplus (p+q)} \\ \star \end{bmatrix} \underset{\Sigma_p \times \Sigma_q}{\otimes} A(c)^{\otimes p} \otimes Y^{\otimes q} \\ \Big\| \\ \mathcal{O}_A \begin{bmatrix} a \boxplus c^{\boxplus (p+q)} \\ \star \end{bmatrix} \underset{\Sigma_p \times \Sigma_q}{\otimes} \mathcal{O}_A \begin{bmatrix} \emptyset \\ c \end{bmatrix}^{\otimes p} \otimes Y^{\otimes q} \\ \Big\downarrow \text{operadic composition} \\ \mathcal{O}_A \begin{bmatrix} a \boxplus c^{\boxplus q} \\ \star \end{bmatrix} \underset{\Sigma_q}{\otimes} Y^{\otimes q} \\ \Big\downarrow \text{inductively defined map for } q \\ \mathcal{O}_{A[j],q} \begin{bmatrix} a \\ \star \end{bmatrix} \\ \Big\downarrow g_{t-1} \cdots g_{q+1} \\ \mathcal{O}_{A[j],t-1} \begin{bmatrix} a \\ \star \end{bmatrix} \end{array} .$$

  It suffices to check that each morphism in the composite is compatible with the $\overline{\mathcal{O}}$-action. The first ones does not involve the tensor factor where $\overline{\mathcal{O}}$-acts; the following one is given by right multiplication in $\mathcal{O}$ and hence, by associativity, commutes with left multiplication (which is the one that defines the $\overline{\mathcal{O}}$-action); compatibility holds for the remaining maps by induction.

□

**Lemma 3.22.** *Let $\varphi\colon \mathcal{O} \to \mathcal{P}$ be a morphism of operads. Let $g\colon A \to A[j]$ be the map of $\mathcal{O}$-algebras in Lemma 3.20. Then, the natural map*

$$\varphi_!\, \mathcal{O}_{A[j]} \begin{bmatrix} \underline{a} \\ \star \end{bmatrix} \to \mathcal{P}_{\varphi_\sharp A[j]} \begin{bmatrix} \varphi(\underline{a}) \\ \star \end{bmatrix}$$

*is the transfinite colimit of the $\omega$-ladder*

$$\begin{array}{ccccccccc}
\varphi_!\, \mathcal{O}_A \begin{bmatrix} \underline{a} \\ \star \end{bmatrix} & \longrightarrow & \cdots & \longrightarrow & \varphi_!\, \mathcal{O}_{A[j],n-1} \begin{bmatrix} \underline{a} \\ \star \end{bmatrix} & \longrightarrow & \varphi_!\, \mathcal{O}_{A[j],n} \begin{bmatrix} \underline{a} \\ \star \end{bmatrix} & \longrightarrow & \cdots \\
\downarrow & & & & \downarrow{\scriptstyle \gamma_{n-1}} & & \downarrow{\scriptstyle \gamma_n} & & \\
\mathcal{P}_{\varphi_\sharp A} \begin{bmatrix} \varphi(\underline{a}) \\ \star \end{bmatrix} & \longrightarrow & \cdots & \longrightarrow & \mathcal{P}_{\varphi_\sharp A[j],n-1} \begin{bmatrix} \varphi(\underline{a}) \\ \star \end{bmatrix} & \longrightarrow & \mathcal{P}_{\varphi_\sharp A[j],n} \begin{bmatrix} \varphi(\underline{a}) \\ \star \end{bmatrix} & \longrightarrow & \cdots
\end{array}$$

*where $\gamma_n$ is defined inductively using the filtration obtained in Lemma 3.20.*

*Proof.* It suffices to define $\gamma_n$ since then, using that the functors $\varphi_!$ and evaluation at $\underline{a}$ preserve sequential colimits and Lemma 3.23, one deduces the claim. Again, we assume for ease of notation that $j$ is concentrated in color $c \in O$.

To define $\gamma_n$, one observes that it fits (dashed arrow) into a cube

$$\begin{array}{ccccc}
 & \varphi_!\Big(\mathcal{O}_A \begin{bmatrix} \underline{a} \boxplus c^{\boxplus n} \\ \star \end{bmatrix} \underset{\Sigma_n}{\otimes} s(j^{\square n})\Big) & \xrightarrow{\ \psi\ } & \mathcal{P}_{\varphi_\sharp A} \begin{bmatrix} \varphi(\underline{a} \boxplus c^{\boxplus n}) \\ \star \end{bmatrix} \underset{\Sigma_n}{\otimes} s(j^{\square n}) & \\
\swarrow{\scriptstyle \varphi_!(\mathrm{attach}_n)} & & \downarrow \qquad \swarrow{\scriptstyle \mathrm{attach}_n} & & \downarrow{\scriptstyle \mathrm{id} \otimes_{\Sigma_n} j^{\square n}} \\
\varphi_!\, \mathcal{O}_{A[j],n-1} \begin{bmatrix} \underline{a} \\ \star \end{bmatrix} & \xrightarrow{\ \gamma_{n-1}\ } & \mathcal{P}_{\varphi_\sharp A[j],n-1} \begin{bmatrix} \varphi(\underline{a}) \\ \star \end{bmatrix} & & \\
 & \varphi_!\Big(\mathcal{O}_A \begin{bmatrix} \underline{a} \boxplus c^{\boxplus n} \\ \star \end{bmatrix} \underset{\Sigma_n}{\otimes} Y^{\otimes n}\Big) & \xrightarrow{\ \phi\ } & \mathcal{P}_{\varphi_\sharp A} \begin{bmatrix} \varphi(\underline{a} \boxplus c^{\boxplus n}) \\ \star \end{bmatrix} \underset{\Sigma_n}{\otimes} Y^{\otimes n} & \\
\downarrow \qquad \swarrow & & \downarrow \qquad \swarrow & & \\
\varphi_!\, \mathcal{O}_{A[j],n} \begin{bmatrix} \underline{a} \\ \star \end{bmatrix} & \overset{\gamma_n}{\dashrightarrow} & \mathcal{P}_{\varphi_\sharp A[j],n} \begin{bmatrix} \varphi(\underline{a}) \\ \star \end{bmatrix} & &
\end{array}$$

whose left and right faces are pushout squares. Hence, the commutativity of the back square implies the existence of $\gamma_n$. This commutativity follows easily from the definition of $\psi$ and $\phi$; both are induced from the map of $\overline{\mathcal{P}}$-algebras

$$\varphi_!\, \mathcal{O}_A \begin{bmatrix} \underline{a} \\ \star \end{bmatrix} \longrightarrow \mathcal{P}_{\varphi_\sharp A} \begin{bmatrix} \varphi(\underline{a}) \\ \star \end{bmatrix}.$$

□

For the record, we spell out the arity 0 version of Lemmas 3.20 and 3.22.

**Lemma 3.23.** *Let*

$$\begin{array}{ccc}
\mathcal{O} \circ X & \longrightarrow & A \\
{\scriptstyle \mathrm{id} \circ j}\downarrow & \ulcorner & \downarrow{\scriptstyle g} \\
\mathcal{O} \circ Y & \longrightarrow & A[j]
\end{array}$$

*be a pushout in* $\mathtt{Alg}_{\mathcal{O}}(\mathtt{V})$*. Then, the underlying arrow associated to* $g$ *in* $\mathtt{Alg}_{\overline{\mathcal{O}}}(\mathtt{V})$ *is the transfinite composition of a sequence* $(g_n)_{n\in\mathbb{N}}$ *defined by pushouts in* $\mathtt{Alg}_{\overline{\mathcal{O}}}(\mathtt{V})$*:*

$$\begin{array}{ccc} \bullet & \longrightarrow & A[j]_{n-1} \\ \coprod_{|\underline{a}|=n}\Big(\mathcal{O}_A\left[\begin{smallmatrix}\underline{a}\\ \star\end{smallmatrix}\right] \otimes_{\mathrm{Aut}(\underline{a})} \square_r\, j(a_r)\Big)\Big\downarrow & \ulcorner & \Big\downarrow g_n \\ \bullet & \longrightarrow & A[j]_n \end{array} .$$

**Lemma 3.24.** *Let* $\varphi\colon \mathcal{O} \to \mathcal{P}$ *be a morphism of operads. Let* $g\colon A \to A[j]$ *be the map of* $\mathcal{O}$*-algebras in Lemma 3.23. Then, the natural map* $\varphi_!\overline{A[j]} \to \overline{\varphi_\sharp A[j]}$ *is the transfinite colimit of the* $\omega$*-ladder*

$$\begin{array}{ccccccccc} \varphi_!\overline{A} & \longrightarrow & \cdots & \longrightarrow & \varphi_!\overline{A[j]}_{n-1} & \xrightarrow{\varphi_!(\overline{g}_n)} & \varphi_!\overline{A[j]}_n & \longrightarrow & \cdots \\ \Big\downarrow & & & & \gamma_{n-1}\Big\downarrow & & \Big\downarrow\gamma_n & & \\ \overline{\varphi_\sharp A} & \longrightarrow & \cdots & \longrightarrow & \big(\overline{\varphi_\sharp A[j]}\big)_{n-1} & \xrightarrow[\overline{(\varphi_\sharp g)}_n]{} & \big(\overline{\varphi_\sharp A[j]}\big)_n & \longrightarrow & \cdots \end{array}$$

*where* $\gamma_n$ *is defined inductively using the natural transformation* $\varphi_!\overline{A} \to \overline{\varphi_\sharp A}$*.*

## 3.3 Homotopical analysis

**Equivariant cofibrancy.** The whole previous discussion about enveloping operads makes use of trees and isomorphisms of trees. In order to understand homotopical properties of enveloping operads, it is then useful to have some control on $\mathrm{Aut}(\Upsilon)$-equivariant homotopy theory over $\mathcal{V}$ for trees $\Upsilon$ or more generally, control over model categories as $[(\mathsf{Tree}^{|v|,\simeq}_{O,\not\downarrow})^{op}, \mathcal{V}]$. We will assume that diagram categories as $\mathcal{V}^G$ all carry the projective-model structure in this section, where $G$ is a group. We will call (acyclic) $G$-cofibrations to the (acyclic) proj-cofibrations in this setting.

*Remark* 3.25. Even though we are working with colored operads, and hence we have to deal with groupoids to encode equivariance, one may fix a skeleton of those groupoids and just focus on groups. We are going to state all the results using groups because of that reason and because they are easier to state. Classical references such as [3, 4] and [14] essentially cover the following material about equivariant cofibrancy.

Following Spitzweck (see [36, Lemma 3.6 and Proof of 3.7]) or [3], we observe that any decomposition of $\Upsilon$ as grafting$(\tau; \Upsilon'_1, \ldots, \Upsilon'_m)$, with $\tau$ being a corolla and the trees $\Upsilon'_j$ being of types $\{\Upsilon_1, \ldots \Upsilon_k\}$ where $\Upsilon_i$ appears $n_i$ times, yields an identification

$$\mathrm{Aut}(\Upsilon) \cong \Big(\prod_i \Sigma_{n_i}\Big) \ltimes \Big(\prod_i \mathrm{Aut}(\Upsilon_i)^{\times n_i}\Big).$$

This formula can be applied inductively on the depth of a tree to decompose its group of automorphisms in terms of (semi)direct products of symmetric groups.

Also note that $M \otimes L^{\otimes n}$ carries an action of the wreath product $\Sigma_n \ltimes G^{\times n}$ when $\Sigma_n$ acts on $M$ and $G$ acts on $L$.

These two simple observations together with the following list of lemmas will be all what we need about $\mathrm{Aut}(\Upsilon)$-equivariant homotopy theory over $\mathcal{V}$ in the sequel.

**Lemma 3.26.** *Let* $\mathsf{H} \subseteq \mathsf{G}$ *be a subgroup. Then the restriction* $\mathsf{rest}\colon \mathcal{V}^{\mathsf{G}} \to \mathcal{V}^{\mathsf{H}}$ *is left Quillen.*

*Proof.* It suffices to check that rest sends generating (acyclic) G-cofibrations to H-cofibrations. These generating maps are of the form $\mathsf{G} \cdot \mathsf{i}$, for i a generating (acyclic) cofibration in $\mathcal{V}$. Then, just observe that $\mathsf{rest}(\mathsf{G} \cdot \mathsf{i}) = \coprod_{\mathsf{G}/\mathsf{H}} \mathsf{H} \cdot \mathsf{i}$ is a coproduct of H-cofibrations. □

**Lemma 3.27.** *Let* $1 \to \mathsf{G}_1 \to \mathsf{G} \to \mathsf{G}_2 \to 1$ *be a short exact sequence of groups. Then, taking* $\mathsf{G}_1$*-coinvariants,* $(\star)_{\mathsf{G}_1}\colon \mathcal{V}^{\mathsf{G}} \to \mathcal{V}^{\mathsf{G}_2}$*, is a left Quillen functor.*

*Proof.* Similar to Lemma 3.26. Note that $(\mathsf{G} \cdot \mathsf{i})_{\mathsf{G}_1} = (\mathsf{G}/\mathsf{G}_1) \cdot \mathsf{i} = \mathsf{G}_2 \cdot \mathsf{i}$. □

**Lemma 3.28.** *The tensor product* $\otimes\colon \mathcal{V}^{\mathsf{G}_1} \times \mathcal{V}^{\mathsf{G}_2} \to \mathcal{V}^{\mathsf{G}_1 \times \mathsf{G}_2}$ *is a left Quillen bifunctor.*

*Proof.* By an argument similar to that of Lemma 3.26, one should just notice

$$(\mathsf{G}_1 \cdot \mathsf{i}) \square (\mathsf{G}_2 \cdot \mathsf{j}) = (\mathsf{G}_1 \times \mathsf{G}_2) \cdot (\mathsf{i} \square \mathsf{j}).$$

□

**Lemma 3.29.** *Let* $\mathsf{H} \subseteq \mathsf{G}$ *be a subgroup. Then, the* H*-invariant internal-hom*

$$\underline{\mathsf{Hom}}(\star, \star)^{\mathsf{H}}\colon (\mathcal{V}^{\mathsf{G}})^{\mathsf{op}} \times \mathcal{V}^{\mathsf{G}} \to \mathcal{V}$$

*is a right Quillen bifunctor.*

*Proof.* Take a G-cofibration i and a G-fibration p. We need to show that the pullback product map $\widehat{\underline{\mathsf{Hom}}}(\mathsf{i}, \mathsf{p})^{\mathsf{H}}$ is a fibration. Let us see that this map has rlp against cofibrations in $\mathcal{V}$. Considering a cofibration j in $\mathcal{V}$, one just has to use the equivalence between lifting properties

$$\mathsf{j} \boxslash \widehat{\underline{\mathsf{Hom}}}(\mathsf{i}, \mathsf{p})^{\mathsf{H}} \quad \Leftrightarrow \quad \mathsf{i} \boxslash \widehat{\underline{\mathsf{Hom}}}\big((\mathsf{G}/\mathsf{H}) \cdot \mathsf{j}, \mathsf{p}\big)$$

since $(\mathsf{G}/\mathsf{H}) \cdot \mathsf{j}$ is a coproduct of cofibrations and $\underline{\mathsf{Hom}}(\star, \star)$ is a right Quillen bifunctor. The acyclic cases are similar. □

**Lemma 3.30.** *Any finite tensor product functor* $\otimes^{(m)}\colon \prod_{\mathsf{i}} \mathcal{V}^{\mathsf{G}_{\mathsf{i}}} \to \mathcal{V}^{\prod_{\mathsf{i}} \mathsf{G}_{\mathsf{i}}}$ *preserves core cofibrations. If all the maps are core acyclic cofibrations, their tensor product is so.*

*Proof.* By induction and symmetry, it suffices to see that $\mathsf{j} \otimes \mathsf{X}$ is (acyclic) cofibration if j is a core (acyclic) $\mathsf{G}_1$-cofibration and X is $\mathsf{G}_2$-cofibrant. It is enough to check the llp against acyclic fibrations (resp. fibrations) using the equivalence between lifting properties

$$(\mathsf{j} \otimes \mathsf{X}) \boxslash \mathsf{p} \quad \Leftrightarrow \quad \mathsf{j} \boxslash \underline{\mathsf{Hom}}(\mathsf{X}, \mathsf{p})^{\mathsf{G}_2}.$$

Hence, to finish the proof, use Lemma 3.29. □

**Lemma 3.31.** *Let* i, j *be* fgt*-cofibrations in* $\mathcal{V}^{\mathsf{G}}$*. Then, the pushout product* $\mathsf{i} \square \mathsf{j}$ *is a* G*-cofibration if* i *or* j *is a* G*-cofibration. Furthermore, if additionally* i *or* j *is a weak equivalence,* $\mathsf{i} \square \mathsf{j}$ *is an acyclic* G*-cofibration.*

*Proof.* Again, we should check the llp against acyclic fibrations. Use

$$(i \square j) \boxslash p \quad \Leftrightarrow \quad i \boxslash \underline{\widehat{\mathrm{Hom}}}(j, p)$$

assuming that i is the G-cofibration, since acyclic G-fibrations are the same as acyclic fgt-fibrations. The acyclic case is analogous. □

**Lemma 3.32.** *Let $1 \to G_1 \to G \to G_2 \to 1$ be a short exact sequence of groups. Consider a $G_2$-cofibration $i$ and a map $j$ in $\mathcal{V}^G$ which is a $G_1$-cofibration. Then, the pushout product $i \square j$ is a G-cofibration and it is acyclic if $i$ or $j$ is so.*

*Proof.* Use the equivalence between lifting properties

$$(i \square j) \boxslash p \quad \Leftrightarrow \quad i \boxslash \underline{\widehat{\mathrm{Hom}}}(j, p)^{G_2}.$$

□

*Remark* 3.33. By a simple induction procedure, it is possible to generalize some of the results above for finite families of maps.

**Homotopical analysis.** The filtrations obtained in §3 can be employed to show fundamental properties of $(\mathcal{O}, A) \mapsto \mathcal{O}_A$. Let us start exploring cofibrancy properties of $\mathcal{O}_A$ depending on those of $\mathcal{O}$ and A.

**Proposition 3.34.** *The enveloping operad $\mathcal{O}_A$ satisfies:*

(i) *$\mathcal{O}_A$ is $\Sigma$-cofibrant (resp. well-pointed) if $\mathcal{O}$ is $\Sigma$-cofibrant (resp. well-pointed) and A is proj-cofibrant.*

(ii) *$\mathcal{O}_A$ is well-pointed if $\mathcal{O}$ is cofibrant and A is fgt-cofibrant.*

The proof of this result makes use of the following pair of lemmas.

**Lemma 3.35.** *Consider a pushout square in $\mathcal{A}\mathrm{lg}_{\mathcal{O}}(\mathcal{V})$*

$$\begin{array}{ccc} \mathcal{O} \circ X & \longrightarrow & A \\ {\scriptstyle \mathcal{O} \circ j}\downarrow & \ulcorner & \downarrow{\scriptstyle g} \\ \mathcal{O} \circ Y & \longrightarrow & A[j] \end{array}$$

*with $j \colon X \to Y$ being (acyclic) cofibration in $\mathcal{V}^{\mathrm{O}}$. Then, $\mathcal{O}_A \to \mathcal{O}_{A[j]}$ is a(n acyclic) $\Sigma$-cofibration if $\mathcal{O}_A$ is $\Sigma$-cofibrant or well-pointed.*

*Proof.* Using the filtration obtained in Proposition 3.17 and since (acyclic) $\Sigma$-cofibrations are closed under pushouts, it suffices to see that

$$\begin{bmatrix} c \\ d \end{bmatrix} \longmapsto \underset{\tau}{\mathrm{colim}}\, (\mathcal{O}_A; j)(\tau),$$

with $\tau$ being corollas in $\mathsf{Tree}^{\cong}_{O,\not\downarrow}\left[\begin{smallmatrix}\underline{c}\\ d\end{smallmatrix}\right]$ with $t$ straight leaves, is a(n acyclic) $\Sigma^{op}_{\underline{c}}$-cofibration. We can assume without loss of generality that $j$ is a generating (acyclic) cofibration in $\mathcal{V}^{O}$ concentrated in one color $b \in O$. Thus, we find

$$\operatorname*{colim}_{\tau}\,(\mathcal{O}_A;\, j)(\tau) \cong \left(\mathcal{O}_A \begin{bmatrix} s(\tau) \\ d \end{bmatrix} \otimes j \begin{bmatrix} \mathbb{0} \\ b \end{bmatrix}^{\square t}\right)_{\mathrm{Aut}(\tau)^{op}}$$

since there is only one isomorphism class in the subgroupoid of corollas $\tau$ in $\mathsf{Tree}^{\cong}_{O,\not\downarrow}\left[\begin{smallmatrix}\underline{c}\\ d\end{smallmatrix}\right]$ with $t$ straight leaves labelled by $b \in O$. Furthermore, $\mathrm{Aut}(\tau)^{op} \cong \Sigma^{op}_{t}$. Since $j$ is a(n acyclic) cofibration, $j\left[\begin{smallmatrix}\mathbb{0}\\ b\end{smallmatrix}\right]^{\square t}$ is a(n acyclic) cofibration by the pushout product axiom. Its tensor product with $\mathcal{O}_A\left[\begin{smallmatrix}s(\tau)\\ d\end{smallmatrix}\right]$ is a(n acyclic) $\Sigma^{op}_{\underline{c}} \times \Sigma^{op}_{t}$-cofibration because $\mathcal{O}_A$ is $\Sigma$-cofibrant by Lemma 3.28. For the well-pointed case, notice that $X \otimes \star\colon \mathcal{V} \to \mathcal{V}$ is left Quillen for any $\mathbb{I}$-cofibrant object $X \in \mathcal{V}$ (see Definition A.12 and Remark A.14). We conclude the proof using Lemma 3.27. $\square$

**Lemma 3.36.** *Let*

$$\begin{array}{ccc} \mathcal{F}(\mathcal{X}) & \longrightarrow & \mathcal{O} \\ {\scriptstyle \mathcal{F}j}\big\downarrow & \ulcorner & \big\downarrow{\scriptstyle g} \\ \mathcal{F}(\mathcal{Y}) & \longrightarrow & \mathcal{O}[j] \end{array}$$

*be a pushout of operads with* $j\colon \mathcal{X} \to \mathcal{Y}$ *being (acyclic) cofibration in* $\Sigma\,\mathcal{C}\mathsf{oll}_O(\mathcal{V})$ *and let* $A$ *be an* $\mathcal{O}[j]$*-algebra. Then,* $\mathcal{O}_A \to \mathcal{O}[j]_A$ *is a(n acyclic)* $\Sigma$*-cofibration if* $\mathcal{O}_A$ *is* $\Sigma$*-cofibrant or well-pointed.*

*Proof.* Using the filtration obtained in Proposition 3.14 and since (acyclic) $\Sigma$-cofibrations are closed under pushouts, it suffices to see that

$$\begin{bmatrix} \underline{c} \\ d \end{bmatrix} \longmapsto \operatorname*{colim}_{\Lambda}\,(\mathcal{O}_A; j; A)(\Lambda),$$

with $\Lambda$ running over the groupoid of trees in Proposition 3.14, i.e. trees with all external vertices being even, $t$ odd vertices (all of them without leaves), all leaves snaky and containing no forbidden subtrees, is (acyclic) $\Sigma^{op}_{\underline{c}}$-cofibration. We can assume without loss of generality that $j$ is a generating (acyclic) cofibration in $\Sigma\,\mathcal{C}\mathsf{oll}_O(\mathcal{V})$ concentrated in one profile $\left[\begin{smallmatrix}\underline{a}\\ b\end{smallmatrix}\right]$. Thus, we find

$$\operatorname*{colim}_{\Lambda}\,(\mathcal{O}_A;\, j\,;A)(\Lambda) \cong \coprod_{[\Lambda]} \left(\bigotimes_{\nu \in V_{\mathrm{even}}} \mathcal{O}_A \begin{bmatrix} s(\nu) \\ t(\nu) \end{bmatrix} \otimes j \begin{bmatrix} \underline{a} \\ b \end{bmatrix}^{\square t}\right)_{\mathrm{Aut}(\Lambda)^{op}}.$$

In this case, $j\left[\begin{smallmatrix}\underline{a}\\ b\end{smallmatrix}\right]^{\square t}$ is an (acyclic) $(\mathrm{Aut}(\underline{a})^{op})^{\times t}$-cofibration and the object $\bigotimes_{\nu} \mathcal{O}_A\left[\begin{smallmatrix}s(\nu)\\ t(\nu)\end{smallmatrix}\right]$ is $\prod_{\nu} \mathrm{Aut}(\nu)^{op}$-cofibrant by Lemma 3.28. For the well-pointed case, use Remark A.14.

Applying Lemma 3.32, we deduce that the morphism $\bigotimes_{\nu} \mathcal{O}_A\left[\begin{smallmatrix}s(\nu)\\ t(\nu)\end{smallmatrix}\right] \otimes j\left[\begin{smallmatrix}\underline{a}\\ b\end{smallmatrix}\right]^{\square t}$ is a(n acyclic) $G^{op}$-cofibration for $G$ the group of automorphisms of $\Lambda$ which do not have to fix the leaves labeled by $\underline{c}$. Then, taking coinvariants over $\mathrm{Aut}(\Lambda)^{op}$, i.e. automorphisms of $\Lambda$ that do fix the leaves labeled by $\underline{c}$, we obtain a(n acyclic) $\Sigma^{op}_{\underline{c}}$-cofibration due to Lemma 3.27. $\square$

Actually, a simplified version of the argument given in Lemma 3.36 yields:

**Corollary 3.37.** *If $\mathcal{O} \in \mathcal{O}\mathrm{pd}_{\mathrm{O}}(\mathcal{V})$ is cofibrant, the map from the initial* O*-operad $\mathcal{I}_{\mathrm{O}} \to \mathcal{O}$ is a $\Sigma$-cofibration. In other words, $\mathcal{O}$ is well-pointed (and hence $\Sigma$-cofibrant if the monoidal unit $\mathbb{I} \in \mathcal{V}$ is cofibrant).*

*Proof of Proposition 3.34.* By the retract closure of $\Sigma$-cofibrations, in (i) one can take A to be cellular and in (ii) one can take $\mathcal{O}$ to be cellular. By the closure under transfinite composites, in both cases we are reduced to study just the step in the construction of cellular objects associated to attaching a cell. For (i) use Lemma 3.35 and for (ii) use Lemma 3.36. □

Let us now move on to study when $\mathcal{O}_\star \colon \mathrm{A} \mapsto \mathcal{O}_{\mathrm{A}}$ preserves equivalences.

**Proposition 3.38.** *The functor $\mathcal{O}_\star$ sends equivalences of algebras to equivalences of operads if one of the following conditions holds:*

(i) *$\mathcal{O}$ is $\Sigma$-cofibrant or well-pointed and we restrict $\mathcal{O}_\star$ to proj-cofibrant algebras.*

(ii) *$\mathcal{O}$ is cofibrant and we restrict $\mathcal{O}_\star$ to* fgt*-cofibrant algebras.*

*Proof.* We prove each statement separatedly.

(i) By Proposition 3.34 and the fact that equivalences in $\mathcal{O}\mathrm{pd}_{\mathrm{O}}(\mathcal{V})$ are defined objectwise, we may apply Ken Brown's lemma ([10, Lemma 12.1.6])[ii] to the composite

$$\begin{array}{ccccc} \mathcal{A}\mathrm{lg}_{\mathcal{O}}(\mathcal{V}) & \xrightarrow{\mathcal{O}_\star} & \mathcal{O}\mathrm{pd}_{\mathrm{O}}(\mathcal{V}) & \xrightarrow{\mathrm{fgt}_\Sigma} & \Sigma\,\mathcal{C}\mathrm{oll}_{\mathrm{O}}(\mathcal{V}) \\ \mathrm{A} & \longmapsto & \mathcal{O}_{\mathrm{A}} & & \end{array}$$

in order to prove the result. Hence, we should check if $\mathcal{O}_\star$ sends core acyclic proj-cofibrations to acyclic $\Sigma$-cofibrations. The usual line of argument reduces this problem to the case of a core acyclic proj-cofibration $\mathrm{A} \to \mathrm{A}[\mathrm{j}]$ which is a cobase change of $\mathcal{O} \circ \mathrm{j}$, with j core acyclic cofibration in $\mathcal{V}^{\mathrm{O}}$ (acyclic cell attachment). The acyclic statement in Lemma 3.35 yields the result.

(ii) Now, one makes use of the filtration associated to a cellular operad instead. By the retract closure of equivalences and since sequential colimits of equivalences along cofibrations between cofibrant objects are equivalences (applied to $\Sigma\,\mathcal{C}\mathrm{oll}_{\mathrm{O}}(\mathcal{V})$ here), we should study the step in the construction of cellular objects associated to attaching a cell. We prove a slightly more general statement in Lemma 3.39. We are implicitly applying Remark A.9 and the fact that fgt preserves core cofibrations in this situation.

□

[ii] To be able to apply Ken's Brown lemma, we are implicitly using Proposition 4.2. It is also possible to produce a direct proof avoiding this lemma, but the one we present is much shorter.

**Lemma 3.39.** *Let*

$$\begin{array}{ccc} \mathcal{F}(\mathcal{X}) & \longrightarrow & \mathcal{O} \\ {\scriptstyle \mathcal{F}j}\downarrow & \ulcorner & \downarrow{\scriptstyle g} \\ \mathcal{F}(\mathcal{Y}) & \longrightarrow & \mathcal{O}[j] \end{array}$$

*be a pushout of operads with* $j\colon \mathcal{X} \to \mathcal{Y}$ *being a core cofibration in* $\Sigma\,\mathcal{C}\mathrm{oll}_{\mathrm{O}}(\mathcal{V})$ *and let* $A \to B$ *be a morphism of* $\mathcal{O}[j]$*-algebras. Then, if* $\mathcal{O}_A \to \mathcal{O}_B$ *is an equivalence of* $\Sigma$*-cofibrant or well-pointed operads, so is* $\mathcal{O}[j]_A \to \mathcal{O}[j]_B$.

*Proof.* Using Proposition 3.14, we may filter the commutative square

$$\begin{array}{ccc} \mathcal{O}_A & \longrightarrow & \mathcal{O}[j]_A \\ \downarrow & & \downarrow \\ \mathcal{O}_B & \longrightarrow & \mathcal{O}[j]_B \end{array}$$

in $\Sigma\,\mathcal{C}\mathrm{oll}_{\mathrm{O}}(\mathcal{V})$ by the ladder

$$\begin{array}{ccccccc} \mathcal{O}_A = \mathcal{O}[j]_{A,0} & \longrightarrow & \mathcal{O}[j]_{A,1} & \longrightarrow \cdots \longrightarrow & \mathcal{O}[j]_{A,\omega} = \mathcal{O}[j]_A \\ \downarrow & & \downarrow & & \downarrow \\ \mathcal{O}_B = \mathcal{O}[j]_{B,0} & \longrightarrow & \mathcal{O}[j]_{B,1} & \longrightarrow \cdots \longrightarrow & \mathcal{O}[j]_{B,\omega} = \mathcal{O}[j]_B \end{array} ,$$

where each step $(n-1) \Rightarrow (n)$ (evaluated at $\left[\begin{smallmatrix} c \\ d \end{smallmatrix}\right]$) fits into a commutative cube whose black faces are pushouts. These pushouts are moreover hopushouts since $\mathcal{O}[j]_{A,t}$ (resp. $\mathcal{O}[j]_{B,t}$) is $\Sigma$-cofibrant for all $t \geqslant 0$ plus the fact that

$$\begin{bmatrix} c \\ d \end{bmatrix} \longmapsto \operatorname*{colim}_{\Lambda}(\mathcal{O}_A; j; A)(\Lambda)$$

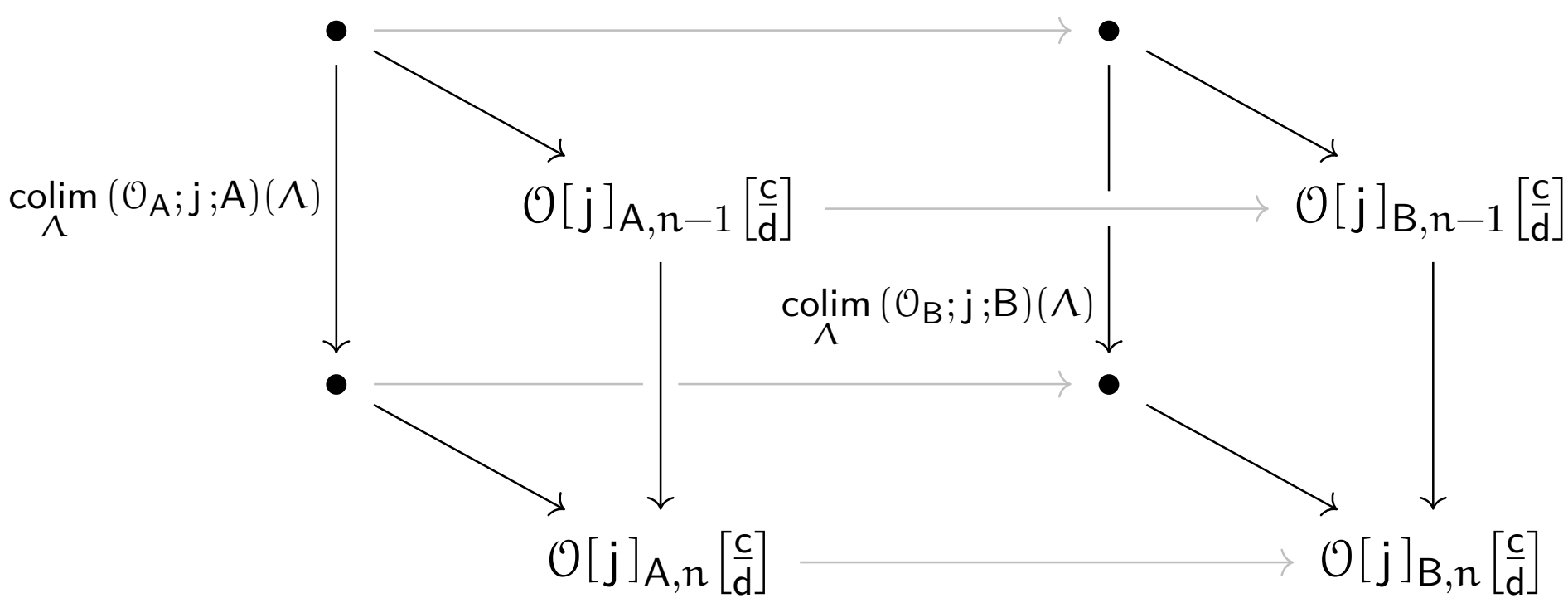

(resp. for $(\mathcal{O}_B; j; B)$) is a core $\Sigma_{\underline{c}}$-cofibration (see the proof of Lemma 3.36 applying additionally that j is a core $\Sigma$-cofibration). Thus, we should just check that the gray horizontal maps in the cube (except the lower one in the front face) are equivalences. This follows by induction and because $(\mathcal{O}_A; j; A)(\Lambda) \to (\mathcal{O}_B; j; B)(\Lambda)$ is an equivalence between $\mathrm{Aut}(\Lambda)^{\mathrm{op}}$-cofibrant objects since $\mathcal{O}_A \to \mathcal{O}_B$ is so and j is a core $\Sigma$-cofibration. If one the operads is well-pointed, let us say $\mathcal{O}_A$, $\Lambda \mapsto (\mathcal{O}_A; j; A)(\Lambda) \to (\mathcal{O}_B; j; B)(\Lambda)$ is still $\mathrm{Aut}(\Lambda)^{\mathrm{op}}$-cofibrant by Remark A.14. □

It remains to study when the assignment $\mathcal{O} \mapsto \mathcal{O}_A$, suitably understood, preserves equivalences.

**Proposition 3.40.** *Let $(\mathcal{O}, A) \to (\mathcal{P}, B)$ be a map of operadic algebras, i.e. a map in* $\mathtt{OpdAlg(V)}$*, such that its component $\phi\colon \mathcal{O} \to \mathcal{P}$ is an equivalence in $\mathcal{O}\mathrm{pd}_O(\mathcal{V})$. Then, the morphism $\mathcal{O}_A \to \mathcal{P}_B$ is an equivalence provided one of the following conditions holds:*

(i) *the map between $\mathcal{P}$-algebras is $\varphi_\sharp A = \varphi_\sharp A$ (equiv. the unit $A \to \phi^*\phi_\sharp A$ in $\mathcal{O}$-algebras) for a proj-cofibrant $\mathcal{O}$-algebra $A$ and either:*

- *$\mathcal{O}$, $\mathcal{P}$ are $\Sigma$-cofibrant operads, or*
- *$\mathcal{O}, \mathcal{P} \in \{\Sigma$-cofibrant or well-pointed operads$\}$ and the $\mathbb{I}$-strong unit axiom holds.*

(ii) *the map between $\mathcal{O}$-algebras is $\phi^* B = \phi^* B$ (equiv. the counit $\phi_\sharp\phi^* B \to B$ in $\mathcal{P}$-algebras) for a* fgt*-cofibrant $\mathcal{P}$-algebra $B$ and $\mathcal{O}$, $\mathcal{P}$ are cofibrant operads.*

*Proof.* We prove each statement separatedly.

(i) One makes use of the filtration associated to a cellular algebra and the induced filtration of the enveloping operad. By the retract closure of equivalences and since sequential colimits of equivalences along cofibrations between cofibrant objects are equivalences[iii] (applied to $\Sigma\,\mathcal{C}\mathrm{oll}_O(\mathcal{V})$ here), we are reduced to study just the step in the construction of cellular algebras associated to attaching a cell. Lemma 3.41 covers such step.

(ii) We prove a slightly more general statement. Let $\mathcal{Q}$ be an operad under $\mathcal{O} \to \mathcal{P}$ and $B$ a fgt-cofibrant $\mathcal{Q}$-algebra. By Proposition 3.34, we may apply Ken Brown's lemma[iv] to the composite

$$\begin{array}{ccccc} \mathcal{O}\mathrm{pd}_O(\mathcal{V})\downarrow\mathcal{Q} & \xrightarrow{(\star)_B} & \mathcal{O}\mathrm{pd}_O(\mathcal{V}) & \xrightarrow{\mathrm{fgt}_\Sigma} & \Sigma\,\mathcal{C}\mathrm{oll}_O(\mathcal{V}) \\ (\mathcal{O} \to \mathcal{Q}) & \longmapsto & \mathcal{O}_B & & \end{array}$$

in order to prove the claim. Hence, it suffices to check that $(\star)_B$ sends core acyclic proj-cofibrations to acyclic $\Sigma$-cofibrations. The only non-trivial part is the analysis of acyclic cell attachments and such case is covered by the acyclic part of Lemma 3.36. We are implicitly applying Remark A.9.

∎

**Lemma 3.41.** *Let $\phi\colon \mathcal{O} \to \mathcal{P}$ be a map in $\mathcal{O}\mathrm{pd}_O(\mathcal{V})$ and consider a pushout square in $\mathcal{A}\mathrm{lg}_{\mathcal{O}}(\mathcal{V})$*

$$\begin{array}{ccc} \mathcal{O}\circ X & \longrightarrow & A \\ {\scriptstyle \mathcal{O}\circ j}\downarrow & \ulcorner & \downarrow{\scriptstyle g} \\ \mathcal{O}\circ Y & \longrightarrow & A[j] \end{array}$$

[iii] In the well-pointed case, similar statements hold true; see Lemma A.17.

[iv] To be able to apply this lemma, we are implicitly applying Proposition 4.2 to the colored operad in $\mathtt{Set}$ whose algebras are O-operads (alternatively use [10, Theorem 12.2.A]). It is also possible to produce a direct proof of the statement, but this one is much shorter.

*with* $j\colon X \to Y$ *being a core cofibration in* $\mathcal{V}^{O}$. *Then,* $\mathcal{O}_{A[j]} \to \mathcal{P}_{\phi_\sharp(A[j])}$ *is an equivalence of operads if so is* $\mathcal{O}_A \to \mathcal{P}_{\phi_\sharp A}$ *and either:*

- $\mathcal{O}$, $\mathcal{P}$ *are* $\Sigma$*-cofibrant operads, or*

- $\mathcal{O}, \mathcal{P} \in \{\Sigma$*-cofibrant or well-pointed operads*$\}$ *and the* $\mathbb{I}$*-strong unit axiom holds.*

*Proof.* Let us argue the $\Sigma$-cofibrant case.

First notice that the pushout in the statement implies that $\phi_\sharp A \to \phi_\sharp(A[j])$ is the cobase change of $\mathcal{P} \circ j$. Using Proposition 3.17, we may filter the commutative square

$$\begin{array}{ccc} \mathcal{O}_A & \longrightarrow & \mathcal{O}_{A[j]} \\ \downarrow & & \downarrow \\ \mathcal{P}_{\phi_\sharp A} & \longrightarrow & \mathcal{P}_{\phi_\sharp A[j]} \end{array}$$

in $\Sigma\,\mathcal{C}oll_O(\mathcal{V})$ by the ladder

$$\begin{array}{ccccccc} \mathcal{O}_A = \mathcal{O}_{A[j],0} & \longrightarrow & \mathcal{O}_{A[j],1} & \longrightarrow \cdots \longrightarrow & \mathcal{O}_{A[j],\omega} = \mathcal{O}_{A[j]} \\ \downarrow & & \downarrow & & \downarrow \\ \mathcal{P}_{\phi_\sharp A} = \mathcal{P}_{\phi_\sharp A[j],0} & \longrightarrow & \mathcal{P}_{\phi_\sharp A[j],1} & \longrightarrow \cdots \longrightarrow & \mathcal{P}_{\phi_\sharp A[j],\omega} = \mathcal{P}_{\phi_\sharp A[j]} \end{array}.$$

Each step $(n-1) \Rightarrow (n)$ (evaluated at $\left[\begin{smallmatrix} c \\ d \end{smallmatrix}\right]$) fits into a commutative cube whose black faces are pushouts. These pushouts are moreover homotopy pushouts since $\mathcal{O}_{A[j],t}$ (resp. $\mathcal{P}_{\phi_\sharp A[j],t}$) is $\Sigma$-cofibrant for all $t \geqslant 0$ plus the fact that

$$\begin{bmatrix} c \\ d \end{bmatrix} \longmapsto \operatorname*{colim}_{\Lambda}(\mathcal{O}_A; j)(\tau)$$

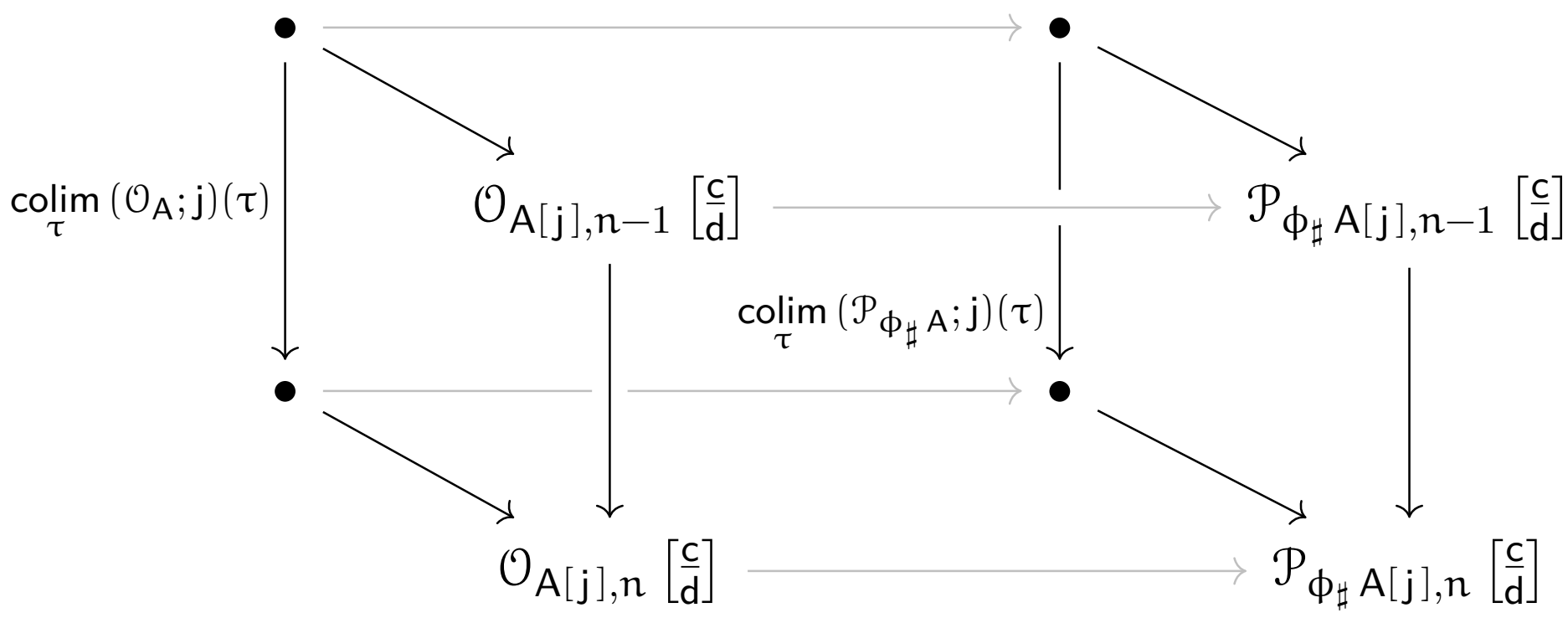

(resp. $(\mathcal{P}_{\phi_\sharp A}; j)$) is a core $\Sigma_{\underline{c}}$-cofibration (see the proof of Lemma 3.36 applying additionally that $j$ is a core cofibration). Thus, we should just check that the gray horizontal maps in the cube (except the lower one in the front face) are equivalences. This follows by induction and because $(\mathcal{O}_A; j)(\tau) \to (\mathcal{P}_{\phi_\sharp A}; j)(\tau)$ is an equivalence between $\mathrm{Aut}(\tau)^{\mathrm{op}}$-cofibrant objects since $\mathcal{O}_A \to \mathcal{P}_{\phi_\sharp A}$ is so and $j$ is a core cofibration.

In the $\Sigma$-cofibrant case this is enough, since the previous ladder starts with $\Sigma$-cofibrant objects and the horizontal maps are $\Sigma$-cofibrations. For the case where one of the operads might be well-pointed, use Lemma A.17. $\square$

The invariance properties of $(\mathcal{O}, A) \mapsto \mathcal{O}_A$ for $\Sigma$-cofibrant operads and proj-cofibrant algebras can be combined to obtain a last result in this respect.

**Corollary 3.42.** *$\mathcal{O}_{Q(\phi^* B)} \simeq \mathcal{P}_B$ when $B$ is a proj-cofibrant $\mathcal{P}$-algebra, $\mathcal{O} \to \mathcal{P}$ in $\mathcal{O}\mathrm{pd}_O(\mathcal{V})$ is an equivalence of operads and either:*

- *$\mathcal{O}$, $\mathcal{P}$ are $\Sigma$-cofibrant operads, or*
- *$\mathcal{O}, \mathcal{P} \in \{\Sigma$-cofibrant or well-pointed operads$\}$ and the $\mathbb{I}$-strong unit axiom holds.*

*Proof.* By the invariance of $B \mapsto \mathcal{P}_B$ on proj-cofibrant algebras, we can assume without loss of generality that $B$ is proj-bifibrant. Then, the (derived) counit $\phi_\sharp\phi^* B \to B$ is an equivalence ($\phi$ induces a Quillen equivalence between algebras). Since $\phi^*$ preserves equivalences, one may form the following commutative diagram of equivalences in $\mathtt{Alg}_{\mathcal{O}}$

$$\begin{array}{ccccc} \phi^* B & \longrightarrow & \phi^*\phi_\sharp\phi^* B & \xrightarrow{\sim} & \phi^* B \\ \wr\uparrow & & \uparrow & & \\ Q\,\phi^* B & \xrightarrow[\sim]{} & \phi^*\phi_\sharp\, Q\,\phi^* B & & \end{array}$$

by choosing a cofibrant replacement $Q\,\phi^* B \to \phi^* B$. Note that the decorated arrows are the ones that we already know are equivalences, the rest are so by 2-out of-3. Now apply $(i)$ in Propositions 3.38 and 3.40 to get the chain of equivalences of operads

$$\mathcal{O}_{Q\,\phi^* B} \xrightarrow{\sim} \mathcal{P}_{\phi_\sharp\, Q\,\phi^* B} \xrightarrow{\sim} \mathcal{P}_B\,.$$

□

*Remark* 3.43. There is no dual version for cofibrant operads and fgt-cofibrant algebras because one cannot ensure that $\mathbb{L}\phi_\sharp\, A \simeq \phi_\sharp\, A$ if $A$ is not proj-cofibrant. In other words, $\phi_\sharp$ does not preserve equivalences between fgt-cofibrant algebras in general.

In summary, the functor $(\mathcal{O}, A) \mapsto \mathcal{O}_A$ satisfies:

- **Cofibrancy**:
  - $(i)$ If the monoidal unit $\mathbb{I} \in \mathcal{V}$ is cofibrant, $\mathcal{O}_A$ is $\Sigma$-cofibrant if one of the following conditions holds:
    - $\mathcal{O}$ is $\Sigma$-cofibrant and $A$ is proj-cofibrant,
    - $\mathcal{O}$ is cofibrant and $A$ is fgt-cofibrant
  - $(ii)$ In general,
    - $\mathcal{O}_A$ is $\Sigma$-cofibrant (resp. well-pointed) if $\mathcal{O}$ is $\Sigma$-cofibrant (resp. well-pointed) and $A$ is proj-cofibrant,
    - $\mathcal{O}_A$ is well-pointed if $\mathcal{O}$ is cofibrant and $A$ is fgt-cofibrant
- **Preservation of equivalences**:

(i) $\mathcal{O}_A \to \mathcal{O}_B$ is an equivalence provided $A \to B$ is so and one of the following conditions holds:

- $\mathcal{O}$ is $\Sigma$-cofibrant or well-pointed and $A, B$ are proj-cofibrant,
- $\mathcal{O}$ is cofibrant and $A, B$ are fgt-cofibrant

(ii) $\mathcal{O}_A \to \mathcal{P}_{A'}$, associated to a map of operadic algebras $(\mathcal{O}, A) \to (\mathcal{P}, A')$, is an equivalence provided $\phi\colon \mathcal{O} \to \mathcal{P}$ is so and one of the following conditions holds:

- $\mathcal{O}, \mathcal{P}$ are $\Sigma$-cofibrant, $A$ is proj-cofibrant and the map of $\mathcal{O}$-algebras is the unit $A \to \phi^*\phi_\sharp A$,
- $\mathcal{O}, \mathcal{P}$ are cofibrant, $A'$ is fgt-cofibrant and the map of $\mathcal{P}$-algebras is the counit $\phi_\sharp\phi^* A' \to A'$.

# 4 Applications

## 4.1 Admissibility and rectification

**(Strong semi-)admissibility.** Recall that an operad $\mathcal{O}$ always produces a sifted monad, and hence, one can try to apply [10, Theorem 11.3.2, Theorem 12.1.4] or [40, Theorem 2.2.2] to the adjunction $F_\mathcal{O}\colon \mathtt{V}^{\mathrm{O}} \rightleftarrows \mathtt{Alg}_\mathcal{O}(\mathtt{V}) : \mathrm{fgt}$ to get a transferred (semi)model structure on $\mathcal{O}$-algebras. The difficult point is always to check that any cobase change of a map in $F_\mathcal{O}(J)$ is a weak equivalence. There are mainly two ways to see that this is the case: (1) analyzing those pushouts, i.e. acyclic cell attachments

$$\begin{array}{ccc} \mathcal{O}\circ X & \longrightarrow & A \\ {\scriptstyle \mathcal{O}\circ j}\downarrow & \ulcorner & \downarrow \\ \mathcal{O}\circ Y & \longrightarrow & A[j] \end{array}, \qquad \text{(cell attach.)}$$

explicitly; or (2) directly ensuring that $F_\mathcal{O}(J)$-cell $\subseteq \mathtt{Eq}$.

We first discuss alternative (1), which has been followed by a long list of authors: Spitzweck [36], Batanin-Berger [2], Berger-Moerdijk [3], Fresse [10], Harper [14], Hinich [15, 17], Muro [25, 26, 29], White-Yau [40]... The following results are direct applications of the homotopical analysis of enveloping operads developed in §3. A comparison with the literature is given at the end of the subsection.

*Remark* 4.1. A simple way to justify the connection between enveloping operads and cell attachments, without prior knowledge of the filtrations discussed in §3, is as follows. Rewrite (cell attach.) as a reflexive coequalizer in the usual way

$$A \amalg F_\mathcal{O}(X) \amalg F_\mathcal{O}(Y) \rightrightarrows A \amalg F_\mathcal{O}(Y) \xrightarrow{\text{colim}} A[j] \ .$$

One observes that the universal property of $\mathcal{O}_A$ (recall that $\mathtt{Alg}_{\mathcal{O}_A}(\mathtt{V}) \simeq A\downarrow \mathtt{Alg}_\mathcal{O}(\mathtt{V})$) can be used to absorb the factors $A \amalg \star$ into a free algebra functor; more concretely, $A \amalg F_\mathcal{O}(Z) \cong F_{\mathcal{O}_A}(Z)$. Since the free algebra functor is a left adjoint, we have that $F_\mathcal{O}(X) \amalg F_\mathcal{O}(Y) \cong F_\mathcal{O}(X \amalg Y)$, and the previous reflexive coequalizer can be rewritten as

$$F_{\mathcal{O}_A}(X \amalg Y) \rightrightarrows F_{\mathcal{O}_A}(Y) \xrightarrow{\text{colim}} A[j] \ .$$

Without assuming any additional hypothesis on $\mathcal{V}$, all $\Sigma$-cofibrant operads have a decent homotopy theory of algebras over $\mathcal{V}$ (see Definition 4.19) by the following results.

**Proposition 4.2.** *Let $\mathcal{O}$ be a $\Sigma$-cofibrant or well-pointed operad. Then, the free-forgetful adjunction $\mathtt{V}^{\mathsf{O}} \rightleftarrows \mathtt{Alg}_{\mathcal{O}}(\mathtt{V})$ induces by left-transfer the projective left semimodel structure $\mathcal{A}\mathrm{lg}_{\mathcal{O}}(\mathcal{V})$. That is, $\Sigma$-cofibrant and well-pointed operads are semi-admissible.*

*Proof.* Follows from the acyclic part of Lemma 3.35. □

**Proposition 4.3.** *The forgetful functor* $\mathrm{fgt}\colon \mathcal{A}\mathrm{lg}_{\mathcal{O}}(\mathcal{V}) \to \mathcal{V}^{\mathsf{O}}$ *preserves core (acyclic) cofibrations if $\mathcal{O}$ is $\Sigma$-cofibrant or well-pointed. In other words, $\mathcal{O}$ is strongly semi-admissible.*

*Proof.* The fact that proj-cofibrant $\mathcal{O}$-algebras are fgt-cofibrant is a direct consequence of Proposition 3.34 (i). The general claim follows from Lemma 3.35 using the same strategy. □

*Example* 4.4. Let $\mathcal{Q} \in \mathtt{Opd}(\mathtt{Set})$ be a $\Sigma$-free operad in sets. Then, its image through $\mathtt{Opd}(\mathtt{Set}) \to \mathtt{Opd}(\mathtt{V})$, induced by the sm-functor $S \mapsto \coprod_S \mathbb{I}$, is strongly semi-admissible in $\mathcal{V}$ (it is well-pointed). This applies to a variety of examples. For instance, the $\mathtt{Set}$-operad encoding: (1) $\mathtt{Opd}_{\mathsf{O}}(\mathtt{V})$ (O-colored operads), (2) $\mathtt{cfProp}_{\mathsf{O}}(\mathtt{V})$ (O-colored props without operations with domain $\mathbb{0}$; constant free) or (3) $\mathtt{afProp}_{\mathsf{O}}(\mathtt{V})$ (O-colored props without operations with target $\mathbb{0}$; augmentation free) is $\Sigma$-free. This already produces a homotopical context to deal with these classes of props. The outcome is less powerful than [11, Theorem 4.9], but for some applications is enough to have just a semimodel structure in the sense of Definition A.3. For example, one can recover the *homotopy transfer theorem* (see [11, Theorem C]) with the structure just mentioned, avoiding a non-negligible amount of non-trivial work.

Surprisingly, the $\mathtt{Set}$-operad encoding all O-colored props is not $\Sigma$-free.

For cofibrant operads, we get a better result. Let us first recall the relevant enhanced version of semimodel structure one obtains (see [1, 10, 36]).

**Definition 4.5.** Let $\mathcal{M}$ be a model category and consider an adjunction $\mathrm{F}\colon \mathtt{M} \rightleftarrows \mathtt{N} : \mathrm{R}$, where $\mathtt{N}$ is a bicomplete category. The transferred structure on $\mathtt{N}$ along $\mathrm{F} \dashv \mathrm{R}$

$$\mathcal{N} = \left(\mathtt{N},\ \mathtt{Eq}_{\mathcal{N}} := \mathrm{R}^{-1}\mathtt{Eq}_{\mathcal{M}},\ \mathtt{Fib}_{\mathcal{N}} := \mathrm{R}^{-1}\mathtt{Fib}_{\mathcal{M}},\ \mathtt{Cof}_{\mathcal{N}} := {}^{\boxslash}\mathtt{AFib}_{\mathcal{N}}\right)$$

is a *left* R*-semimodel structure* if the pair $(\mathtt{ACof}_{\mathcal{N}} := \mathtt{Cof}_{\mathcal{N}} \cap \mathtt{Eq}_{\mathcal{N}}, \mathtt{Fib}_{\mathcal{N}})$ satisfies:

(b1) (Core lifting axiom) maps in $\mathtt{ACof}_{\mathcal{N}}$ with R-cofibrant domain[i] satisfy the left lifting property (llp) against $\mathtt{Fib}_{\mathcal{N}}$.

(b2) (Core factorization axiom) any map $\mathrm{f}\colon \mathrm{X} \to \mathrm{Y}$ in $\mathcal{N}$ with R-cofibrant domain can be factored as

$$\mathrm{f}\colon \mathrm{X} \xrightarrow{\mathtt{ACof}_{\mathcal{N}}} \bullet \xrightarrow{\mathtt{Fib}_{\mathcal{N}}} \mathrm{Y}.$$

*Remark* 4.6. Note that the previous definition is essentially Definition A.3, but replacing any appearance of the word "cofibrant" by "R-cofibrant". All the other axioms are automatically satisfied because $\mathcal{M}$ is already a model category.

[i] Recall that a map f in $\mathtt{N}$ is an R-cofibration if Rf is a cofibration.

**Proposition 4.7.** *Let $\mathcal{O}$ be a cofibrant operad. Then, $\mathcal{A}\mathrm{lg}_{\mathcal{O}}(\mathcal{V})$ is a left* fgt*-semimodel structure. Moreover,* fgt: $\mathcal{A}\mathrm{lg}_{\mathcal{O}}(\mathcal{V}) \to \mathcal{V}^{O}$ *preserves (acyclic) cofibrations with* fgt*-cofibrant domain.*

*Proof.* The acyclic part of Lemma 3.35 implies the semi-admissibility result, when combined with Proposition 3.34 (ii). The non-acyclic part of this lemma together with the same proposition imply the other claim. □

Assuming very little on $\mathcal{V}$, it is also possible to obtain full admissibility for cofibrant operads. Recall that $X \in \mathcal{V}$ is *pseudo-cofibrant* if $X \otimes \star : \mathcal{V} \to \mathcal{V}$ preserves cofibrations (see [26, Appendix A], but note that we do not impose the monoid axiom on $\mathcal{V}$).

**Lemma 4.8.** *Any pseudo-cofibrant object $X \in \mathcal{V}$ satisfies that $X \otimes \star \colon \mathcal{V}^G \to \mathcal{V}^G$ preserves $G$-cofibrations, for any discrete group $G$.*

*Proof.* Since the generating cofibrations in $\mathcal{V}^G$ are of the form $G \cdot i = \coprod_G i$, where $i$ is a generating cofibration in $\mathcal{V}$, the result just follows from:

$$X \otimes (\textstyle\coprod_G i) \cong \coprod_G X \otimes i.$$

□

**Proposition 4.9.** *If all objects in $\mathcal{V}$ are pseudo-cofibrant, all cofibrant $\mathcal{V}$-operads are admissible. Furthermore, in this case* fgt: $\mathcal{A}\mathrm{lg}_{\mathcal{O}}(\mathcal{V}) \to \mathcal{V}^{O}$ *preserves all (acyclic) cofibrations if $\mathcal{O}$ is cofibrant.*

*Proof.* Let $\mathcal{O}$ be a cofibrant operad. For both statements, it suffices to check if for any pushout in $\mathcal{A}\mathrm{lg}_{\mathcal{O}}(\mathcal{V})$

$$\begin{array}{ccc} \mathcal{O} \circ X & \longrightarrow & A \\ {\scriptstyle \mathcal{O}\circ j}\downarrow & \ulcorner & \downarrow{\scriptstyle g} \\ \mathcal{O} \circ Y & \longrightarrow & A[j] \end{array}$$

where $j \colon X \to Y$ is a core (acyclic) cofibration in $\mathcal{V}^{O}$ (always possible by Remark A.9), $g$ is a(n acyclic) fgt-cofibration. Due to Proposition 3.17, we must check that $\mathrm{colim}_{\tau}\,(\mathcal{O}_A; j)(\tau)$ is a(n acyclic) cofibration. This follows from Lemma 3.27 if $\tau \mapsto (\mathcal{O}_A; j)(\tau)$ is a(n acyclic) $\mathrm{Aut}(\tau)^{op}$-cofibration. To show this, we assume that $\mathcal{O}$ is cellular and we run an induction over its cellular filtration.

First, consider the initial case $\mathcal{O} = \mathcal{I}_O$. Then, we know that the enveloping operad for $A$ is pretty simple:

$$(\mathcal{I}_O)_A \begin{bmatrix} \underline{c} \\ d \end{bmatrix} = \begin{cases} A(d) & \text{if } \underline{c} = \mathbb{0} \\ \mathbb{I} & \text{if } \underline{c} = d \\ \mathbb{0} & \text{else.} \end{cases}$$

To notice what is going on, suppose $j$ is concentrated on one color $o \in O$; the general case is analogous. Then, for $\tau$ in the relevant subgroupoid of $\mathrm{Tree}^{\cong}_{O,\not\downarrow}\left[\begin{smallmatrix} \underline{c} \\ d \end{smallmatrix}\right]$ (with $t \geqslant 1$ straight leaves)

$$\tau \longmapsto \big((\mathcal{I}_O)_A; j\big)(\tau) \cong (\mathcal{I}_O)_A \begin{bmatrix} \underline{c} \boxplus o^{\boxplus t} \\ d \end{bmatrix} \otimes j \begin{bmatrix} \mathbb{0} \\ o \end{bmatrix}^{\square t} = \begin{cases} j \left[\begin{smallmatrix} \mathbb{0} \\ o \end{smallmatrix}\right] & \text{if } \underline{c} = \mathbb{0},\ t = 1,\ o = d \\ \mathrm{id}_{\mathbb{0}} & \text{else.} \end{cases}$$

Hence, this map is a(n acyclic) $\mathrm{Aut}(\tau)^{\mathrm{op}}$-cofibration by trivial reasons.

For the induction step, we address what happens for a cell attachment $\mathcal{P} \rightarrowtail \mathcal{P}[k]$ of operads, where k is a core $\Sigma$-cofibration and $(\mathcal{I}_O)_A \to \mathcal{P}_A$ is a $\Sigma$-cofibration. Let us show that $\mathcal{P}_A \to \mathcal{P}[k]_A$ is also a $\Sigma$-cofibration. By Proposition 3.14, it suffices to see that

$$\begin{bmatrix} \underline{c} \\ d \end{bmatrix} \longmapsto \operatorname*{colim}_{\Lambda} (\mathcal{P}_A; k; A)(\Lambda) = \operatorname*{colim}_{\Lambda} \left( \bigotimes_{v \in V_{\mathrm{even}}(\Lambda)} \mathcal{P}_A \begin{bmatrix} s(v) \\ t(v) \end{bmatrix} \otimes \mathop{\Box}_{u \in V_{\mathrm{odd}}(\Lambda)} k \begin{bmatrix} s(u) \\ t(u) \end{bmatrix} \right)$$

is a $\Sigma_{\underline{c}}^{\mathrm{op}}$-cofibration (recall that $|V_{\mathrm{odd}}(\Lambda)| = t \geqslant 1$). Using Lemma 4.8, Remark A.14 and induction hypothesis, one observes that $\Lambda \mapsto (\mathcal{P}_A; k; A)(\Lambda)$ is a $G^{\mathrm{op}}$-cofibration for G the group of automorphims of $\Lambda$ which do not have to fix the leaves labeled by $\underline{c}$ (here we use that all objects are pseudo-cofibrant, in particular $\{A(o)\}_o$). See the proof of Lemma 3.36. Then, taking coinvariants over $\mathrm{Aut}(\Lambda)^{\mathrm{op}}$, i.e. automorphisms of $\Lambda$ which do fix the leaves labeled by $\underline{c}$, we obtain a $\Sigma_{\underline{c}}^{\mathrm{op}}$-cofibration due to Lemma 3.27.

Therefore, we obtain a $\Sigma$-cofibration $(\mathcal{I}_O)_A \to \mathcal{O}_A$ by induction on the cellular filtration of $\mathcal{O}$. Combining this fact with $\tau \mapsto \big((\mathcal{I}_O)_A; j\big)(\tau)$ being a(n acyclic) $\mathrm{Aut}(\tau)^{\mathrm{op}}$-cofibration, we find that so is $\tau \mapsto (\mathcal{O}_A; j)(\tau)$, finishing the proof. □

*Remark* 4.10. Proposition 4.7 was first claimed in Spitzweck's PhD thesis, [36, Theorem 4.3]. To the best of our understanding, the proof of [36, Theorem 4.3] appear to be incorrect. It makes crucial use of [36, Proposition 4.5], which seems to be erroneous. The problem is that some relations coming from the action of $\mathcal{O}$ on A are not taken into account in the map whose pushout is $B_{(i,j,k)-1} \to B_{(i,j,k)}$. To illustrate this point, under the notation of loc.cit., consider the unique map $g\colon \mathbb{0} \to \mathbb{I}$ from the initial object to the monoidal unit. Then, $B = A \amalg F_{\mathcal{O}}(\mathbb{I}) \cong \mathcal{O}_A \circ \mathbb{I}$ and $B' \cong \mathcal{O}'_A \circ \mathbb{I}$ by [14, Proposition 7.6]. From Proposition 3.14, we have a filtration on the map $\mathcal{O}_A \to \mathcal{O}'_A$ whose stages involve the enveloping operad $\mathcal{O}_A$. Applying $\star \circ \mathbb{I}$, we get a filtration on the map

$$B \cong \bigoplus_{n \geqslant 0} \big(\mathcal{O}_A(n)\big)_{\Sigma_n} \longrightarrow \bigoplus_{n \geqslant 0} \big(\mathcal{O}'_A(n)\big)_{\Sigma_n} \cong B'.$$

On the other hand, [36, Proposition 4.5] claims that there is a filtration on $B \to B'$ whose stages do not quotient by any relation involving the $\mathcal{O}$-action on A. This cannot work in general due to the previous explicit description of B and $B'$.

For a more concrete statement about the missing relations in [36, Proposition 4.5], it is written there:

"[...] $B_{(i,j,k)}$ *is a pushout of* $B_{(i,j,k)-1}$ *by the quotient of the map*

$$\coprod_T \Big( \bigotimes_v \mathcal{O}(\mathrm{val}(v)) \Big) \otimes A^{\otimes(a(T))} \otimes e^{\Box(U_{\mathrm{old}}(T))} \,\Box\, g^{\Box(m(T))} \,\Box \mathop{\Box}_{v'} f(\mathrm{val}(v')),$$

[...], *with respect to the equivalence relation which identifies for every isomorphism of directed graphs* $\varphi\colon T \to T'$ [...] *as in Proposition 3.4.3 and on the other parts by the identification of the indexing sets via* $\varphi$."

The alluded quotient from Proposition 3.4.3 only involves symmetric group actions.

A similar issue affects [36, Proposition 4.8]. See [29], in particular Example 1.1 in its arxiv version, for a counterexample even in the non-symmetric case. However, this other filtration has been corrected by other authors; e.g. [10, 14, 40].

*Remark* 4.11. In Proposition 4.9, one cannot replace the pseudo-cofibrancy hypothesis by the monoid axiom on $\mathcal{V}$. The reason is that the monoid axiom forces

$$\mathcal{O}_A \begin{bmatrix} s(\tau) \\ t(\tau) \end{bmatrix} \otimes \mathop{\square}_{\ell \in L_{\downarrow}(\tau)} j \begin{bmatrix} \mathbb{0} \\ \ell \end{bmatrix}$$

to be an equivalence when j is a core acyclic cofibration in $\mathcal{V}^{\mathrm{O}}$, but its quotient

$$\left( \mathcal{O}_A \begin{bmatrix} s(\tau) \\ t(\tau) \end{bmatrix} \otimes \mathop{\square}_{\ell \in L_{\downarrow}(\tau)} j \begin{bmatrix} \mathbb{0} \\ \ell \end{bmatrix} \right)_{\mathrm{Aut}(\tau)^{\mathrm{op}}}$$

may not be so (we lose control on $\mathcal{O}_A$ for a general $A \in \mathcal{A}lg_{\mathcal{O}}(\mathcal{V})$). However, if quotients by finite group actions are homotopical in $\mathcal{V}$, the previous quotient is still an equivalence and hence, we would obtain full admissibility of cofibrant operads if $\mathcal{V}$ satisfies the monoid axiom.

This fact contrasts strongly with the non-symmetric setting (see [25, 29]). However, despite of being different than a classical axiom ([35]), the pseudo-cofibrant condition on objects is potentially easy to check due to [26, Lemma A.3].

**Admissibility via path objects.** Let us collect a very simple criterion to check admissibility of operads through alternative (2) (see the first paragraph of this subsection), just for completeness. This is based on [3, 11, 35].

Assume that $\mathtt{Alg}_{\mathcal{O}}(\mathtt{V})$ admits the following structure, usually referred as *admitting path objects*:

- there is a fibrant replacement functor $A \xrightarrow{\sim} R(A)$ in $\mathtt{Alg}_{\mathcal{O}}(\mathtt{V})$, and
- every fibrant $\mathcal{O}$-algebra B admits a path object $B \xrightarrow{\sim} \mathrm{Path}(B) \twoheadrightarrow B \times B$.

Note that the previous structure only refers to equivalences and fibrations, which are detected by fgt: $\mathtt{Alg}_{\mathcal{O}}(\mathtt{V}) \to \mathtt{V}^{\mathrm{O}}$.

**Lemma 4.12.** *If* $\mathtt{Alg}_{\mathcal{O}}(\mathtt{V})$ *admits path objects, then* $F_{\mathcal{O}}(J)\text{-cell} \subseteq \mathtt{Eq}$.

*Proof.* See [35, Lemma 2.3]. Take any $j\colon A \to B$ in $F_{\mathcal{O}}(J)$-cell. We only know that j has left lifting property against fibrations in $\mathtt{Alg}_{\mathcal{O}}(\mathtt{V})$, but this suffices to show that j is an equivalence by the following argument due to Quillen. Choose a lifting for

$$\begin{array}{ccc} A & \xrightarrow{r} & R(A) \\ {\scriptstyle j}\downarrow & \overset{\rho}{\nearrow} & \\ B & & \end{array} .$$

Then, a solution for the lifting problem

$$\begin{array}{ccccccc} A & \xrightarrow{j} & B & \xrightarrow{\sim} & R(B) & \xrightarrow{\sim} & \mathrm{Path}(R(B)) \\ {\scriptstyle j}\downarrow & & & & & \nearrow & \downarrow \\ B & & & \xrightarrow{(r,\, R(j)\cdot\rho)} & & & R(B) \times R(B) \end{array}$$

yields a right homotopy $r \sim R(j) \cdot \rho$. Combined with the identity $\rho \cdot j = r$ and the fact that $r$ is an equivalence, we deduce that $\rho$ and $j$ are equivalences. □

The problem consists of finding hypotheses on $\mathcal{V}$ so that categories of operadic algebras admit path objects. Berger-Moerdijk [3] and Fresse [11] studied this problem giving the following result:

**Lemma 4.13.** *If $\mathcal{V}$ admits $\otimes$-path objects, any category $\mathtt{Alg}_{\mathcal{O}}(\mathtt{V})$ admits path objects. In particular, $\mathcal{V}$ admits $\otimes$-path objects if it has a lax monoidal fibrant replacement functor and it contains a cocommutative comonoidal interval object (see [3, §3]).*

Let us briefly recall the notion of *admitting $\otimes$-path objects* (see [11, §5]) to illustrate the simplicity of Lemma 4.13 and some examples.

We say that $\mathcal{V}$ admits $\otimes$-path objects if we can find the following structure:

- a lax monoidal fibrant replacement functor $R\colon \mathcal{V} \to \mathcal{V}$ whose monoidal structure makes the following triangles commute

$$\begin{array}{ccc} & X \otimes Y & \\ {}^{r\otimes r}\swarrow & & \searrow^{r} \\ R(X) \otimes R(Y) & \longrightarrow & R(X \otimes Y) \end{array} \quad ;$$

- a lax monoidal path functor $(\mathrm{Path}, \mathrm{id}_{\mathcal{V}} \to \mathrm{Path} \rightrightarrows \mathrm{id}_{\mathcal{V}})$ which yields path objects when evaluated on fibrant objects and whose monoidal structure makes the following diagrams commute

$$\begin{array}{ccc} & X \otimes Y & \\ \swarrow & & \searrow \\ \mathrm{Path}(X) \otimes \mathrm{Path}(Y) & \longrightarrow & \mathrm{Path}(X \otimes Y) \\ \rightrightarrows & & \rightrightarrows \\ & X \otimes Y & \end{array} \quad .$$

These conditions imply that the fibrant replacement functor and the path functor lift to categories of operadic algebras by simply applying the functors objectwise.

**Corollary 4.14.** *If $\mathcal{V}$ is combinatorial and admits $\otimes$-path objects, any operad in $\mathcal{V}$ is admissible, i.e. $\mathtt{Alg}_{\mathcal{O}}(\mathtt{V})$ admits the projective model structure for any $\mathcal{O}$.*

*Examples* 4.15. Simplicial sets, simplicial modules over an arbitrary ring, symmetric spectra and chain complexes over a commutative ring $\Bbbk$ such that $\mathbb{Q} \subseteq \Bbbk$ have combinatorial model structures and admit $\otimes$-path objects. The Quillen model structure on (cgwh) topological spaces also admits $\otimes$-path objects and sufficient smallness conditions (ensuring that Quillen's small object argument works) so that all operads are admissible. See [11, §5].

**Rectification and DK-equivalences.** We now address when an "equivalence" of operads induces a Quillen equivalence between algebras.

**Proposition 4.16.** *Let $\phi\colon \mathcal{O} \to \mathcal{P}$ be a map of $\Sigma$-cofibrant $\mathcal{V}$-operads. Then, the Quillen adjunction*

$$\phi_\sharp\colon \mathcal{A}\mathsf{lg}_{\mathcal{O}}(\mathcal{V}) \rightleftarrows \mathcal{A}\mathsf{lg}_{\mathcal{P}}(\mathcal{V}) : \phi^*$$

*is a Quillen equivalence provided that $\phi\colon \mathcal{O} \to \mathcal{P}$ is a DK-equivalence. If moreover $\mathcal{V}$ satisfies the $\mathbb{I}$-strong unit axiom, the operad $\mathcal{O}$ and/or $\mathcal{P}$ can be well-pointed instead of $\Sigma$-cofibrant for the conclusion to hold.*

*Proof.* First, note that any map of operads $\mathcal{O} \to \mathcal{O}'$ can be factored as a map in $\mathtt{Opd}_{\mathsf{col}(\mathcal{O})}(\mathtt{V})$ followed by a fully-faithful map of operads. This fact reduces the proof to two cases: $(a)$ $\phi$ is the identity on colors; and $(b)$ $\phi$ is fully-faithful on the nose.

Let us start with $(a)$. Since $\mathcal{A}\mathsf{lg}_{\mathcal{O}}(\mathcal{V}) \leftarrow \mathcal{A}\mathsf{lg}_{\mathcal{P}}(\mathcal{V}) : \phi^*$ creates equivalences, $\phi_\sharp \dashv \phi*$ yields a Quillen equivalence iff the unit transformation

$$\mathsf{A} \longrightarrow \phi^*\phi_\sharp\,\mathsf{A} \cong \mathcal{P} \underset{\mathcal{O}}{\circ} \mathsf{A}$$

is an equivalence for any cellular proj-cofibrant $\mathsf{A} \in \mathcal{A}\mathsf{lg}_{\mathcal{O}}(\mathcal{V})$. This fact follows directly from Proposition 3.40 $(i)$.

Let us now discuss $(b)$. Recall that we are assuming that $\mathcal{O} \to \mathcal{P}$ is fully-faithful on the nose and that $\pi_0(\overline{\mathcal{O}}) \to \pi_0(\overline{\mathcal{P}})$ is essentially surjective. Hence, we are reduced to show that the derived unit and counit of $\phi_\sharp \dashv \phi^*$ are equivalences under these assumptions. The underived unit is an isomorphism since $\mathcal{O} \to \mathcal{P}$ is fully-faithful on the nose, so it is enough to analyze the derived counit.

Just observe that there is a commutative diagram

$$\begin{array}{ccc}
\mathsf{Ho}\,\mathcal{A}\mathsf{lg}_{\mathcal{O}}(\mathcal{V}) & \xleftarrow{\quad\phi^*\quad} & \mathsf{Ho}\,\mathcal{A}\mathsf{lg}_{\mathcal{P}}(\mathcal{V}) \\
\downarrow & & \downarrow \\
\mathsf{Ho}\,\left[\overline{\mathcal{O}}, \mathcal{V}\right] & \xleftarrow{\quad\phi^*\quad} & \mathsf{Ho}\,\left[\overline{\mathcal{P}}, \mathcal{V}\right] \\
\downarrow & & \downarrow \\
\left[\pi_0(\overline{\mathcal{O}}), \pi_0(\mathcal{V})\right] & \xleftarrow[\pi_0(\phi)^*]{} & \left[\pi_0(\overline{\mathcal{P}}), \pi_0(\mathcal{V})\right]
\end{array}$$

and that equivalences in $\mathcal{A}\mathsf{lg}_{\mathcal{O}}(\mathcal{V})$ and $\mathcal{A}\mathsf{lg}_{\mathcal{P}}(\mathcal{V})$ are detected by the vertical functors, almost by definition. Then, whether or not $\mathbb{L}\phi_\sharp\phi^*\,\mathsf{B} \to \mathsf{B}$ is an equivalence can be detected by going down in the diagram. Since the lower horizontal functor yields an equivalence of categories by hypothesis, one deduces that $\mathbb{L}\phi_\sharp\phi^*\,\mathsf{B} \to \mathsf{B}$ is an equivalence. □

*Remark* 4.17. An elegant way to rephrase part of Proposition 4.16 is to assert that the functor

$$\left(\mathtt{Opd}(\mathtt{V})_{\mathsf{semiadm}}, \mathtt{Eq}_{\mathsf{DK}}\right) \longrightarrow \left(\mathtt{semiMODEL}^{\mathtt{R}}, \mathtt{Eq}_{\mathsf{Q}}\right), \quad \mathcal{O} \mapsto \mathcal{A}\mathsf{lg}_{\mathcal{O}}(\mathcal{V})$$

is homotopical on $\Sigma$-cofibrant operads, where the right hand side denotes the very large relative category of left semimodel categories with right Quillen functors between them.

**Corollary 4.18.** *Assume $\mathcal{V}$ satisfies the $\mathbb{I}$-strong unit axiom and let $\phi\colon \mathcal{O} \to \mathcal{P}$ be a DK-equivalence of cofibrant $\mathcal{V}$-operads. Then, the Quillen adjunction*

$$\phi_\sharp\colon \mathcal{A}lg_{\mathcal{O}}(\mathcal{V}) \rightleftarrows \mathcal{A}lg_{\mathcal{P}}(\mathcal{V}) : \phi^*$$

*is a Quillen equivalence.*

**Summary and relations with the literature.** A nice way to encapsulate the previous results (Propositions 4.2, 4.3, 4.9 and 4.16) is via the following working definition (based on [15, 17, 19]).

**Definition 4.19.** Let $\mathfrak{O} \subseteq \mathtt{Opd}(\mathtt{V})$ be a class of colored operads. Then, we say that $\mathcal{V}$ admits a *good homotopy theory of algebras over operads in* $\mathfrak{O}$ if:

- every operad $\mathcal{O} \in \mathfrak{O}$ is strongly admissible, and
- any DK-equivalence $\varphi\colon \mathcal{O} \to \mathcal{P}$ between operads in $\mathfrak{O}$ induces a Quillen equivalence

$$\varphi_\sharp\colon \mathcal{A}lg_{\mathcal{O}}(\mathcal{V}) \simeq_Q \mathcal{A}lg_{\mathcal{P}}(\mathcal{V}) : \varphi^*.$$

If in the first condition elements in $\mathfrak{O}$ are only strongly semi-admissible, we say that $\mathcal{V}$ admits a *decent homotopy theory of algebras over operads in* $\mathfrak{O}$.

*Examples* 4.20. Let $\Bbbk$ be a field with $\mathrm{char}(\Bbbk) = 0$ and R be a general commutative ring. Then, the projective model structure $\mathcal{C}h(\Bbbk)$ (resp. $\mathcal{C}h(R)$) admits a good homotopy theory of algebras over all operads (resp. all $\Sigma$-split and locally cofibrant operads); see [15, Theorem 2.6.1]. The Kan-Quillen model structure on simplicial sets $\mathcal{S}pc$ admits a good homotopy theory of algebras over all $\Sigma$-cofibrant operads.

From the results mentioned before, we obtain:

**Corollary 4.21.** *Let $\mathcal{V}$ be a cofibrantly generated closed sm-model category. Then,*

- *$\mathcal{V}$ admits a decent homotopy theory of algebras over $\Sigma$-cofibrant operads.*
- *If $\mathcal{V}$ satisfies the $\mathbb{I}$-strong unit axiom, $\mathcal{V}$ admits a decent homotopy theory of algebras over well-pointed, as well as $\Sigma$-cofibrant, operads.*
- *If all objects in $\mathcal{V}$ are pseudo-cofibrant and $\mathcal{V}$ satisfies the $\mathbb{I}$-strong unit axiom, $\mathcal{V}$ admits a good homotopy theory of algebras over cofibrant operads.*

*Remark* 4.22. In all the results in this subsection, one essential ingredient is that for any pushout (cell attach.), the map $A \to A[j]$ admits a filtration

$$A = A[j]_0 \longrightarrow \cdots \longrightarrow A[j]_{t-1} \xrightarrow{g_t} A[j]_t \longrightarrow \cdots$$

when seen in $\mathcal{V}^O$ (it is just the 0-arity component of the filtration for $\mathcal{O}_A \to \mathcal{O}_{A[j]}$), where $g_t$ is a pushout of a map $\widetilde{g}_t$ constructed out of $\mathcal{O}_A$ and j. More precisely, and assuming j is concentrated in a single color $c \in O$ for simplicity, the map $\widetilde{g}_t$ is given by

$$\widetilde{g}_t(d) = \mathcal{O}_A \begin{bmatrix} c^{\boxplus t} \\ d \end{bmatrix} \underset{\Sigma_t}{\otimes} j \begin{bmatrix} \mathbb{0} \\ c \end{bmatrix}^{\square t}.$$

Therefore, to obtain more general admissibility and rectification results for operadic algebras, one can look for axioms over $\mathcal{V}$ ensuring good homotopical properties of expressions of this kind, i.e. $\mathcal{X}\left[\begin{smallmatrix}\underline{c}\\ d\end{smallmatrix}\right] \otimes_{\mathrm{Aut}(\underline{c})} \square_i f \left[\begin{smallmatrix}\mathbb{0}\\ c_i\end{smallmatrix}\right]$ or its relative version

$$\left( x \begin{bmatrix}\underline{c}\\ d\end{bmatrix} \square \bigsquare_i f \begin{bmatrix}\mathbb{0}\\ c_i\end{bmatrix} \right)_{\mathrm{Aut}(\underline{c})^{op}}$$

for a morphism $x \colon \mathcal{X} \to \mathcal{Y}$ in $\Sigma\mathcal{C}oll_O(\mathcal{V})$ and a morphism f in $\mathcal{V}^O$. This is the strategy followed by Pavlov-Scholbach in [30] and White-Yau in [40] and the *raison d'être* for axioms such as *symmetric h-monoidality* or *symmetroidality* in [31].

**Related work:** Let us comment some connections with the existing literature.

Probably, one of the first references about this kind of problematic (admissibility and rectification of operadic algebras in general model categories) is Spitzweck's thesis [36]. While inspiring, this document contains some technical gaps (see Remark 4.10). Other authors, like Berger-Moerdijk [3, 4] or Pavlov-Scholbach [30], worked on cleaning and expanding Spitzweck's original ideas. However, in these references, it is assumed that operads are well-pointed instead of just Σ-cofibrant for some of their results. In particular, for strong-admissibility and rectification theorems (e.g. [30, Theorem 1.1] and [4, Theorem 4.1][ii]). We want to point out that our proofs work in both cases (Σ-cofibrant and well-pointed) and the reason for this is that the filtrations developed in §3 are designed to avoid the unit inclusion; ultimately due to Proposition A.11. Also, the axioms that we impose on $\mathcal{V}$ are minimal. Because of that, our work complements [30], where the authors impose strong assumptions on $\mathcal{V}$ to work with all operads.

On the other hand, Fresse [10, Theorem 12.3.A] and White-Yau [40], probably inspired by Spitzweck's thesis as well, showed strong semi-admissibility for Σ-cofibrant operads (as done in Propositions 4.2 and 4.3). The well-pointed version is new.

Notice that our rectification result, Proposition 4.16, deals with general DK-equivalences, while the rectification results in [4, 30] (see also [10, Theorem 12.3.4]) apply to morphisms of operads which are the identity on colors.

Regarding cofibrant operads, Proposition 4.7 appeared as [36, Theorem 4.3]. The proof there was based on Proposition 4.5 in loc.cit., which has not being replicated in the literature and seems to contain a major gap (see Remark 4.10). Proposition 4.9, which is also inspired by [36, Theorem 4.3], is new and generalizes the main theorem of [16] (see also [17, Proposition 2.3.2]) to a quite general $\mathcal{V}$. Notice that the monoid axiom is not enough to get full admissibility of cofibrant operads by Remark 4.11.

Admissibility and rectification of uncolored non-symmetric operads in general model categories were addressed in [25, 26, 29].

## 4.2 Change of homotopy cosmos

Given a colax-lax symmetric monoidal adjunction $F \colon V \rightleftarrows V' : R$, [7, Theorem 4.5.6] provides liftings to monadic categories over V and V′ of $F \dashv R$. In particular, we have

$$F^{opd} \colon \mathrm{Opd}_O(V) \rightleftarrows \mathrm{Opd}_O(V') : R$$

[ii]Notice that several hypotheses are assumed on $\mathcal{V}$, while we virtually assume no condition.

and for any map of O-operads $\mathcal{O} \to R\,\mathcal{P}$, where $\mathcal{O} \in \mathtt{Opd}_O(\mathtt{V})$, $\mathcal{P} \in \mathtt{Opd}_O(\mathtt{V}')$,

$$F^{alg}\colon \mathtt{Alg}_{\mathcal{O}}(\mathtt{V}) \rightleftarrows \mathtt{Alg}_{\mathcal{P}}(\mathtt{V}') :R \ .$$

Note that the right adjoint R automatically lifts to operads and algebras since it is a lax sm-functor. Of course, we are interested in adding a homotopical flavor to this.

One of the main results in [41] shows that, under some hypotheses, a Quillen equivalence $F\colon \mathcal{V} \rightleftarrows \mathcal{V}' :R$ lifts to Quillen equivalences at the level of operadic algebras. Recall the notion of weak sm-Quillen adjunction/equivalence ([41, Definition 2.7]), which is a useful homotopical enhancement of a colax-lax sm-adjunction.

The following theorem is a slight generalization of this, which drops one hypothesis on the generating cofibrations of $\mathcal{V}$. See Remark 4.24 for another slight generalization for well-pointed operads.

**Theorem 4.23.** *[4.3.1 in [41]] Let $F\colon \mathcal{V} \rightleftarrows \mathcal{V}' :R$ be a weak sm-Quillen equivalence between homotopy cosmoi and $\mathcal{O} \in \mathcal{O}pd(\mathcal{V})$, $\mathcal{P} \in \mathcal{O}pd(\mathcal{V}')$ be $\Sigma$-cofibrant operads. Then, the Quillen pair associated to a map $\mathcal{O} \to R\,\mathcal{P}$ in $\mathcal{O}pd_O(\mathcal{V})$,*

$$F^{alg}\colon \mathcal{A}lg_{\mathcal{O}}(\mathcal{V}) \rightleftarrows \mathcal{A}lg_{\mathcal{P}}(\mathcal{V}') :R \ ,$$

*is a Quillen equivalence if the pointwise adjoint map $F\left(\mathcal{O}\left[\begin{smallmatrix}c\\d\end{smallmatrix}\right]\right) \to \mathcal{P}\left[\begin{smallmatrix}c\\d\end{smallmatrix}\right]$ is an equivalence for any O-corolla $\left[\begin{smallmatrix}c\\d\end{smallmatrix}\right]$.*

*Proof.* Let us provide a sketch of the proof resembling [41, Theorem 4.3.1], to also stress why it is possible to drop hypothesis (3) in loc.cit.. First note that it suffices to show that, for any proj-cofibrant $A \in \mathcal{A}lg_{\mathcal{O}}(\mathcal{V})$, the canonical comparison map $F(A(c)) \to (F^{alg}A)(c)$ is an equivalence $\forall c \in O$, as explained in the proof of [41, Theorem 4.2.1] (here one uses that $F \dashv R$ is a Quillen equivalence). Such map is the arity zero component of a map relating enveloping operads

$$F\left(\mathcal{O}_A\left[\begin{smallmatrix}c\\d\end{smallmatrix}\right]\right) \longrightarrow \mathcal{P}_{F^{alg}A}\left[\begin{smallmatrix}c\\d\end{smallmatrix}\right].$$

Therefore, we are done if these morphisms are weak equivalences for any O-corolla. Without loss of generality, assume that A is cellular proj-cofibrant and note that the cellular filtration of A yields a map of transfinite composites

$$\begin{array}{ccccccccc}
F\left(\mathcal{O}\left[\begin{smallmatrix}c\\d\end{smallmatrix}\right]\right) = F\left(\mathcal{O}_{A_0}\left[\begin{smallmatrix}c\\d\end{smallmatrix}\right]\right) & \longrightarrow & \cdots & \longrightarrow & F\left(\mathcal{O}_{A_\alpha}\left[\begin{smallmatrix}c\\d\end{smallmatrix}\right]\right) & \longrightarrow & \cdots & \longrightarrow & F\left(\mathcal{O}_{A}\left[\begin{smallmatrix}c\\d\end{smallmatrix}\right]\right)\\
\downarrow & & & & \downarrow & & & & \downarrow\\
\mathcal{P}\left[\begin{smallmatrix}c\\d\end{smallmatrix}\right] = \mathcal{P}_{F^{alg}A_0}\left[\begin{smallmatrix}c\\d\end{smallmatrix}\right] & \longrightarrow & \cdots & \longrightarrow & \mathcal{P}_{F^{alg}A_\alpha}\left[\begin{smallmatrix}c\\d\end{smallmatrix}\right] & \longrightarrow & \cdots & \longrightarrow & \mathcal{P}_{F^{alg}A}\left[\begin{smallmatrix}c\\d\end{smallmatrix}\right]
\end{array}$$

such that: (a) the map that we want to analyze, the right most vertical map, is the colimit of the rest of the vertical maps; (b) the horizontal maps are $\Sigma$-cofibrations by Lemma 3.35 (since F and $F^{alg}$ are left Quillen); (c) the first two terms $F\left(\mathcal{O}\left[\begin{smallmatrix}c\\d\end{smallmatrix}\right]\right)$ and $\mathcal{P}\left[\begin{smallmatrix}c\\d\end{smallmatrix}\right]$ are cofibrant when $\mathcal{O}$ and $\mathcal{P}$ are $\Sigma$-cofibrant. Thus, we are reduced to show:

$$\left(\begin{array}{c} F\left(\mathcal{O}_{A_\alpha}\left[\begin{smallmatrix}c\\d\end{smallmatrix}\right]\right) \longrightarrow \mathcal{P}_{F^{alg}A_\alpha}\left[\begin{smallmatrix}c\\d\end{smallmatrix}\right] \\ \text{is an equivalence } \forall \left[\begin{smallmatrix}c\\d\end{smallmatrix}\right] \end{array}\right) \implies \left(\begin{array}{c} F\left(\mathcal{O}_{A_{\alpha+1}}\left[\begin{smallmatrix}c\\d\end{smallmatrix}\right]\right) \longrightarrow \mathcal{P}_{F^{alg}A_{\alpha+1}}\left[\begin{smallmatrix}c\\d\end{smallmatrix}\right] \\ \text{is an equivalence } \forall \left[\begin{smallmatrix}c\\d\end{smallmatrix}\right] \end{array}\right),$$

since the hypothesis in the statement corresponds to $\alpha = 0$. This induction step $(\alpha) \Rightarrow (\alpha+1)$ is nothing more than the analysis of a cell attachment

$$\begin{array}{ccc} \mathcal{O}\circ X & \longrightarrow & A' \\ {\scriptstyle \mathcal{O}\circ j}\downarrow & \ulcorner & \downarrow{\scriptstyle g} \\ \mathcal{O}\circ Y & \longrightarrow & A'[j] \end{array},$$

where $j$ is a cofibration in $\mathcal{V}^{O}$. White-Yau assume that generating cofibrations of $\mathcal{V}$ have cofibrant domains to ensure that the domain of $\mathcal{O}\circ j$ is cofibrant, but this can be always arranged in this situation due to Remark A.9 (i.e. without imposing such condition on $\mathcal{V}$). In other words, we can always assume that $j$ is a core cofibration.

Now, one runs another induction to deal with the filtration on enveloping operads associated to this cell attachment (see Proposition 3.17). The relevant step, which is the induction step $(t-1) \Rightarrow (t)$, sits into a cubical diagram whose black faces are pushouts (note that we used that $F$ is a left adjoint functor and that $F^{alg}(\mathcal{O}\circ j) \cong \mathcal{P}\circ F j$). By induction hypothesis, the black faces are homotopy pushouts if we check that

$$\operatorname*{colim}_{\tau}\,(\mathcal{O}_{A'};\, j)(\tau) \qquad \text{and} \qquad \operatorname*{colim}_{\tau}\,(\mathcal{P}_{F^{alg}A'};\, F j)(\tau)$$

are core cofibrations $(1)$. If that is the case, we would finish the proof by showing that the horizontal arrows in the square

$$\operatorname*{colim}_{\tau}\,(F\,\mathcal{O}_{A'};\, F j)(\tau) \Longrightarrow \operatorname*{colim}_{\tau}\,(\mathcal{P}_{F^{alg}A'};\, F j)(\tau),$$

sitting in the back of the cube, are equivalences $(2)$. For $(1)$, use the proof of Lemma 3.35. For $(2)$, note that $F\,\mathcal{O}_{A'} \to \mathcal{P}_{F^{alg}A'}$ is an equivalence in $\Sigma\mathcal{C}oll_{O}(\mathcal{V}')$ by induction hypothesis and hence, the map in $\mathcal{V}^{2}$

$$(F\,\mathcal{O}_{A'}; F j)(\tau) \Longrightarrow (\mathcal{P}_{F^{alg}A'}; F j)(\tau)$$

is an equivalence $\forall\tau$ since $F \dashv R$ is a weak sm-Quillen adjunction. Applying $\operatorname{colim}_{\tau}$ in the arrow category $\mathcal{V}^{2}$ preserves this equivalence by $\operatorname{Aut}(\tau)$-cofibrancy of the factors (see the proof of Lemma 3.35). □

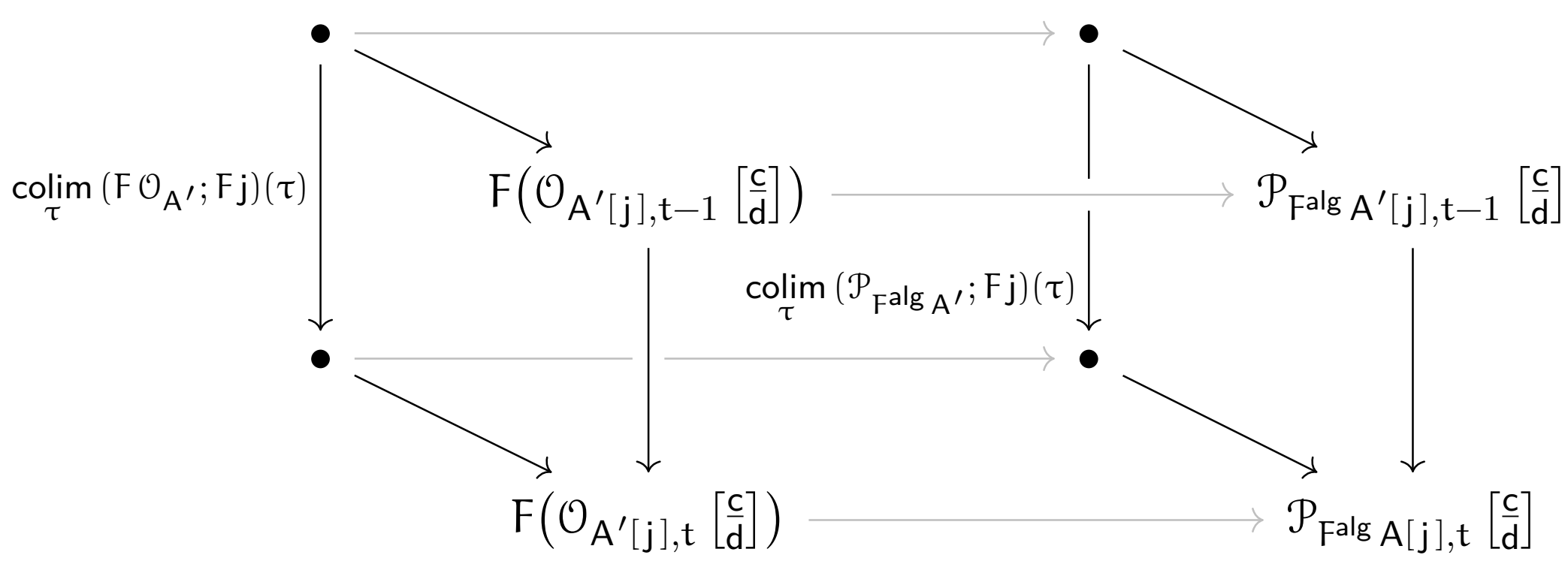

*Remark* 4.24. Using results from [26, Appendix A, B], it is possible to obtain a mild generalization of Theorem 4.23 where $\mathcal{O}$ and $\mathcal{P}$ can be considered well-pointed instead of $\Sigma$-cofibrant. When the monoidal units $\mathbb{I}_{\mathtt{V}}$ and $\mathbb{I}_{\mathtt{V}'}$ are cofibrant, the statement is obviously true. If this cofibrancy condition is not satisfied, the same conclusion can be obtained by adding the hypotheses:

- strong unit axiom on $\mathcal{V}$ and $\mathcal{V}'$ [26, Definition A.9],
- $\mathrm{F}$ satisfies pseudo-cofibrant axiom [26, Definition B.6 (1)],
- $\mathrm{F}$ satisfies $\mathbb{I}$-cofibrant axiom [26, Definition B.6 (2)].

Checking the condition: $\mathrm{F}\left(\mathcal{O}\left[\begin{smallmatrix}c\\d\end{smallmatrix}\right]\right) \to \mathcal{P}\left[\begin{smallmatrix}c\\d\end{smallmatrix}\right]$ is an equivalence $\forall\left[\begin{smallmatrix}c\\d\end{smallmatrix}\right]$, is not easy in practice. However, when the monoidal units $\mathbb{I}_{\mathtt{V}}$, $\mathbb{I}_{\mathtt{V}'}$ are cofibrant, there is a class of examples for which the condition is automatic. Assume that $\mathcal{O}$ and $\mathcal{P}$ are the images of a colored $\mathtt{Set}$-operad $\mathcal{Q} \in \mathtt{Opd}_{\mathrm{O}}(\mathtt{Set})$ via the functors

$$\mathtt{Opd}_{\mathrm{O}}(\mathtt{V}) \longleftarrow \mathtt{Opd}_{\mathrm{O}}(\mathtt{Set}) \longrightarrow \mathtt{Opd}_{\mathrm{O}}(\mathtt{V}')$$

induced by the sm-functor $\mathtt{Set} \to \mathtt{V}$, $\mathrm{S} \mapsto \coprod_{\mathrm{S}} \mathbb{I}$. Then, there is an obvious map $\mathcal{O} \to \mathrm{R}\,\mathcal{P}$ and it satisfies

$$\mathrm{F}\left(\mathcal{O}\left[\begin{smallmatrix}c\\d\end{smallmatrix}\right]\right) \cong \coprod_{\mathcal{Q}\left[\begin{smallmatrix}c\\d\end{smallmatrix}\right]} \mathrm{F}(\mathbb{I}_{\mathtt{V}}) \xrightarrow{\ \sim\ } \coprod_{\mathcal{Q}\left[\begin{smallmatrix}c\\d\end{smallmatrix}\right]} \mathbb{I}_{\mathtt{V}'} = \mathcal{P}\left[\begin{smallmatrix}c\\d\end{smallmatrix}\right]$$

by definition of weak sm-Quillen equivalence (and since taking coproducts of equivalences between cofibrant objects yield equivalences). Most of the examples provided in [41] are of this form.

The next result shows that the two most natural examples of maps of $\mathcal{V}$-operads that fit into Theorem 4.23, namely $\mathrm{id}\colon \mathrm{R}\,\mathcal{P} \to \mathrm{R}\,\mathcal{P}$ and $\eta\colon\ \mathcal{O} \to \mathrm{RF}^{\mathrm{opd}}\,\mathcal{P}$, satisfy the conclusion under natural assumptions.

**Proposition 4.25.** *Let* $\mathrm{F}\colon \mathcal{V} \rightleftarrows \mathcal{V}' : \mathrm{R}$ *be a weak sm-Quillen equivalence between homotopy cosmoi and* $\mathcal{O} \in \mathcal{O}\mathrm{pd}_{\mathrm{O}}(\mathcal{V})$, $\mathcal{P} \in \mathcal{O}\mathrm{pd}_{\mathrm{O}}(\mathcal{V}')$ *be two operads. Then,*

- *An equivalence of* $\mathrm{O}$*-operads* $\mathcal{O} \to \mathrm{R}\,\mathcal{P}$ *induces a Quillen equivalence*

  $$\mathrm{F}^{\mathrm{alg}}\colon \mathcal{A}\mathrm{lg}_{\mathcal{O}}(\mathcal{V}) \rightleftarrows \mathcal{A}\mathrm{lg}_{\mathcal{P}}(\mathcal{V}') : \mathrm{R}$$

  *if* $\mathcal{O}$ *is* $\Sigma$*-cofibrant and* $\mathcal{P}$ *is fibrant and* $\Sigma$*-cofibrant. In particular, the identity map* $\mathrm{id}\colon \mathrm{R}\,\mathcal{P} \to \mathrm{R}\,\mathcal{P}$ *satisfies the claim if* $\mathcal{P}$ *and* $\mathrm{R}\,\mathcal{P}$ *are fibrant and* $\Sigma$*-cofibrant.*[iii]

- *The unit map* $\eta\colon\ \mathcal{O} \to \mathrm{RF}^{\mathrm{opd}}\,\mathcal{O}$ *induces a Quillen equivalence*

  $$\mathrm{F}^{\mathrm{alg}}\colon \mathcal{A}\mathrm{lg}_{\mathcal{O}}(\mathcal{V}) \rightleftarrows \mathcal{A}\mathrm{lg}_{\mathrm{F}^{\mathrm{opd}}\,\mathcal{O}}(\mathcal{V}') : \mathrm{R}$$

  *if* $\mathcal{O}$ *is cofibrant as an operad and* $\mathbb{I}_{\mathtt{V}}$, $\mathbb{I}_{\mathtt{V}'}$ *are also cofibrant*[iv].

[iii] There is a similar statement for well-pointed operads that you can find out using Remark 4.24.
[iv] The hypothesis on the monoidal units can be relaxed; see Remark 4.24.

*Proof.* For the first item, just note that $F \dashv R$ being a Quillen equivalence implies that:

$$f \colon Fv \to v' \text{ is an equivalence} \qquad \Longleftrightarrow \qquad f^\flat \colon v \to Rv' \text{ is an equivalence},$$

for any $v \in \mathcal{V}$ cofibrant and $v' \in \mathcal{V}'$ fibrant. Apply this remark to the maps

$$\left\{ \mathcal{O}\begin{bmatrix} c \\ d \end{bmatrix} \to R\left(\mathcal{P}\begin{bmatrix} c \\ d \end{bmatrix}\right) : \quad \begin{bmatrix} c \\ d \end{bmatrix} \text{ is an } O\text{-corolla} \right\}.$$

The second item requires more work, since it is not just a direct consequence of Theorem 4.23. However, we are reduced to check that for any cofibrant operad $\mathcal{O}$, the pointwise comparison map $F\left(\mathcal{O}\begin{bmatrix} c \\ d \end{bmatrix}\right) \to \left(F^{\mathrm{opd}}\,\mathcal{O}\right)\begin{bmatrix} c \\ d \end{bmatrix}$ is an equivalence for any O-corolla. Observe that both $\mathcal{O}$ and $F^{\mathrm{opd}}\,\mathcal{O}$ are $\Sigma$-cofibrant by Corollary 3.37 and the hypothesis on the monoidal units. This justifies the application of Theorem 4.23 to this case.

Not surprisingly, we assume without loss of generality that $\mathcal{O}$ is cellular and we argue by induction over its cellular filtration, i.e. we analyze cell attachments:

$$\begin{array}{ccc} \mathcal{F}(\mathcal{X}) & \longrightarrow & \mathcal{O} \\ {\scriptstyle \mathcal{F}j}\big\downarrow & \ulcorner & \big\downarrow{\scriptstyle g} \\ \mathcal{F}(\mathcal{Y}) & \longrightarrow & \mathcal{O}[j] \end{array}, \tag{4.1}$$

where j is a core cofibration (always possible by Remark A.9). Note that the initial step in the induction over the cellular filtration is automatically satisfied since the initial operad $\mathcal{I}_O^{\mathcal{V}} \in \mathcal{O}\mathrm{pd}_O(\mathcal{V})$ is sent to the initial operad $\mathcal{I}_O^{\mathcal{V}'} \in \mathcal{O}\mathrm{pd}_O(\mathcal{V}')$ via $F^{\mathrm{opd}}$, and we assume: $F \dashv R$ is a weak sm-Quillen adjunction. For a cell attachment such as (4.1), we obtain a decomposition of $g\colon \mathcal{O}\begin{bmatrix} c \\ d \end{bmatrix} \to \mathcal{O}[j]\begin{bmatrix} c \\ d \end{bmatrix}$ as a $\omega$-transfinite composite

$$\mathcal{O}\begin{bmatrix} c \\ d \end{bmatrix} = \mathcal{O}[j]_0\begin{bmatrix} c \\ d \end{bmatrix} \longrightarrow \cdots \longrightarrow \mathcal{O}[j]_{t-1}\begin{bmatrix} c \\ d \end{bmatrix} \xrightarrow{g_t} \mathcal{O}[j]_t\begin{bmatrix} c \\ d \end{bmatrix} \longrightarrow \cdots,$$

by Proposition 3.14. Each $g_t$ is part of a $\Sigma$-cofibration by the same argument applied in Lemma 3.36, and induction hypothesis. Arguing similarly with $F^{\mathrm{opd}}\,\mathcal{O} \to F^{\mathrm{opd}}\,\mathcal{O}[Fj]$, we get a map of transfinite composites

$$\begin{array}{ccccccc} F\left(\mathcal{O}\begin{bmatrix} c \\ d \end{bmatrix}\right) = F\left(\mathcal{O}[j]_0\begin{bmatrix} c \\ d \end{bmatrix}\right) & \longrightarrow & \cdots & \xrightarrow{F g_t} & F\left(\mathcal{O}[j]_t\begin{bmatrix} c \\ d \end{bmatrix}\right) & \longrightarrow & \cdots \\ \big\downarrow & & & & \big\downarrow & & \\ \left(F^{\mathrm{opd}}\,\mathcal{O}\right)\begin{bmatrix} c \\ d \end{bmatrix} = \left(F^{\mathrm{opd}}\,\mathcal{O}[Fj]\right)_0\begin{bmatrix} c \\ d \end{bmatrix} & \longrightarrow & \cdots & \xrightarrow[\widetilde{g}_t]{} & \left(F^{\mathrm{opd}}\,\mathcal{O}[Fj]\right)_t\begin{bmatrix} c \\ d \end{bmatrix} & \longrightarrow & \cdots \end{array}$$

such that: (a) the map that we want to analyze is the colimit of the vertical maps in this diagram; (b) the horizontal maps are cofibrations; (c) the first two terms $F\left(\mathcal{O}\begin{bmatrix} c \\ d \end{bmatrix}\right)$ and $\left(F^{\mathrm{opd}}\,\mathcal{O}\right)\begin{bmatrix} c \\ d \end{bmatrix}$ are cofibrant. Thus, we are left to show:

$$\begin{pmatrix} F\left(\mathcal{O}[j]_{t-1}\begin{bmatrix} c \\ d \end{bmatrix}\right) \longrightarrow \left(F^{\mathrm{opd}}\,\mathcal{O}[Fj]\right)_{t-1}\begin{bmatrix} c \\ d \end{bmatrix} \\ \text{is an equivalence } \forall \begin{bmatrix} c \\ d \end{bmatrix} \end{pmatrix} \Longrightarrow \begin{pmatrix} F\left(\mathcal{O}[j]_t\begin{bmatrix} c \\ d \end{bmatrix}\right) \longrightarrow \left(F^{\mathrm{opd}}\,\mathcal{O}[Fj]\right)_t\begin{bmatrix} c \\ d \end{bmatrix} \\ \text{is an equivalence } \forall \begin{bmatrix} c \\ d \end{bmatrix} \end{pmatrix}.$$

Again by Proposition 3.14, this induction step $(t-1) \Rightarrow (t)$ can be argued via the following cubical diagram

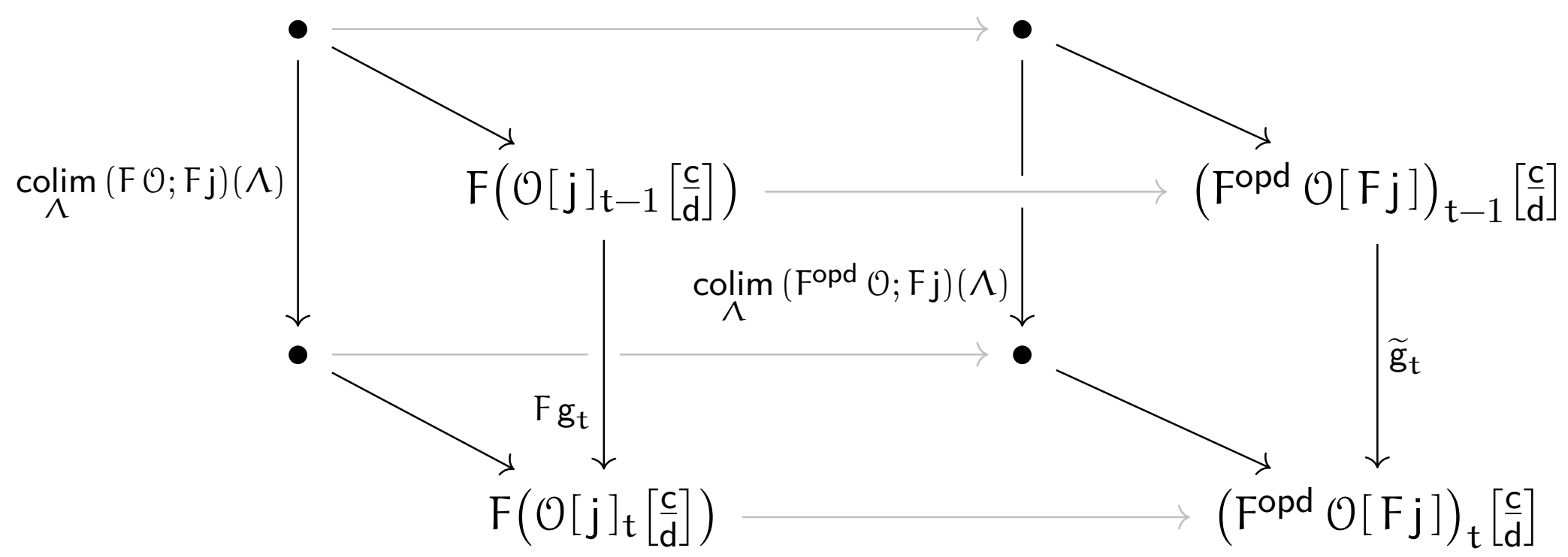

whose black faces are pushouts (note that we used that $F$ is a left adjoint functor and that $F^{\mathrm{opd}}$ commutes with the free operad functor). Observe that we are not decorating the algebra variable in the tagging functor since it is the initial algebra on both sides (and since there are no straight leaves in $\Lambda$). By induction hypothesis, the black faces are homotopy pushouts if we check that

$$\operatorname*{colim}_{\Lambda}(F\,\mathcal{O};\, F j)(\Lambda) \qquad \text{and} \qquad \operatorname*{colim}_{\Lambda}(F^{\mathrm{opd}}\,\mathcal{O};\, F j)(\Lambda)$$

are core cofibrations $(1)$. If that is the case, we would finish the proof by showing that the horizontal arrows in the square

$$\operatorname*{colim}_{\Lambda}(F\,\mathcal{O};\, F j)(\Lambda) \Longrightarrow \operatorname*{colim}_{\Lambda}(F^{\mathrm{opd}}\,\mathcal{O};\, F j)(\Lambda),$$

sitting in the back of the cube, are equivalences $(2)$. For $(1)$, use the proof of Lemma 3.36. For $(2)$, note that $F\,\mathcal{O} \to F^{\mathrm{opd}}\,\mathcal{O}$ is an equivalence in $\Sigma\,\mathcal{C}\mathrm{oll}_{\mathsf{O}}(\mathcal{V}')$ by induction hypothesis and hence, the map in $\mathcal{V}^2$

$$(F\,\mathcal{O}; F j)(\Lambda) \Longrightarrow (F^{\mathrm{opd}}\,\mathcal{O}; F j)(\Lambda)$$

is an equivalence $\forall \Lambda$ since $F \dashv R$ is a weak sm-Quillen adjunction. Applying $\mathrm{colim}_{\Lambda}$ in the arrow category $\mathcal{V}^2$ preserves this equivalence by $\mathrm{Aut}(\Lambda)$-cofibrancy of the factors (see the proof of Lemma 3.36). ☐

From the proof of Proposition 4.25 we can extract:

**Corollary 4.26.** *Let* $F\colon \mathcal{V} \rightleftarrows \mathcal{V}' : R$ *be a weak sm-Quillen equivalence between homotopy cosmoi s.t. the units* $\mathbb{I}_{\mathcal{V}}$ *and* $\mathbb{I}_{\mathcal{V}'}$ *are cofibrant. Then, for any operad* $\mathcal{O} \in \mathcal{O}\mathrm{pd}_{\mathsf{O}}(\mathcal{V})$*, the pointwise comparison map*

$$\mathbb{L}F(\mathcal{O}\begin{bmatrix} c \\ d \end{bmatrix}) \longrightarrow (\mathbb{L}F^{\mathrm{opd}}\,\mathcal{O})\begin{bmatrix} c \\ d \end{bmatrix}$$

*is an equivalence for any* O*-corolla* $\begin{bmatrix} c \\ d \end{bmatrix}$.

*Proof.* Just note that both functors can be derived by picking a cofibrant replacement $Q\mathcal{O}$ of $\mathcal{O}$ in $\mathcal{O}\mathrm{pd}_{\mathsf{O}}(\mathcal{V})$, since $Q\mathcal{O}$ is $\Sigma$-cofibrant when the monoidal unit is cofibrant. ☐

**Related work:** As already mentioned, Theorem 4.23 appeared in [41] with an additional technical condition. Proposition 4.25 generalizes [26, Theorem 1.5] to the colored and symmetric setting (see also [29]). The same applies to Corollary 4.26 and [26, Proposition 1.2]. Also note that a version of Proposition 4.25, under more hypotheses, have appeared in work of Pavlov-Scholbach ([30, Theorem 8.10]) with a more convoluted proof.

## 4.3 (Relative) left properness

**Definition 4.27.** Let $\mathcal{M}$ be a model category and $\mathsf{K}$ a class of objects in $\mathcal{M}$. We say that $\mathcal{M}$ is *left proper relative to* $\mathsf{K}$, or $\mathsf{K}$-*left proper*, if equivalences between objects in $\mathsf{K}$ are stable under cobase changes along cofibrations.

Rezk observed that left/right properness can be characterized by homotopy invariance of slices. More generally, one has:

**Lemma 4.28.** *Let* $\mathcal{M}$ *be a (semi)model category and* $\mathsf{K}$ *a class of objects in* $\mathcal{M}$*. Then,* $\mathcal{M}$ *is* $\mathsf{K}$*-left proper iff for any equivalence* $f\colon X \to Y$ *between objects in* $\mathsf{K}$*, the Quillen pair* $f_!\colon X{\downarrow}\mathcal{M} \rightleftarrows Y{\downarrow}\mathcal{M}\colon f^*$ *is a Quillen equivalence.*

*Proof.* Recall that $f\colon X \to Y$ induces the Quillen pair $f_!\colon X{\downarrow}\mathcal{M} \rightleftarrows Y{\downarrow}\mathcal{M}\colon f^*$ whose right adjoint $f^*$ is just precomposition with $f$ and $f_!$ is cobase change along $f$. Since $f^*$ creates equivalences, it suffices to see that the unit transformation over a cofibrant object in $X{\downarrow}\mathcal{M}$ is an equivalence. Thus, one deduces the claim by observing that the unit for $f_! \dashv f^*$ over $X \rightarrowtail Z$, is just the cobase change $Z \to Y \amalg_X Z$ of $X \to Y$. □

**Proposition 4.29.** *Let* $\mathcal{O}$ *be a cofibrant operad. Then,* $\mathcal{A}\mathsf{lg}_{\mathcal{O}}(\mathcal{V})$ *is left proper relative to* fgt*-cofibrant* $\mathcal{O}$*-algebras if* $\mathcal{V}$ *satisfies* $\mathbb{I}$*-strong unit axiom.*

*Proof.* We must check that any equivalence $A \to B$ between fgt-cofibrant algebras induces a Quillen equivalence $A{\downarrow}\mathcal{A}\mathsf{lg}_{\mathcal{O}}(\mathcal{V}) \rightleftarrows B{\downarrow}\mathcal{A}\mathsf{lg}_{\mathcal{O}}(\mathcal{V})$ by Lemma 4.28. Use Proposition 3.9 to identify this Quillen adjunction with the change of operad adjunction $\mathcal{A}\mathsf{lg}_{\mathcal{O}_A}(\mathcal{V}) \rightleftarrows \mathcal{A}\mathsf{lg}_{\mathcal{O}_B}(\mathcal{V})$ induced by $\mathcal{O}_A \to \mathcal{O}_B$. Then, the result follows from Lemma 3.38 and rectification of algebras over well-pointed operads (see Proposition 4.16). □

*Remark* 4.30. We are not assuming that $\mathcal{V}$ is left proper in the previous statement. Also, the class of cofibrant operads is not sharp. For example, $\mathsf{Ass}_+$ and $\mathsf{Ass}$ satisfy the conclusion as well as any operad $\mathcal{O}$ equipped with a cofibration from $\mathsf{Ass}_+$ or $\mathsf{Ass}$. Such generalizations can also be proven with the results presented here as in [29], but for the lack of applications we decided to avoid such discussion.

The previous proposition has an interesting consequence that generalizes the main result of Rezk in [33] (see also [26, 29]):

**Corollary 4.31.** *Assume that all objects in* $\mathcal{V}$ *are cofibrant. Then, any category of algebras over a* $\Sigma$*-cofibrant or well-pointed operad* $\mathcal{O}$ *can be replaced by a Quillen equivalent one which is left proper.*

*Proof.* Direct application of Proposition 4.29 and Proposition 4.16. Just choose a cofibrant replacement $Q\mathcal{O} \simeq \mathcal{O}$. □

*Remark* 4.32. In the spirit of Proposition 4.9, one can write down a generalization of this result using pseudo-cofibrant objects if $\mathcal{V}$ satisfies the strong unit axiom. We leave this task as an exercise to the reader in order to avoid repeating ourselves.

We close this subsection by showing that $\mathcal{A}\mathrm{lg}_{\mathcal{O}}(\mathcal{V})$ is (relatively) left proper for more general operads if we assume stronger hypotheses on $\mathcal{V}$. First, let us define an algebraic notion that will ensure $\Sigma$-cofibrancy of the relevant enveloping operads. It should be reminiscent of $\Sigma$-splitness in [17].

Forgetting the action of symmetric groups allows us to pass from operads to non-symmetric operads. In fact, there is an adjunction

$$(\star)^{\Sigma}\colon \mathtt{nsOpd}_{\mathsf{O}}(\mathtt{V}) \rightleftarrows \mathtt{Opd}_{\mathsf{O}}(\mathtt{V}) :\mathrm{fgt}_{\mathrm{ns}} \ ,$$

whose left adjoint is explicitly described, at the level of underlying collections, by

$$\mathcal{Q}^{\Sigma}\begin{bmatrix}\{c_1,\dots,c_m\}\\ d\end{bmatrix} = \coprod_{\sigma\in\Sigma_m}\mathcal{Q}\begin{bmatrix}(c_{\sigma(1)},\dots,c_{\sigma(m)})\\ d\end{bmatrix}.$$

We will denote the composition $(\mathrm{fgt}_{\mathrm{ns}}(\star))^{\Sigma}$ simply by $(\star)^{\Sigma}$.

**Definition 4.33.** An operad $\mathcal{O}\in\mathtt{Opd}_{\mathsf{O}}(\mathtt{V})$ is said to carry a $\Sigma$-*splitting* if it comes with a section $s\colon \mathcal{O}\to\mathcal{O}^{\Sigma}$ for the counit map $\epsilon\colon \mathcal{O}^{\Sigma}\to\mathcal{O}$ (in $\mathtt{Opd}_{\mathsf{O}}(\mathtt{V})$).

*Remark* 4.34. Carrying a $\Sigma$-splitting seems a bit stronger than being $\Sigma$-split in the sense of [17], although we have not compared both notions carefully. Our choice comes from its simplicity in the subsequent exposition. Note the similarity of $\mathcal{O}^{\Sigma}$ with the Hadamard product of operads $\mathcal{O}\boxtimes\mathrm{Ass}_{+}$.

With $\Sigma$-splittings already in the game, we have a simple preliminary result:

**Lemma 4.35.** *Let $\mathcal{O}\in\mathtt{Opd}_{\mathsf{O}}(\mathtt{V})$ be an operad and $A\in\mathtt{Alg}_{\mathcal{O}}(\mathtt{V})$ an algebra. Then,*

- *There is a natural isomorphism $(\mathcal{O}^{\Sigma})_A\cong(\mathcal{O}_A)^{\Sigma}$.*
- *If $\mathcal{O}$ carries a $\Sigma$-splitting, then $\mathcal{O}_A$ carries a $\Sigma$-splitting.*
- *If $\mathcal{O}$ carries a $\Sigma$-splitting and it is locally cofibrant, i.e. $\mathcal{O}\begin{bmatrix}c\\ d\end{bmatrix}$ is cofibrant for any corolla, then $\mathcal{O}$ is $\Sigma$-cofibrant.*

*Proof.* Let us check item by item.

- Comparing the universal properties for both objects (see [10, §4.1.1]) we get

$$\begin{array}{ccc}
\mathrm{Hom}_{\mathcal{O}^{\Sigma}\downarrow\mathtt{Opd}}\left((\mathcal{O}^{\Sigma})_{\epsilon^{*}A},\mathcal{P}\right) & & \mathrm{Hom}_{\mathcal{O}^{\Sigma}\downarrow\mathtt{Opd}}\left((\mathcal{O}_A)^{\Sigma},\mathcal{P}\right)\\
\wr\!\downarrow & & \downarrow\!\wr\\
\mathrm{Hom}_{\mathcal{O}^{\Sigma}\text{-alg}}\left(\epsilon^{*}A,\varphi^{*}\mathbb{0}_{\mathcal{P}}\right) & \xleftrightarrow{\ 1:1\ } & \mathrm{Hom}_{\mathcal{O}_{\mathrm{ns}}\text{-alg}}\left(\epsilon^{*}_{\mathrm{ns}}A,\varphi^{*}_{\mathrm{ns}}\mathbb{0}_{\mathcal{P}}\right)
\end{array},$$

  where $\varphi\colon\mathcal{O}^{\Sigma}\to\mathcal{P}$ (equivalently $\varphi_{\mathrm{ns}}\colon\mathrm{fgt}_{\mathrm{ns}}\,\mathcal{O}\to\mathrm{fgt}_{\mathrm{ns}}\,\mathcal{P}$) is given, $\mathbb{0}_{\mathcal{P}}\in\mathtt{Alg}_{\mathcal{P}}$ denotes the initial $\mathcal{P}$-algebra and when we say morphism of $\mathcal{O}_{\mathrm{ns}}$-algebras we mean that one forgets the $\Sigma$-equivariance.

- Using the splitting $s\colon \mathcal{O} \to \mathcal{O}^{\Sigma}$, one can simply define vertical maps in the following diagram

$$\begin{array}{ccccc} \mathcal{O}^1_A\left[\begin{smallmatrix} c \\ d \end{smallmatrix}\right] & \rightrightarrows & \mathcal{O}^0_A\left[\begin{smallmatrix} c \\ d \end{smallmatrix}\right] & \xrightarrow{\text{colim}} & \mathcal{O}_A\left[\begin{smallmatrix} c \\ d \end{smallmatrix}\right] \\ \downarrow & & \downarrow & & \\ (\mathcal{O}^{\Sigma})^1_A\left[\begin{smallmatrix} c \\ d \end{smallmatrix}\right] & \rightrightarrows & (\mathcal{O}^{\Sigma})^0_A\left[\begin{smallmatrix} c \\ d \end{smallmatrix}\right] & \xrightarrow[\text{colim}]{} & (\mathcal{O}^{\Sigma})_A\left[\begin{smallmatrix} c \\ d \end{smallmatrix}\right] \end{array}$$

inducing a map $s\colon \mathcal{O}_A \to (\mathcal{O}^{\Sigma})_A$ between colimits. This map is a section of $(\mathcal{O}^{\Sigma})_A \to \mathcal{O}_A$ by definition and yields a $\Sigma$-splitting of $\mathcal{O}_A$ by the first item.

- This is a direct consequence of the explicit description of the underlying symmetric collection of $\mathcal{O}^{\Sigma}$ and the proof of Lemma 3.26.

□

**Proposition 4.36.** *Assume that $\mathcal{V}$ is left proper, tractable, h-monoidal ([2, Definition 1.11]), pretty small ([31, Definition 2.1]), finite coproducts in $\mathcal{V}$ are homotopical and that all cofibrant objects in $\mathcal{V}$ are flat*[v]*. Let $\mathcal{O} \in \mathcal{O}\mathrm{pd}(\mathcal{V})$ be an operad which carries a $\Sigma$-splitting. Then, $\mathcal{A}\mathrm{lg}_{\mathcal{O}}(\mathcal{V})$ is left proper relative to $\mathrm{K}_{\mathcal{O}}$, where $\mathrm{K}_{\mathcal{O}}$ is the biggest full subcategory of $\mathcal{A}\mathrm{lg}_{\mathcal{O}}(\mathcal{V})$ such that $\mathcal{O}_{\star}\colon \mathrm{K}_{\mathcal{O}} \to \mathcal{O}\mathrm{pd}(\mathcal{V})$ preserves equivalences.*

*Proof.* We have to show that given a pushout square in $\mathcal{A}\mathrm{lg}_{\mathcal{O}}(\mathcal{V})$

$$\begin{array}{ccc} A & \xrightarrow{f} & B \\ {\scriptstyle g}\downarrow & \ulcorner & \downarrow{\scriptstyle \widetilde{g}} \\ \widetilde{A} & \xrightarrow[\widetilde{f}]{} & \widetilde{B} \end{array} \tag{4.2}$$

where $f$ is an equivalence with $A, B \in \mathrm{K}_{\mathcal{O}}$ and $g$ is a proj-cofibration, then $\widetilde{f}$ is an equivalence.

First note that, as equivalences in $\mathcal{V}$ are closed under retracts, we can consider without loss of generality that $g$ is a cellular proj-cofibration, i.e. $g$ is a transfinite composite of maps $g_{\alpha}\colon A_{\alpha} \to A_{\alpha+1}$ obtained as pushouts of generating cofibrations $j_{\alpha}$ in $\mathcal{V}^{O}$,

$$\begin{array}{ccc} \bullet & \longrightarrow & A_{\alpha} \\ {\scriptstyle \mathcal{O}\circ j_{\alpha}}\downarrow & \ulcorner & \downarrow{\scriptstyle g_{\alpha}} \\ \bullet & \longrightarrow & A_{\alpha+1} \end{array} .$$

Note that the maps $j_{\alpha}$ are core cofibrations by tractability. Due to Proposition 3.17 or equiv. [40, Proposition 4.3.17], $g_{\alpha}$, when viewed in $\mathcal{V}^{O}$ via the forgetful functor, can be described as a $\omega$-transfinite composite:

$$g_{\alpha}\colon A_{\alpha} = A^0_{\alpha} \xrightarrow{g^1_{\alpha}} A^1_{\alpha} \xrightarrow{g^2_{\alpha}} \cdots \longrightarrow A^{\omega}_{\alpha} = A_{\alpha+1} .$$

[v]An object $X \in \mathcal{V}$ is flat if $X \otimes \star\colon \mathcal{V} \to \mathcal{V}$ is homotopical, i.e. preserves equivalences.

Arranging the two filtrations together, $g$ can be written as a long transfinite composite of maps in $\mathcal{V}^{O}$. Using this, the pushout square (4.2) is decomposed into the following commutative diagram in $\mathcal{V}^{O}$:

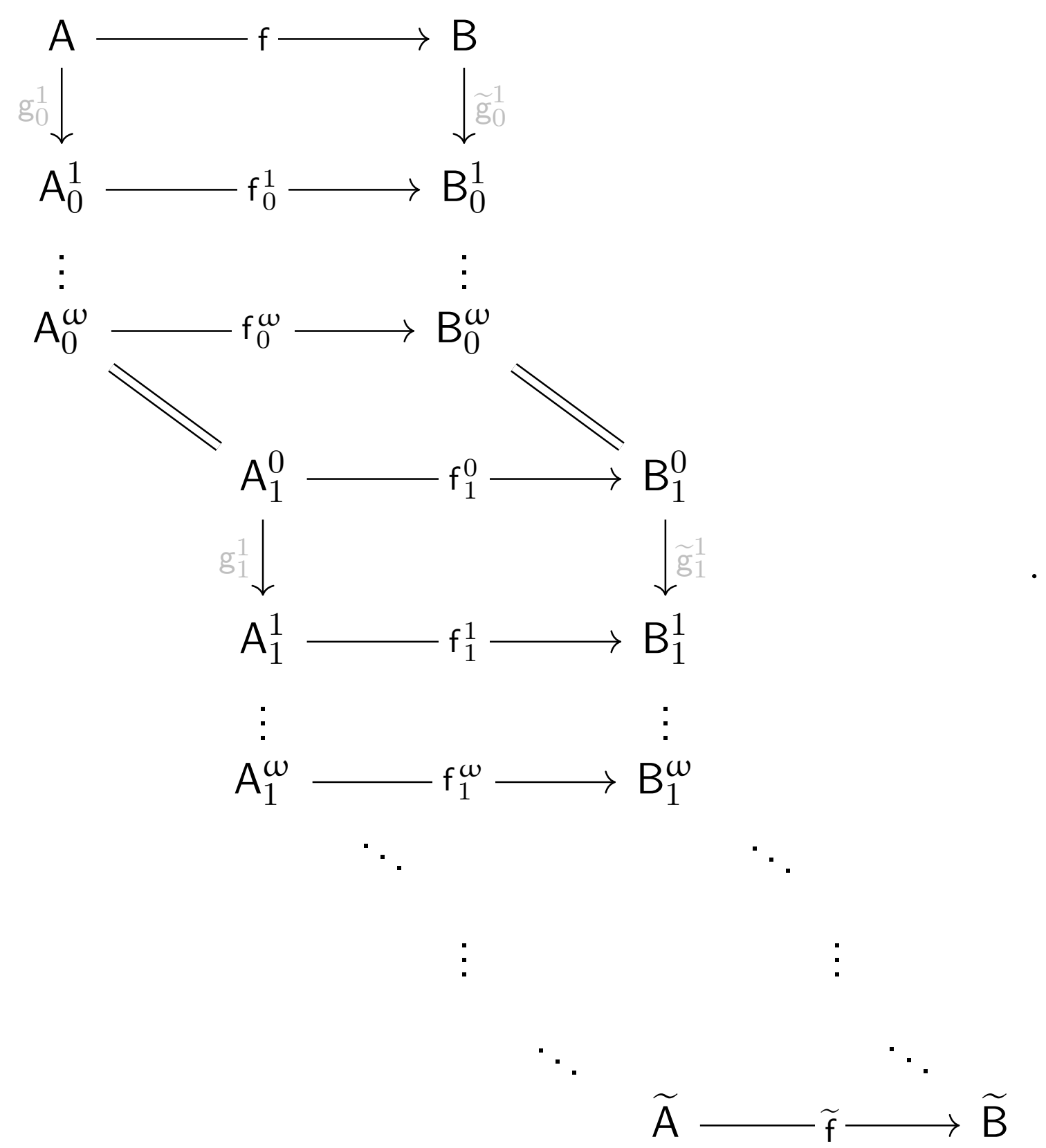

Using that $\mathcal{V}$ is pretty small (and more concretely [31, Lemma 2.2 (iii)]), we obtain that $\widetilde{f}\colon \widetilde{A} \to \widetilde{B}$ is an equivalence (in $\mathcal{V}^{O}$) if we show: $f_\alpha^{t-1}$ is an equivalence $\Rightarrow$ $f_\alpha^{t}$ is an equivalence.

Let us assume for simplicity that $j_\alpha$ is concentrated on a single color $c \in O$ (this simplification does not affect the result). In that case, from Proposition 3.17, we know that $g_\alpha^t$ fits into a pushout in $\mathcal{V}^{O}$ of the form

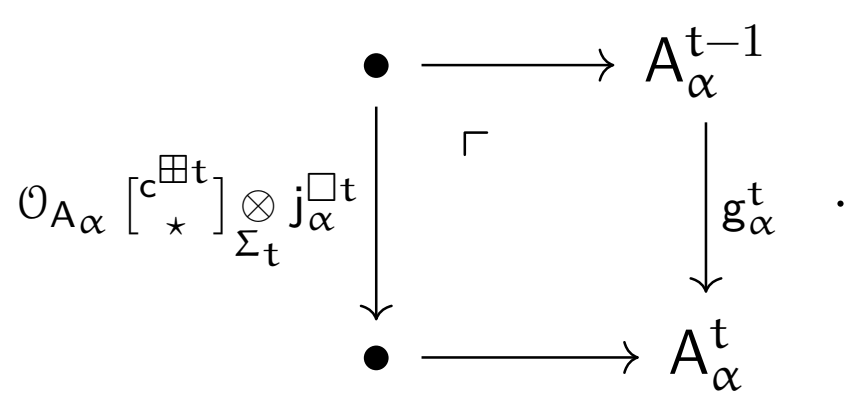

An analogous statement holds for $\widetilde{g}_\alpha^t\colon B_\alpha^{t-1} \to B_\alpha^t$. Now, observe that the original square

([4.2](#)) induces a cube comparing these pushouts

$$
\begin{array}{ccccccc}
\mathcal{O}_{A_\alpha}\begin{bmatrix} c^{\boxplus t} \\ \star \end{bmatrix} \underset{\Sigma_t}{\otimes} s(j_\alpha^{\square t}) & & & \longrightarrow & \mathcal{O}_{B_\alpha}\begin{bmatrix} c^{\boxplus t} \\ \star \end{bmatrix} \underset{\Sigma_t}{\otimes} s(j_\alpha^{\square t}) & & \\
& \searrow & & & & \searrow & \\
{\scriptstyle \mathrm{id} \underset{\Sigma_t}{\otimes} j_\alpha^{\square t}}\downarrow & & A_\alpha^{t-1} & \xrightarrow{\quad f_\alpha^{t-1} \quad} & & & B_\alpha^{t-1} \\
& & & & \downarrow{\scriptstyle \mathrm{id} \underset{\Sigma_t}{\otimes} j_\alpha^{\square t}} & & \\
\mathcal{O}_{A_\alpha}\begin{bmatrix} c^{\boxplus t} \\ \star \end{bmatrix} \underset{\Sigma_t}{\otimes} t(j_\alpha^{\square t}) & & & \longrightarrow & \mathcal{O}_{B_\alpha}\begin{bmatrix} c^{\boxplus t} \\ \star \end{bmatrix} \underset{\Sigma_t}{\otimes} t(j_\alpha^{\square t}) & & \\
& & {\scriptstyle g_\alpha^t}\downarrow & & & & \downarrow{\scriptstyle \tilde{g}_\alpha^t} \\
& \searrow & & & & \searrow & \\
& & A_\alpha^t & \xrightarrow{\quad f_\alpha^t \quad} & & & B_\alpha^t
\end{array}
\qquad . \tag{4.3}
$$

By [[2](#), Proposition 1.8 (b)], to conclude that $f_\alpha^t$ is an equivalence if so is $f_\alpha^{t-1}$, we must show: (a) the vertical maps on the back, $\mathrm{id} \otimes_{\Sigma_t} j_\alpha^{\square t}$, are h-cofibrations; and (b) the square on the back of the cube

$$
\mathcal{O}_{A_\alpha}\begin{bmatrix} c^{\boxplus t} \\ \star \end{bmatrix} \underset{\Sigma_t}{\otimes} j_\alpha^{\square t} \Longrightarrow \mathcal{O}_{B_\alpha}\begin{bmatrix} c^{\boxplus t} \\ \star \end{bmatrix} \underset{\Sigma_t}{\otimes} j_\alpha^{\square t}
$$

is an equivalence in $(\mathcal{V}^{\mathrm{O}})^{2}$, i.e. the gray horizontal maps in that square are equivalences.

Let us start with (a). First, note that $\mathrm{id} \otimes j_\alpha^{\square t}$ is an h-cofibration since $\mathcal{V}$ is h-monoidal. We must show that after taking $\Sigma_t$-coinvariants, we still have an h-cofibration, and for this we use the $\Sigma$-splitting of $\mathcal{O}$. By Lemma [4.35](#), we obtain a retract diagram

$$
\begin{array}{ccccc}
\bullet & \xrightarrow{\quad s\otimes \mathrm{id} \quad} & \bullet & \xrightarrow{\quad \epsilon\otimes \mathrm{id} \quad} & \bullet \\
{\scriptstyle \mathcal{O}_{A_\alpha}\left[\begin{smallmatrix} c^{\boxplus t} \\ \star \end{smallmatrix}\right] \otimes j_\alpha^{\square t}}\downarrow & & \downarrow{\scriptstyle \mathcal{O}^{\Sigma}_{A_\alpha}\left[\begin{smallmatrix} c^{\boxplus t} \\ \star \end{smallmatrix}\right] \otimes j_\alpha^{\square t}} & & \downarrow{\scriptstyle \mathcal{O}_{A_\alpha}\left[\begin{smallmatrix} c^{\boxplus t} \\ \star \end{smallmatrix}\right] \otimes j_\alpha^{\square t}} \\
\bullet & \xrightarrow[\quad s\otimes \mathrm{id} \quad]{} & \bullet & \xrightarrow[\quad \epsilon\otimes \mathrm{id} \quad]{} & \bullet
\end{array}
\; .
$$

The vertical map in the middle is still an h-cofibration ([[2](#), Lemma 1.3]), but it also satisfies that taking $\Sigma_t$-coinvariants

$$
\mathcal{O}^{\Sigma}_{A_\alpha}\begin{bmatrix} c^{\boxplus t} \\ \star \end{bmatrix} \underset{\Sigma_t}{\otimes} j_\alpha^{\square t} \cong \Big( \coprod_{\sigma \in \Sigma_t} \mathcal{O}_{A_\alpha}\begin{bmatrix} c^{\boxplus t} \\ \star \end{bmatrix} \otimes j_\alpha^{\square t} \Big)_{\Sigma_t}
$$

is an h-cofibration. Take $\Sigma_t$-coinvariants over the retract diagram above to deduce that $\mathcal{O}_{A_\alpha}\left[\begin{smallmatrix} c^{\boxplus t} \\ \star \end{smallmatrix}\right] \otimes_{\Sigma_t} j_\alpha^{\square t}$ is an h-cofibration (again by [[2](#), Lemma 1.3]). The same holds if we replace $A_\alpha$ by $B_\alpha$.

Therefore, it remains to check (b). Assume that we already know that $\mathcal{O}_{A_\alpha} \to \mathcal{O}_{B_\alpha}$ is an equivalence of operads (fact which follows from Lemma [4.37](#)). Then, since cofibrant objects in $\mathcal{V}$ are assumed to be flat by hypothesis and since $j_\alpha$ is a core-cofibration, we obtain that before taking $\Sigma_t$-coinvariants

$$\mathcal{O}_{A_\alpha}\begin{bmatrix} c^{\boxplus t} \\ \star \end{bmatrix} \otimes j_\alpha^{\square t} \Longrightarrow \mathcal{O}_{B_\alpha}\begin{bmatrix} c^{\boxplus t} \\ \star \end{bmatrix} \otimes j_\alpha^{\square t}$$

is an equivalence in $(\mathcal{V}^{\mathrm{O}})^{2}$. Using a retract argument similar to the one employed to check $(a)$, one observes that after applying $\Sigma_t$-coinvariants we still have an equivalence (here is where we use that finite coproducts in $\mathcal{V}$ are homotopical). □

**Lemma 4.37.** *Consider a pasting of pushout squares in $\mathcal{A}\mathrm{lg}_{\mathcal{O}}(\mathcal{V})$*

$$\begin{array}{ccccc} \mathcal{O}\circ X & \longrightarrow & A & \longrightarrow & B \\ {\scriptstyle \mathcal{O}\circ j}\downarrow & \ulcorner & \downarrow & \ulcorner & \downarrow \\ \mathcal{O}\circ Y & \longrightarrow & A[j] & \longrightarrow & B[j] \end{array},$$

*where $j$ is a core cofibration in $\mathcal{V}^{\mathrm{O}}$. Under the hypothesis of Proposition 4.36, if the map of operads $\mathcal{O}_A \to \mathcal{O}_B$ is an equivalence, so is $\mathcal{O}_{A[j]} \to \mathcal{O}_{B[j]}$.*

*Proof.* We run the usual induction argument thanks to Proposition 3.17. As in the previous proof (4.36), pretty smallness on $\mathcal{V}$ (and more concretely [31, Lemma 2.2 (iii)]) allows us to focus on the induction steps. For the step $(t-1) \Rightarrow (t)$, we have a cube whose black faces are pushouts. By [2, Proposition 1.8 (b)] and induction hypothesis, we should check: $(1)$ the vertical maps in the back,

$$\operatorname*{colim}_{\tau}\,(\mathcal{O}_A;\, j)(\tau) \qquad \text{and} \qquad \operatorname*{colim}_{\tau}\,(\mathcal{O}_B;\, j)(\tau),$$

are h-cofibrations; and $(2)$ the horizontal arrows in the square

$$\operatorname*{colim}_{\tau}\,(\mathcal{O}_A;\, j)(\tau) \Longrightarrow \operatorname*{colim}_{\tau}\,(\mathcal{O}_B;\, j)(\tau),$$

sitting in the back of the cube, are equivalences. For $(1)$ (resp. $(2)$), use the argument to check $(a)$ (resp. $(b)$) within the proof of Proposition 4.36. In other words, play with $\Sigma$-splittings to pass from the statement at the level of $(\mathcal{O}_A;\, j)(\tau)$ for corollas $\tau$ to the statement after taking colimit over $\tau$. □

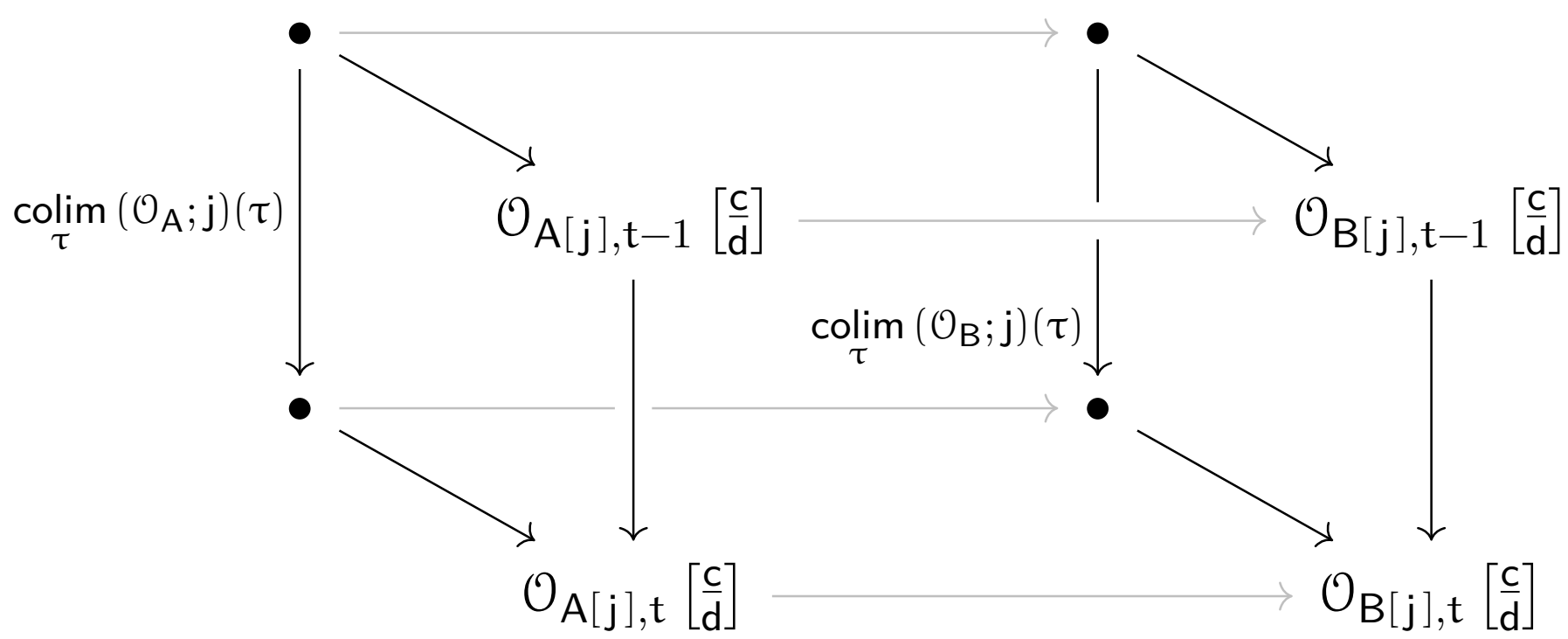

It is possible to replace most of the technical conditions in Proposition 4.36 by a single strong hypothesis:

**Proposition 4.38.** *Assume that all objects in $\mathcal{V}$ are cofibrant. Let $\mathcal{O} \in \mathcal{O}\mathrm{pd}(\mathcal{V})$ be an operad which carries a $\Sigma$-splitting. Then, $\mathcal{A}\mathrm{lg}_{\mathcal{O}}(\mathcal{V})$ is left proper relative to $\mathsf{K}_{\mathcal{O}}$, where $\mathsf{K}_{\mathcal{O}}$ is the biggest full subcategory of $\mathcal{A}\mathrm{lg}_{\mathcal{O}}(\mathcal{V})$ so that $\mathcal{O}_\star \colon \mathsf{K}_{\mathcal{O}} \to \mathcal{O}\mathrm{pd}(\mathcal{V})$ preserves equivalences.*

*Proof.* The proof follows the same pattern as the one of Proposition 4.36. The modifications are:

- Instead of pretty smallness to argue that it suffices to show:
$$(f_\alpha^{t-1} \text{ is an equivalence} \quad \Rightarrow \quad f_\alpha^t \text{ is an equivalence}),$$
we observe that the vertical arrows $g_\alpha^t$, $\widetilde{g}_\alpha^t$ are cofibrations $\forall \alpha, t$ (and that A, B are automatically fgt-cofibrant by hypothesis).

- To check that $g_\alpha^t$, $\widetilde{g}_\alpha^t$ are cofibrations, we show that $\mathrm{id} \otimes_{\Sigma_t} j_\alpha^{\square t}$ are cofibrations. By the third item in Lemma 4.35 (again, all objects in $\mathcal{V}$ are cofibrant), this is just an application of Lemma 3.35.

- By the previous item and since all objects in $\mathcal{V}$ are cofibrant, the black faces of (4.3) are homotopy pushouts. Hence, it remains to prove that
$$\mathcal{O}_{A_\alpha}\begin{bmatrix} c^{\boxplus t} \\ \star \end{bmatrix} \underset{\Sigma_t}{\otimes} j_\alpha^{\square t} \Longrightarrow \mathcal{O}_{B_\alpha}\begin{bmatrix} c^{\boxplus t} \\ \star \end{bmatrix} \underset{\Sigma_t}{\otimes} j_\alpha^{\square t}$$
is an equivalence in $(\mathcal{V}^{\mathsf{O}})^{\mathbf{2}}$. Use again the third item in Lemma 4.35 to see that the claim reduces to have an equivalence of operads $\mathcal{O}_{A_\alpha} \to \mathcal{O}_{B_\alpha}$.

- The analog of Lemma 4.37 can be proven by using the replacements given in the previous items.

□

*Example* 4.39. The Kan-Quillen model structure $\mathcal{S}\mathrm{pc}$ and the projective model structure $\mathcal{C}\mathrm{h}(\mathbb{F})$ on chain complexes over an arbitrary field $\mathbb{F}$ satisfy the hypotheses in Proposition 4.36 and 4.38. More examples can be extracted from [31].

*Remark* 4.40. In Proposition 4.36, one can drop the assumption on $\Sigma$-splittings if taking coinvariants in $\mathcal{V}$ for finite group actions preserves equivalences and if $\mathcal{V}$ is symmetric h-monoidal instead of just h-monoidal; for instance, when $\mathcal{V} = \mathcal{C}\mathrm{h}(\Bbbk)$ and $\Bbbk$ is a field of characteristic zero. Similarly, in Proposition 4.38, $\Sigma$-splittings can be replaced by the following assumption: taking coinvariants in $\mathcal{V}$ for finite group actions preserve equivalences and injective-cofibrations.

**Related work:** First, it is important to point out that categories of operadic algebras are not always left proper (see [9, 13, 29, 33]).

Proposition 4.29 and Corollary 4.31 are new and actually generalize: (1) the main result in [33] to any homotopy cosmos $\mathcal{V}$ where all objects are cofibrant; (2) [26, Theorem 1.13 and Corollary 1.14] to the colored and symmetric setting.

An even stronger result than Proposition 4.29 was claimed as part of [36, Theorem 4.3], but the proof contains a major gap that invalidates it. The issue concerns [36, Proposition 4.5], and the reason is the same as in Remark 4.10.

Propositions 4.36 and 4.38 are new. Both results are designed to show how one uses enveloping operads to deal with left-properness. They are inspired by Batanin-Berger's work on algebras over polynomial monads (e.g. [2, Theorem 0.1]). However, notice that our proofs do not go beyond classical operadic techniques and they apply to honest $\mathcal{V}$-operads.

## 4.4 Moduli of algebra structures

One of the central results in Rezk's thesis, [34, Theorem 1.1.5], relates the two natural ways to define a moduli space of derived $\mathcal{O}$-algebra structures on a fixed object when $\mathtt{V} = \mathtt{sSet}$ or $\mathtt{V} = \mathtt{sMod}(\mathtt{R})$. In this subsection, we show that the coincidence of both "moduli spaces" of $\mathcal{O}$-algebra structures holds essentially for all $\mathcal{V}$.

We assume in this subsection that $\mathcal{V}$ satisfies the $\mathbb{I}$-strong unit axiom (recall that this holds for instance if $\mathbb{I}$ is cofibrant). By [28], this is really a mild condition. We require this axiom because we will use rectification of cofibrant operads (Corollary 4.18).

The crux of the argument is a simple homotopical analysis of endomorphism operads combined with strong (semi)admissibility results of operadic algebras.

**Endomorphism operads.** Recall that the internal hom $\underline{\mathrm{Hom}}_{\mathcal{V}}\colon \mathcal{V}^{\mathrm{op}} \times \mathcal{V} \to \mathcal{V}$ yields an enriched hom bifunctor between diagrams $\underline{\mathrm{Hom}}_{\mathcal{V}^{\mathrm{D}}}\colon (\mathcal{V}^{\mathrm{D}})^{\mathrm{op}} \times \mathcal{V}^{\mathrm{D}} \to \mathcal{V}$; more concretely, if $\mathrm{Y}, \mathrm{Y}'\colon \mathrm{D} \rightrightarrows \mathcal{V}$ are two diagrams, one defines the enriched hom via the following end

$$\underline{\mathrm{Hom}}_{\mathcal{V}^{\mathrm{D}}}(\mathrm{Y}, \mathrm{Y}') := \int_{i \in \mathrm{D}} \underline{\mathrm{Hom}}_{\mathcal{V}}(\mathrm{Y}(i), \mathrm{Y}'(i)).$$

Using this gadget, it is straightforward to define endomorphism operads of diagrams: given $\mathrm{Y}\colon \mathrm{D} \to \mathcal{V}^{\mathrm{O}}$, there is an operad $\mathfrak{End}_{\mathrm{D}}(\mathrm{Y}) \in \mathcal{O}\mathrm{pd}_{\mathrm{O}}(\mathcal{V})$ with

$$\mathfrak{End}_{\mathrm{D}}(\mathrm{Y})\begin{bmatrix} \underline{c} \\ d \end{bmatrix} := \underline{\mathrm{Hom}}_{\mathcal{V}^{\mathrm{D}}}\Big(\bigotimes_k \mathrm{Y}(c_k), \mathrm{Y}(d)\Big).$$

Moreover, for any $\mathcal{O} \in \mathcal{O}\mathrm{pd}_{\mathrm{O}}(\mathcal{V})$ there is a natural bijection

$$\left\{ \begin{matrix} \text{morphism in } \mathcal{O}\mathrm{pd}_{\mathrm{O}}(\mathcal{V}) \\ \mathcal{O} \to \mathfrak{End}_{\mathrm{D}}(\mathrm{Y}) \end{matrix} \right\} \xleftrightarrow{1:1} \left\{ \begin{matrix} \text{diagram of } \mathcal{O}\text{-algebras} \\ \mathrm{D} \to \mathcal{A}\mathrm{lg}_{\mathcal{O}}(\mathcal{V}) \end{matrix} \right\}. \tag{4.4}$$

A functor $z\colon \mathrm{D} \to \mathrm{D}'$ between indexing categories induces a $\mathcal{V}$-natural transformation $z^*\colon \underline{\mathrm{Hom}}_{\mathcal{V}^{\mathrm{D}}}(*, \star) \to \underline{\mathrm{Hom}}_{\mathcal{V}^{\mathrm{D}'}}(*z, \star z)$ and in particular, a morphism of operads

$$z^*\colon \mathfrak{End}_{\mathrm{D}}(\mathrm{Y}) \longrightarrow \mathfrak{End}_{\mathrm{D}'}(\mathrm{Y}z)$$

for any diagram $Y : D \to \mathcal{V}^O$. The previous natural bijection can be specialized to yield a 1-to-1 correspondence between $\mathcal{O}$-algebra maps $f\colon (X, \gamma_X) \to (Y, \gamma_Y)$ and commutative diagrams in $\mathcal{O}pd_O(\mathcal{V})$

$$\begin{tikzcd} & \mathcal{O} \arrow[d] \arrow[ddl, bend right, "\gamma_X"'] \arrow[ddr, bend left, "\gamma_Y"] & \\ & \mathfrak{End}_{[1]}\left(X \xrightarrow{f} Y\right) \arrow[dl, "d_1"'] \arrow[dr, "d_0"] & \\ \mathfrak{End}_{[0]}(X) & & \mathfrak{End}_{[0]}(Y) \end{tikzcd}.$$

Since $\underline{\mathrm{Hom}}_{\mathcal{V}}$ is a right Quillen bifunctor, $\mathfrak{End}_D$ inherits a nice homotopical behavior. For instance, if $Y \in \mathcal{V}^O$ is bifibrant, $\mathfrak{End}_{[0]}(Y)$ is a fibrant operad. Furthermore, due to the following pullback diagram (which comes just from the definition of endomorphism operad as an end)

$$\begin{tikzcd} & \mathfrak{End}_{[1]}\left(X \xrightarrow{f} Y\right)\left[\begin{smallmatrix} \underline{c} \\ d \end{smallmatrix}\right] \arrow[dl, "d_1"'] \arrow[dr, "d_0"] & \\ \mathfrak{End}_{[0]}(X)\left[\begin{smallmatrix} \underline{c} \\ d \end{smallmatrix}\right] \arrow[d, equal] & & \mathfrak{End}_{[0]}(Y)\left[\begin{smallmatrix} \underline{c} \\ d \end{smallmatrix}\right] \arrow[d, equal] \\ \underline{\mathrm{Hom}}_{\mathcal{V}}\left(\bigotimes_k X(c_k), X(d)\right) \arrow[dr, "f(d)_*"'] & & \underline{\mathrm{Hom}}_{\mathcal{V}}\left(\bigotimes_k Y(c_k), Y(d)\right) \arrow[dl, "\left(\bigotimes_k f(c_k)\right)^*"] \\ & \underline{\mathrm{Hom}}_{\mathcal{V}}\left(\bigotimes_k X(c_k), Y(d)\right) & \end{tikzcd},$$

one can show:

**Proposition 4.41.** *[3.6 in [27]] Let $f\colon X \to Y$ be a map in $\mathcal{V}^O$. Then,*

- *If $f$ is a(n acyclic) fibration and $X$ is cofibrant, $d_0$ is a(n acyclic) fibration.*
- *If $f$ is a core (acyclic) cofibration and $Y$ is fibrant, $d_1$ is a core (acyclic) fibration.*
- *If $f$ is an acyclic fibration between bifibrant objects, $d_1$ is a weak equivalence between fibrant objects.*
- *If $f$ is an acyclic cofibration between bifibrant objects, $d_0$ is a weak equivalence between fibrant objects.*

*Remark* 4.42. This simple result is the crux of [11, Part II], at least for its operadic analogue. When working with algebras over props, the analog of $\mathfrak{End}_D(Y)$ requires to consider tensor products also in the codomain, complicating the discussion. That is the reason why Fresse in loc.cit. must deal with the interaction of (acyclic) fibrations and tensor products.

Iterating Proposition 4.41, one obtains:

**Corollary 4.43.** *[3.7 in [27]] Let $X_\bullet \colon [n+1] \to \mathcal{V}^{\mathrm{O}}$ be a chain of acyclic fibrations between bifibrant objects for $n \geqslant 0$. Then, the induced maps*

$$\begin{array}{ccc} & \mathfrak{End}_{[n+1]}(X_0 \to \cdots \to X_n \to X_{n+1}) & \\ \swarrow^{d_{n+1}} & & \searrow^{d_0^{n+1}} \\ \mathfrak{End}_{[n]}(X_0 \to \cdots \to X_n) & & \mathfrak{End}_{[0]}(X_{n+1}) \end{array}$$

*are a weak equivalence ($d_{n+1}$) and an acyclic fibration ($d_0^{n+1}$) between fibrant operads.*

**Classifying spaces and spaces of algebra structures.** In modern terms, Rezk's original result [34, Theorem 1.1.5], and Muro's non-symmetric generalization [27, Theorem 4.6], compute the homotopy fiber of the forgetful ($\infty$-)functor $\mathcal{A}\mathrm{lg}_{\mathcal{O}}(\mathcal{V})^{\simeq}_{\infty} \to \mathcal{V}^{\simeq}_{\infty}$, where the subscript "$\infty$" denotes associated $\infty$-category and the superscript "$\simeq$" refers to the maximal $\infty$-groupoid inside an $\infty$-category. To actually be able to compute this homotopy fiber in terms of strict algebra structures, both authors make use of classifying spaces of model categories.

Let us assume once and for all that $\mathcal{O} \in \mathcal{O}\mathrm{pd}_{\mathrm{O}}(\mathcal{V})$ is $\Sigma$-cofibrant or well-pointed in the rest of this subsection. From Proposition 4.2, we know that $\mathcal{A}\mathrm{lg}_{\mathcal{O}}(\mathcal{V})$ exists as a semimodel category at least, and this will suffice for our purposes. We will also use that $\mathrm{fgt} \colon \mathcal{A}\mathrm{lg}_{\mathcal{O}}(\mathcal{V}) \to \mathcal{V}^{\mathrm{O}}$ preserves core (acyclic) cofibrations (Proposition 4.3). In other words, $\Sigma$-cofibrant operads are *strongly semiadmissible*.

Let $w\mathcal{A}\mathrm{lg}_{\mathcal{O}}(\mathcal{V})$ denote the category of $\mathcal{O}$-algebras with weak equivalences between them and $|w\mathcal{A}\mathrm{lg}_{\mathcal{O}}(\mathcal{V})|$ be the nerve of this category (thought as a space). Then, by general considerations $|w\mathcal{A}\mathrm{lg}_{\mathcal{O}}(\mathcal{V})|$ is weak homotopy equivalent to $\mathcal{A}\mathrm{lg}_{\mathcal{O}}(\mathcal{V})^{\simeq}_{\infty}$ (when interpreting both objects in $\mathtt{sSet}$).[vi] Of course, a similar discussion applies to $\mathcal{V}$. Then, our goal translates into computing the homotopy fiber of the map

$$\mathrm{fgt} \colon |w\mathcal{A}\mathrm{lg}_{\mathcal{O}}(\mathcal{V})| \to |w\mathcal{V}^{\mathrm{O}}|.$$

The idea is to replace fgt by another map suitable for applying a version of Quillen's Theorem B.

First, consider a cosimplicial resolution $\mathcal{O}^\bullet \xrightarrow{\sim} cc_* \mathcal{O}$ in $\mathcal{O}\mathrm{pd}_{\mathrm{O}}(\mathcal{V})$ (see [18, Definition 16.1.2]). Taking categories of algebras we obtain a simplicial category $w\mathcal{A}\mathrm{lg}_{\mathcal{O}^\bullet}(\mathcal{V})$. Secondly, by the homotopy invariance of endomorphism operads, Corollary 4.43, it is convenient to consider the subcategory $fw\mathcal{V}_{\mathsf{cf}} \hookrightarrow w\mathcal{V}$ spanned by bifibrant objects and acyclic fibrations between them. With these two ingredients, we form the following commutative diagram of simplicial categories

$$\begin{array}{ccccc} \mathcal{M}_\bullet & \longrightarrow & w\mathcal{A}\mathrm{lg}_{\mathcal{O}^\bullet}(\mathcal{V}) & \longleftarrow & w\mathcal{A}\mathrm{lg}_{\mathcal{O}}(\mathcal{V}) \\ \downarrow{\scriptstyle \widetilde{\mathrm{fgt}}_\bullet} & \lrcorner & \downarrow{\scriptstyle \mathrm{fgt}_\bullet} & & \downarrow{\scriptstyle \mathrm{fgt}} \\ fw\mathcal{V}^{\mathrm{O}}_{\mathsf{cf}} & \hookrightarrow & w\mathcal{V}^{\mathrm{O}} & = & w\mathcal{V}^{\mathrm{O}} \end{array}, \tag{4.5}$$

[vi]We will avoid discussing size issues, as the original references do, since they can be dealt with standard arguments.

where: (1) the terms without "•" are discrete simplicial categories, (2) the left square is a pullback (i.e. $\mathcal{M}_t$ is the category of fgt-bifibrant $\mathcal{O}^t$-algebras with acyclic fibrations between them).

**Lemma 4.44.** *The horizontal functors in (4.5) induce weak homotopy equivalences*

$$|fw\mathcal{V}^{\mathsf{O}}_{\mathsf{cf}}| \simeq |w\mathcal{V}^{\mathsf{O}}|, \qquad |\mathcal{M}_\bullet| \simeq |w\mathcal{A}\mathsf{lg}_{\mathcal{O}^\bullet}(\mathcal{V})| \simeq |w\mathcal{A}\mathsf{lg}_{\mathcal{O}}(\mathcal{V})|.$$ [vii]

*Proof.* The first weak homotopy equivalence can be proven as in [34]. For the last weak homotopy equivalence, $|w\mathcal{A}\mathsf{lg}_{\mathcal{O}^\bullet}(\mathcal{V})| \simeq |w\mathcal{A}\mathsf{lg}_{\mathcal{O}}(\mathcal{V})|$, use the Quillen equivalences

$$\mathcal{A}\mathsf{lg}_{\mathcal{O}^t}(\mathcal{V}) \underset{\mathsf{Q}}{\simeq} \mathcal{A}\mathsf{lg}_{\mathcal{O}^0}(\mathcal{V}) \underset{\mathsf{Q}}{\simeq} \mathcal{A}\mathsf{lg}_{\mathcal{O}}(\mathcal{V})$$

coming from rectification of $\Sigma$-cofibrant and well-pointed operads, Proposition 4.16 (recall that we assumed the $\mathbb{I}$-strong unit axiom on $\mathcal{V}$). To check the last one, observe that $\mathsf{fgt}_t\colon \mathcal{A}\mathsf{lg}_{\mathcal{O}^t}(\mathcal{V}) \to \mathcal{V}^{\mathsf{O}}$ preserves cofibrant objects, Proposition 4.3. Thus, using essentially the same idea as for the first equivalence, one shows

$$|\mathcal{M}_t| \simeq |fw\mathcal{A}\mathsf{lg}_{\mathcal{O}^t}(\mathcal{V})_{\mathsf{cf}}| \simeq |w\mathcal{A}\mathsf{lg}_{\mathcal{O}^t}(\mathcal{V})|,$$

where the simplicial set in the middle corresponds to the category of bifibrant $\mathcal{O}$-algebras with acyclic fibrations between them. □

**Lemma 4.45.** *[4.5 in [27]] Let $\mathsf{Y} \in \mathcal{V}^{\mathsf{O}}_{\mathsf{cf}}$ be a bifibrant object and $\mathcal{O}^\bullet$ be a cosimplicial resolution of $\mathcal{O}$ in $\mathsf{Opd}_{\mathsf{O}}(\mathcal{V})$. Then, $|\mathcal{M}_\bullet \downarrow \mathsf{Y}|$[viii] is weak homotopy equivalent to the derived mapping space $\mathbb{R}\mathsf{Map}_{\mathsf{Opd}_{\mathsf{O}}(\mathcal{V})}\big(\mathcal{O}, \mathfrak{End}_{[0]}(\mathsf{Y})\big)$.*

*Proof.* Let us sketch the proof of [27, Lemma 4.5] that also applies to this case.

By the natural bijection (4.4), we can identify the set of $(s,t)$-bisimplices of $|\mathcal{M}_\bullet \downarrow \mathsf{Y}|$:

$$|\mathcal{M}_t \downarrow \mathsf{Y}|_s \cong \coprod_{\mathsf{X}_\bullet\colon [s]\to fw\mathcal{V}_{\mathsf{cf}}\downarrow \mathsf{Y}} \mathsf{Hom}_{\mathsf{Opd}_{\mathsf{O}}(\mathcal{V})}\big(\mathcal{O}^t, \mathfrak{End}_{[s]}(\mathsf{X}_0 \to \cdots \to \mathsf{X}_s)\big).$$

Making use of the weak equivalences in Corollary 4.43, we identify $|\mathcal{M}_\bullet \downarrow \mathsf{Y}|$ with $|fw\mathcal{V}_{\mathsf{cf}} \downarrow \mathsf{Y}| \times \mathsf{Hom}_{\mathsf{Opd}_{\mathsf{O}}(\mathcal{V})}\big(\mathcal{O}^\bullet, \mathfrak{End}_{[0]}(\mathsf{Y})\big)$ up to weak homotopy equivalence, i.e. with the bisimplicial set

$$(s,t) \longmapsto \coprod_{\mathsf{X}_\bullet\colon [s]\to fw\mathcal{V}_{\mathsf{cf}}\downarrow \mathsf{Y}} \mathsf{Hom}_{\mathsf{Opd}_{\mathsf{O}}(\mathcal{V})}\big(\mathcal{O}^t, \mathfrak{End}_{[0]}(\mathsf{Y})\big).$$

Note that we have used that the relevant endormorphism operads are all fibrant in this case and that $\mathcal{O}^t$ is cofibrant $\forall t \geqslant 0$. Since $|fw\mathcal{V}_{\mathsf{cf}} \downarrow \mathsf{Y}|$ is weakly contractible (it has a final object) and since $\mathsf{Hom}_{\mathsf{Opd}_{\mathsf{O}}(\mathcal{V})}\big(\mathcal{O}^\bullet, \mathfrak{End}_{[0]}(\mathsf{Y})\big)$ is a model for the derived mapping space in the statement, we conclude the claim. □

Finally,

[vii] The bisimplicial sets in this statement are interpreted as simplicial sets via the diagonal construction.
[viii] Again, we see this bisimplicial set as a simplicial set by taking the diagonal construction.

**Theorem 4.46.** *Let $\mathcal{O}$ be a $\Sigma$-cofibrant or well-pointed operad. Then, the homotopy fiber of*

$$\mathrm{fgt}\colon |w\mathcal{A}\mathrm{lg}_{\mathcal{O}}(\mathcal{V})| \longrightarrow |w\mathcal{V}^{\mathsf{O}}|$$

*over a bifibrant object* $\mathsf{Y} \in \mathcal{V}^{\mathsf{O}}$ *is* $\mathbb{R}\mathrm{Map}_{\mathcal{O}\mathsf{pd}_{\mathsf{O}}(\mathcal{V})}\big(\mathcal{O}, \mathfrak{End}_{[0]}(\mathsf{Y})\big)$.

*Proof.* Due to Lemma 4.44, $\mathrm{hofib}_{\mathsf{Y}}(\mathrm{fgt}) \simeq \mathrm{hofib}_{\mathsf{Y}}(\widetilde{\mathrm{fgt}_{\bullet}})$. To finish the proof, apply [34, Lemma 4.2.2], which is the announced variation of Quillen's Theorem B, to compute the homotopy fiber of $\widetilde{\mathrm{fgt}_{\bullet}}\colon \mathcal{M}_{\bullet} \to \mathrm{f}w\mathcal{V}^{\mathsf{O}}_{\mathrm{cf}}$ at $\mathsf{Y}$ thanks to Lemma 4.45. □

*Remark* 4.47. The analog of Theorem 4.46 is the starting point in [27] to deduce geometric properties of the derived stack of $\mathcal{O}$-algebra structures, $\underline{\mathcal{A}\mathrm{lg}}_{\mathcal{O}}$, in a general HAG context. Muro works with non-symmetric operads with a single color, but it seems that several of his results can be generalized to cover symmetric (colored) operads.

**Related work:** Theorem 4.46 generalizes: (1) [34, Theorem 1.1.5] to any cofibrantly generated closed sm-model category $\mathcal{V}$ satisfying the $\mathbb{I}$-strong unit axiom; (2) [27, Theorem 4.6] to the colored and symmetric setting. Notice that Muro in loc.cit. assumes the monoid axiom and the strong unit axiom on $\mathcal{V}$, while we only assume that $\mathcal{V}$ satisfies the $\mathbb{I}$-strong unit axiom.

Yalin in [42] showed that an analogue of Theorem 4.46 holds for algebras over cofibrant props, but he crucially needs to work on $\mathcal{V} = \mathcal{C}h(\Bbbk)$ over a field $\Bbbk$ (or similar settings).

## 4.5 Homotopy invariance of modules

For any $\mathcal{O}$-algebra A, one can define a category of A-modules (under $\mathcal{O}$), denoted $\mathtt{Mod}^{\mathcal{O}}_{\mathsf{A}}(\mathtt{V})$ (see [5]). Varying the operad $\mathcal{O}$ one gets a global category of operadic algebras $\mathtt{OpdAlg}(\mathtt{V})$ and similarly, varying the pair $(\mathcal{O}, \mathsf{A})$ one obtains a global category of operadic modules $\mathtt{OpdMod}(\mathtt{V})$. We may visualize this situation by the following diagram of categories

$$\begin{array}{ccccccc}
(\mathcal{O},\mathsf{A},\mathsf{M}) \in & \mathtt{OpdMod}(\mathtt{V}) & \hookleftarrow & \mathtt{Mod}^{\mathcal{O}}(\mathtt{V}) & \hookleftarrow & \mathtt{Mod}^{\mathcal{O}}_{\mathsf{A}}(\mathtt{V}) \\
\mapsto\downarrow & \downarrow & & \downarrow & & \downarrow \\
(\mathcal{O},\mathsf{A}) \in & \mathtt{OpdAlg}(\mathtt{V}) & \hookleftarrow & \mathtt{Alg}_{\mathcal{O}}(\mathtt{V}) & \hookleftarrow & \{(\mathcal{O},\mathsf{A})\} \\
\mapsto\downarrow & \downarrow & & \downarrow & & \\
\mathcal{O} \in & \mathtt{Opd}(\mathtt{V}) & \hookleftarrow & \{\mathcal{O}\} & &
\end{array}\cdot$$

In fact, the functor $\mathtt{OpdAlg}(\mathtt{V}) \to \mathtt{Opd}(\mathtt{V})$ is the Grothendieck construction of the assignment $\mathcal{O} \mapsto \mathtt{Alg}_{\mathcal{O}}(\mathtt{V})$ and analogously $\mathtt{Mod}^{\mathcal{O}}(\mathtt{V}) \to \mathtt{Alg}_{\mathcal{O}}(\mathtt{V})$ is the Grothendieck construction of $\mathsf{A} \mapsto \mathtt{Mod}^{\mathcal{O}}_{\mathsf{A}}(\mathtt{V})$. The rectification of operadic algebras asks for the homotopy invariance of $\mathcal{O} \mapsto \mathcal{A}\mathrm{lg}_{\mathcal{O}}(\mathcal{V})$. We want to address now the homotopy invariance of $\mathsf{A} \mapsto \mathcal{M}\mathrm{od}^{\mathcal{O}}_{\mathsf{A}}(\mathcal{V})$.

The key point to deduce such homotopy invariance from our discussion of enveloping operads comes from [5, Theorem 1.10] that we enunciate now.

**Proposition 4.48.** *There is an isomorphism of categories*

$$\mathtt{Mod}^{\mathcal{O}}_{A}(\mathtt{V}) \cong \mathtt{Alg}_{\overline{\mathcal{O}}_{A}}(\mathtt{V}) = \left[\overline{\mathcal{O}}_{A}, \mathtt{V}\right] ,$$

*where $\overline{\mathcal{O}}_A$ denotes the underlying category of the enveloping operad $\mathcal{O}_A$.*

Combining this result with the homotopical analysis in §3.3, we find:

**Proposition 4.49.** *The assignment $A \mapsto \mathcal{M}od^{\mathcal{O}}_{A}(\mathcal{V})$ is homotopical, i.e. it sends equivalences of algebras to Quillen equivalences, if one of the following conditions hold:*

- *the operad $\mathcal{O}$ is $\Sigma$-cofibrant or well-pointed and the assignment is restricted to proj-cofibrant algebras;*
- *the operad $\mathcal{O}$ is cofibrant and the assignment is restricted to* fgt*-cofibrant algebras.*

*Proof.* Since $\mathcal{O}_A$ is $\Sigma$-cofibrant (resp. well-pointed) and $A \mapsto \mathcal{O}_A$ is homotopical in both cases (Propositions 3.34 and 3.38), $A \mapsto \overline{\mathcal{O}}_A$ is homotopical and lands in locally cofibrant (resp. well-pointed) categories. Notice that assigning functor categories out of these categories is a homotopical construction (see [5, Proposition 2.7 (b)]). ☐

We can even consider that the underlying operad varies:

**Proposition 4.50.** *A map $(\mathcal{O}, A) \to (\mathcal{P}, B)$ in $\mathtt{OpdAlg}(\mathtt{V})$ such that $\phi\colon \mathcal{O} \to \mathcal{P}$ is an equivalence in $\mathcal{O}\mathrm{pd}_{O}(\mathcal{V})$ induces a Quillen equivalence*

$$\mathcal{M}od^{\mathcal{O}}_{A}(\mathcal{V}) \simeq_{Q} \mathcal{M}od^{\mathcal{P}}_{B}(\mathcal{V})$$

*if one of the following conditions holds:*

(i) *The operads $\mathcal{O}$, $\mathcal{P}$ are $\Sigma$-cofibrant and either:*

  (a) *$A$ is proj-cofibrant and the map of $\mathcal{P}$-algebras is the identity $\phi_\sharp A \to \phi_\sharp A$;*

  (b) *$B$ is* fgt*-cofibrant and the map of $\mathcal{O}$-algebras is a cofibrant replacement $q\colon Q\phi^* B \xrightarrow{\sim} \phi^* B$.*

(ii) *The $\mathbb{I}$-strong unit axiom holds, $\mathcal{O}, \mathcal{P} \in \{\Sigma$-cofibrant or well-pointed operads$\}$, and (a) or (b) above is satisfied.*

(iii) *The operads $\mathcal{O}$, $\mathcal{P}$ are cofibrant in $\mathcal{O}\mathrm{pd}_{O}(\mathcal{V})$, $B$ is* fgt*-cofibrant and the map of $\mathcal{O}$-algebras is the identity $\phi^* B \to \phi^* B$.*

*Proof.* Same argument as in Proposition 4.49, but applying Proposition 3.40 and Corollary 3.42 instead of Proposition 3.38. ☐

**Related work:** In [5], it was obtained the first item of Proposition 4.49 when $\mathcal{O}$ is well-pointed. To the best of our knowledge, the rest of Proposition 4.49 and the entire Proposition 4.50 are new.

## 4.6 Operadic vs categorical Lan

From a map of operads $\phi\colon \mathcal{O} \to \mathcal{P}$ (possibly not the identity on colors), we obtain a commutative square of Quillen pairs

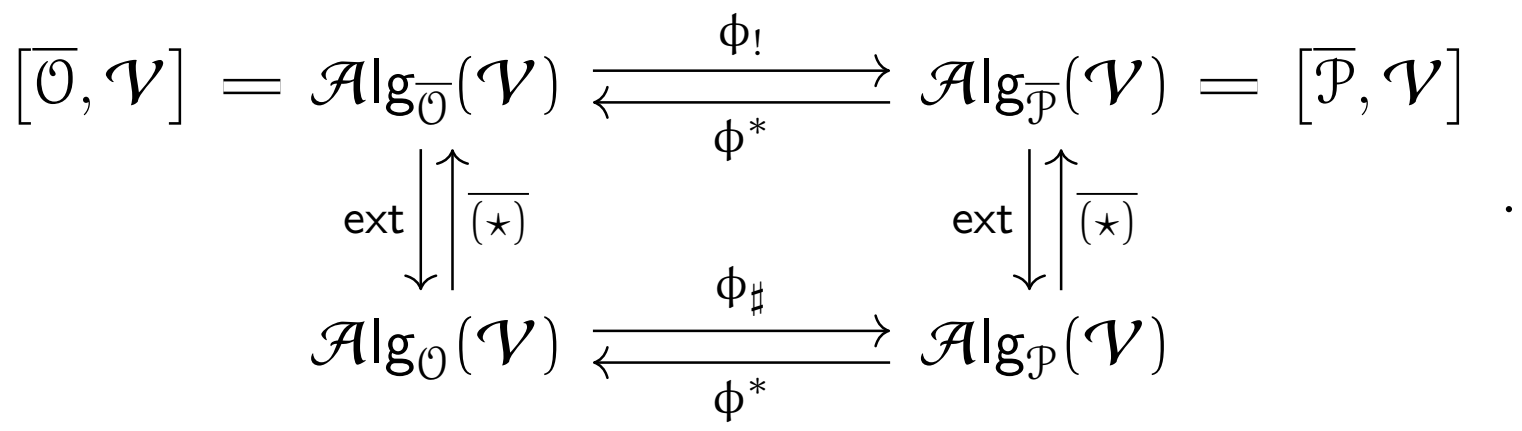

Deriving these functors, we get a natural transformation $\mathbb{L}\phi_!\overline{A} \to \overline{\mathbb{L}\phi_\sharp A}$, called *mate transformation*, by means of the pasting diagram

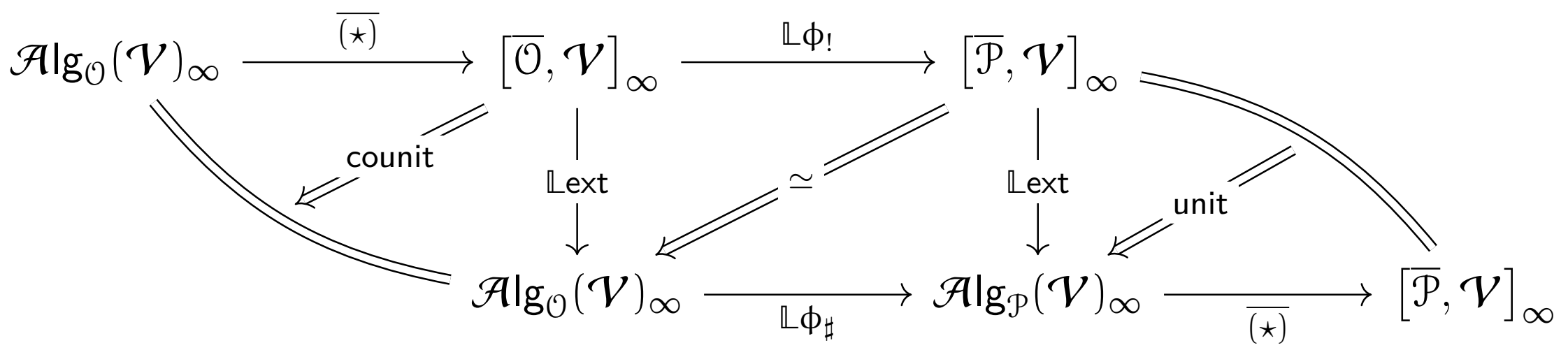

Our goal is to use the enhanced filtrations for the enveloping operad to provide a criterion establishing that $\mathbb{L}\phi_!\overline{A} \to \overline{\mathbb{L}\phi_\sharp A}$ is an equivalence. Note that when A is proj-cofibrant, such canonical comparison map can be computed in $\mathcal{A}\mathrm{lg}_{\overline{\mathcal{P}}}(\mathcal{V})$ as the composite $\phi_!\overline{Q}\,\overline{A} \to \phi_!\overline{A} \to \overline{\phi_\sharp A}$, with $\overline{Q}$ a cofibrant replacement in $\mathcal{A}\mathrm{lg}_{\overline{\mathcal{O}}}(\mathcal{V})$ and the second map being the underived mate transformation.

*Remark* 4.51. It is important to note that the mate natural transformation interwines between the unit and counit of the categorical Lan, $\mathbb{L}\phi_! \dashv \phi^*$, and those of the operadic Lan, $\mathbb{L}\phi_\sharp \dashv \phi^*$. More concretely, commute in $[\overline{\mathcal{O}}, \mathcal{V}]_\infty$ and $[\overline{\mathcal{P}}, \mathcal{V}]_\infty$ respectively. To prove this fact one just have to take into account some more pasting diagrams and triangular identities.

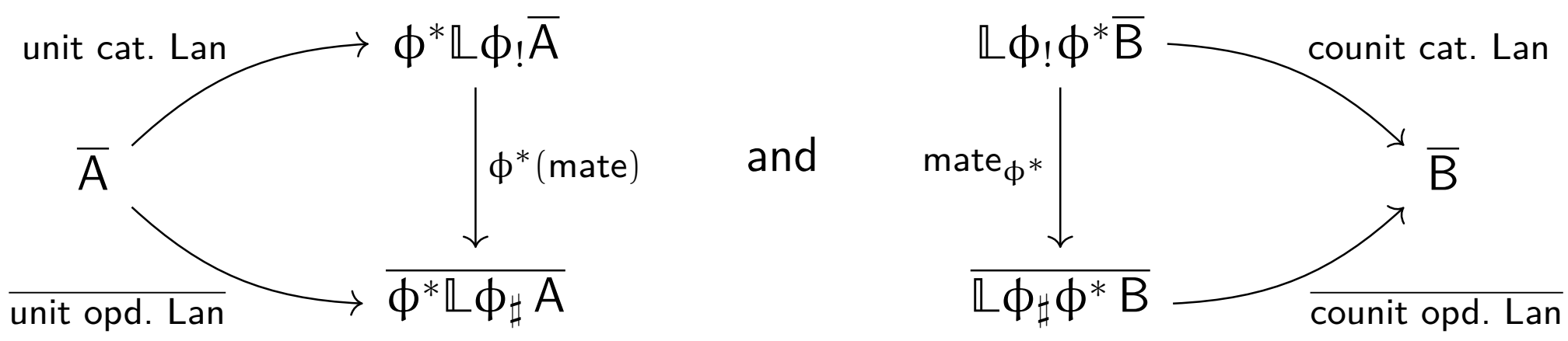

The main result in this subsection is:

**Theorem 4.52.** *Assume that one of the following conditions is satisfied:*

- $\mathcal{O}, \mathcal{P}$ *are* $\Sigma$*-cofibrant operads, or*
- *The* $\mathbb{I}$*-strong unit axiom holds and* $\mathcal{O}, \mathcal{P} \in \{\Sigma$*-cofibrant or well-pointed operads*$\}$.

*Then, the mate natural transformation $\mathbb{L}\phi_! \overline{A} \to \overline{\mathbb{L}\phi_\sharp A}$ is an equivalence for any $\mathcal{O}$-algebra $A$ if*

$$\mathbb{L}\phi_!\,\mathcal{O}\begin{bmatrix}\underline{c}\\ \star\end{bmatrix} \longrightarrow \mathcal{P}\begin{bmatrix}\phi(\underline{c})\\ \star\end{bmatrix} \tag{4.6}$$

*is an equivalence in $[\overline{\mathcal{P}}, \mathcal{V}]$ for all $\underline{c} \in \mathsf{Fin}_{O}^{\simeq}$.*

The proof of this result uses that the mate natural transformation can be computed on proj-cofibrant algebras as $\mathbb{L}\phi_! \overline{A} \to \overline{\phi_\sharp A}$ and that it suffices to check the claim for proj-cellular $\mathcal{O}$-algebras. These objects are constructed inductively by cell attachments, so the following three lemmas control how the construction of cellular algebras allows us to go from the condition in Theorem 4.52 to the claim for general proj-cellular algebras.

*Remark* 4.53. Note that it is possible to assume that proj-cellular algebras are constructed by gluing core-cofibrant cells, as we will assume in Lemma 4.55, because of Remark A.9.

**Lemma 4.54** (Free algebras)**.** *Assume the conditions of Theorem 4.52 hold, in particular (4.6) is an equivalence, then the mate natural transformation is an equivalence on algebras of the form $\mathcal{O} \circ X$ with $X \in \mathcal{V}^{O}$ objectwise-cofibrant.*

*Proof.* Let $X$ be a cofibrant object in $\mathcal{V}$ and $c$ be a color of $\mathcal{O}$. With these data, we can form the collection $X_{(c)}$ in $\mathcal{V}^{O}$, which is concentrated on the color $c \in O$. The free $\mathcal{O}$-algebra on $X_{(c)}$ is

$$\mathcal{O} \circ X_{(c)} = \coprod_{r \geqslant 0} \mathcal{O}\begin{bmatrix} c^{\boxplus r}\\ \star\end{bmatrix} \underset{\Sigma_r}{\otimes} X^{\otimes r}.$$

The mate natural transformation on this free $\mathcal{O}$-algebra becomes

$$\coprod_{r\geqslant 0} \mathbb{L}\phi_!\,\mathcal{O}\begin{bmatrix} c^{\boxplus r}\\ \star\end{bmatrix} \underset{\Sigma_r}{\otimes} X^{\otimes r} \longrightarrow \coprod_{r\geqslant 0} \mathcal{P}\begin{bmatrix} \phi(c)^{\boxplus r}\\ \star\end{bmatrix} \underset{\Sigma_r}{\otimes} X^{\otimes r},$$

where we have applied the commutation of the (homotopy) coproduct with $\mathbb{L}\phi_!$ on the left (which is possible since this is just a derived categorical Lan) and the natural isomorphism $\phi_\sharp \cdot (\mathcal{O} \circ \star) \cong (\mathcal{P} \circ \star) \cdot \phi_!$. The map (4.6) being an equivalence tells us that the mate natural transformation in this case is an equivalence because $X$ is cofibrant and since $\mathcal{O}$ and $\mathcal{P}$ are $\Sigma$-cofibrant. Note that the $\Sigma_r$-quotients here are homotopical in this case and so $\mathbb{L}\phi_!$ commutes with them. In the well-pointed case, use additionally Lemma A.17, to ensure that the arguments above work. □

**Lemma 4.55** (Cell attachments)**.** *Assume the conditions of Theorem 4.52 are met. Let*

$$\begin{array}{ccc} \mathcal{O}\circ X & \longrightarrow & A \\ {\scriptstyle \mathcal{O}\circ j}\downarrow & \ulcorner & \downarrow{\scriptstyle g} \\ \mathcal{O}\circ Y & \longrightarrow & A[j] \end{array}$$

*be a pushout square in $\mathcal{A}\mathsf{lg}_{\mathcal{O}}(\mathcal{V})$ with $j$ core-cofibration in $\mathcal{V}^{O}$ and $A$ proj-cofibrant. Then, the mate natural transformation is an equivalence on $A[j]$ if it is so on $A$.*

*Proof.* We argue the $\Sigma$-cofibrant case. The well-pointed one follows the same lines, but also making use of Lemma A.17. The mate natural transformation evaluated on the pushout in the statement yields a commutative cube

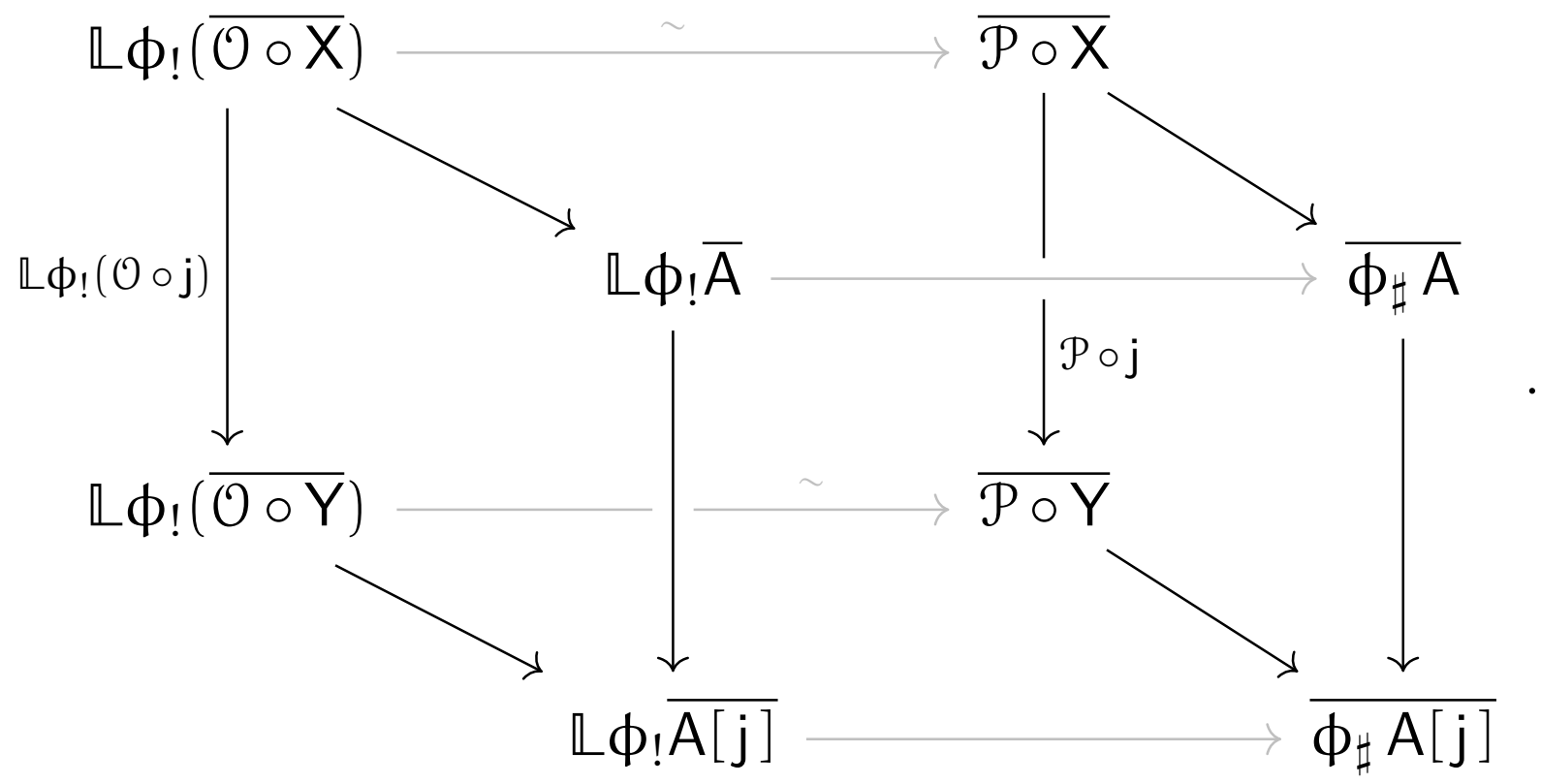

Note that the horizontal maps in the back face are equivalences by Lemma 4.54. It is not true that the black faces are homotopy pushouts in $[\mathcal{P},\mathcal{V}]$, so we cannot apply the gluing lemma for model categories (the forgetful functor $\mathcal{A}lg_{\mathcal{P}}(\mathcal{V}) \to [\mathcal{P},\mathcal{V}]$ does not preserve general (homotopy) colimits, e.g. (ho)pushouts). However, we can filter the front square with the results in §3, see more precisely Lemmas 3.20-3.24, and deduce the claim by an induction argument over that filtration.

Using the pushout square in the statement, one gets a filtration of $\overline{A} \to \overline{A[j]}$ in $[\mathcal{O},\mathcal{V}]$, and so a filtration of $\mathbb{L}\phi_!\overline{A} \to \mathbb{L}\phi_!\overline{A[j]}$. Analogously, applying $\phi_\sharp$ to the original pushout square, one finds a filtration of $\overline{\phi_\sharp A} \to \overline{\phi_\sharp A[j]}$ in $[\mathcal{P},\mathcal{V}]$. Both filtrations may be compared via a $\omega$-ladder (see Lemma 3.24)

$$\begin{array}{ccccccccc}
\mathbb{L}\phi_!\overline{A} & \longrightarrow & \mathbb{L}\phi_!\overline{A[j]}_1 & \longrightarrow & \mathbb{L}\phi_!\overline{A[j]}_2 & \longrightarrow & \cdots & \longrightarrow & \mathbb{L}\phi_!\overline{A[j]} \\
\downarrow & & \downarrow & & \downarrow & & & & \downarrow \\
\overline{\phi_\sharp A} & \longrightarrow & \overline{\phi_\sharp A[j]}_1 & \longrightarrow & \overline{\phi_\sharp A[j]}_2 & \longrightarrow & \cdots & \longrightarrow & \overline{\phi_\sharp A[j]},
\end{array}$$

where each successive step $(t-1) \Rightarrow (t)$ fits into a cube[ix]

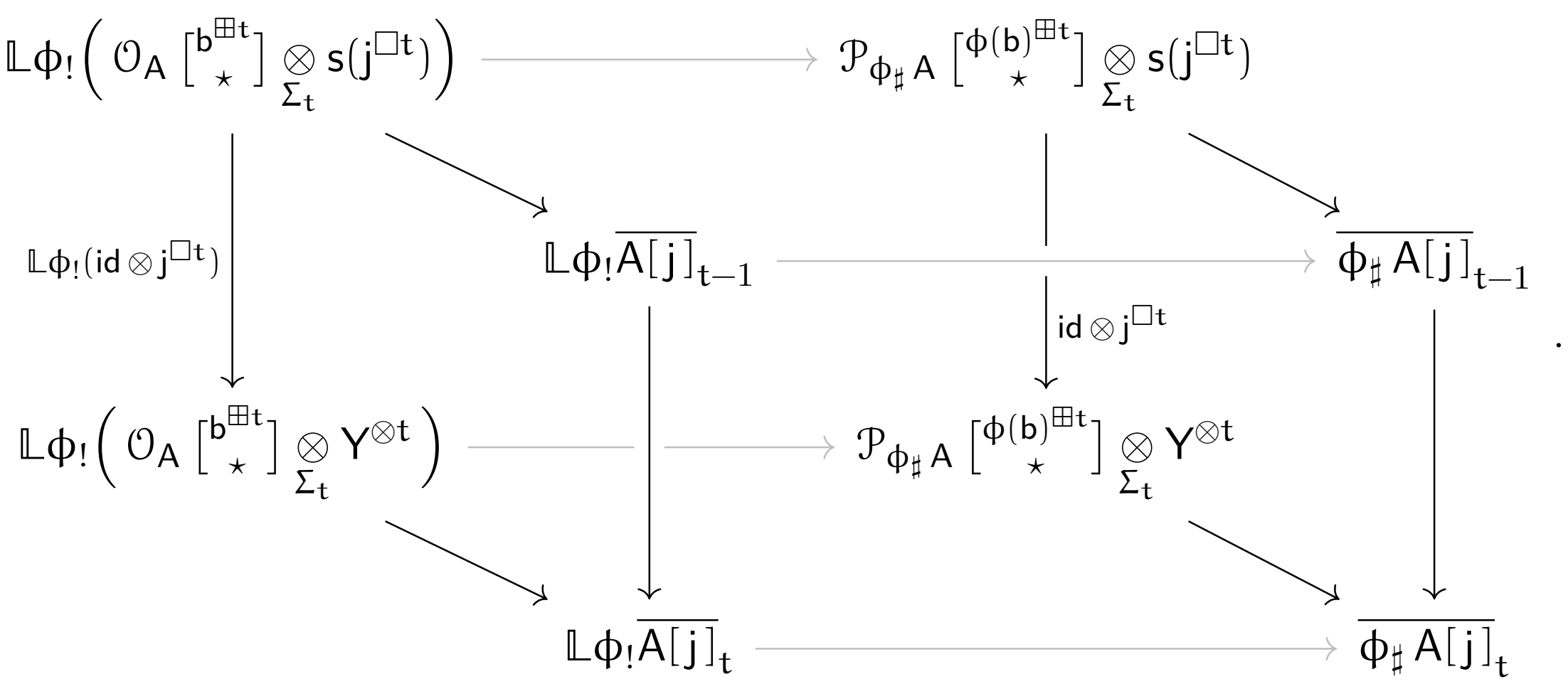

[ix]For notational simplicity, we have assumed that j is concentrated on one color $b \in O$.

Now, the black faces of the cube are homotopy pushouts in $[\mathcal{P}, \mathcal{V}]$. By Proposition 3.34, it suffices to check that

$$\mathbb{L}\phi_! \, \mathcal{O}_A \begin{bmatrix} b^{\boxplus t} \\ \star \end{bmatrix} \longrightarrow \mathcal{P}_{\phi_\sharp A} \begin{bmatrix} \phi(b)^{\boxplus t} \\ \star \end{bmatrix} \tag{4.7}$$

is an objectwise-equivalence to deduce that the horizontal maps in the back square are equivalences. By induction on $t$, we obtain that all vertical maps in the ladder are equivalences.

We may assume without loss of generality that $A$ is proj-cellular and so it comes with a cellular filtration

$$\mathbb{O}_{\mathcal{O}} = A_0 \to \cdots \to A_\alpha \to A_{\alpha+1} \to \cdots \to A_\kappa = A,$$

that is, each step $(\alpha) \Rightarrow (\alpha+1)$ is a cobase change of a map $\mathcal{O} \circ j_\alpha$ with $j_\alpha$ core cofibration in $\mathcal{V}^O$. Running one induction for $\kappa$ and another induction over $\omega$ (for each $A_\alpha \to A_{\alpha+1}$; similar to that of Lemma 3.41), one proves that (4.7) is an equivalence if (4.6) is an equivalence. Just for concreteness, note that the filtration over $\omega$ for $(\alpha) \Rightarrow (\alpha+1)$ yields in step $(m-1) \Rightarrow (m)$ a cube[x] whose black faces are (ho)pushouts. Note that here we are using Lemma 3.22 to run the induction over enveloping operads. □

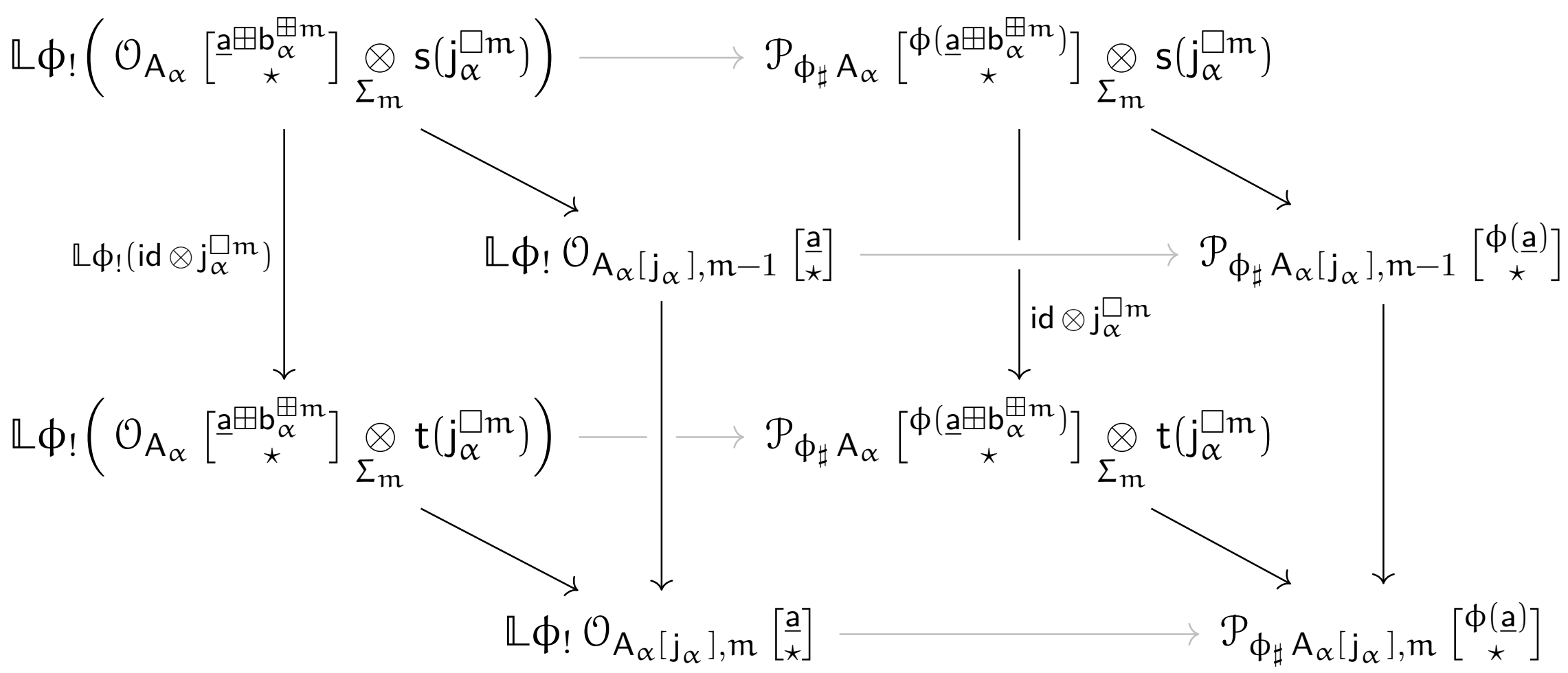

**Lemma 4.56** (Transfinite colimit). *Assume the conditions of Theorem 4.52 are met. Let $\lambda$ be a regular cardinal and*

$$A_0 \to A_1 \to \cdots A_\alpha \to A_{\alpha+1} \to \cdots \to \operatorname*{colim}_{\alpha<\lambda} A_\alpha = A_\lambda$$

*be a $\lambda$-indexed diagram of proj-cofibrations in $\mathcal{A}lg_{\mathcal{O}}(\mathcal{V})$ which is colimit preserving (transfinite composition) and with $A_0$ proj-cofibrant. Then, the mate natural transformation is an equivalence on $A_\lambda$ provided it is so on $A_\alpha$ for any $\alpha < \lambda$.*

[x]Again, we have assumed for ease of notation that $j_\alpha$ is concentrated on color $b_\alpha \in O$.

*Proof.* The $\lambda$-indexed diagram gives rise to a commutative diagram

$$\begin{array}{ccccccccc}
\mathbb{L}\phi_!\overline{A}_0 & \longrightarrow & \cdots & \longrightarrow & \mathbb{L}\phi_!\overline{A}_\alpha & \longrightarrow & \mathbb{L}\phi_!\overline{A}_{\alpha+1} & \longrightarrow & \cdots & \longrightarrow & \mathbb{L}\phi_!\overline{A}_\lambda \\
\downarrow & & & & \downarrow & & \downarrow & & & & \downarrow \\
\overline{\mathbb{L}\phi_\sharp A_0} & \longrightarrow & \cdots & \longrightarrow & \overline{\mathbb{L}\phi_\sharp A_\alpha} & \longrightarrow & \overline{\mathbb{L}\phi_\sharp A_{\alpha+1}} & \longrightarrow & \cdots & \longrightarrow & \overline{\mathbb{L}\phi_\sharp A_\lambda}
\end{array},$$

and we have to prove that the right vertical map is an equivalence provided all vertical maps except it are equivalences. The horizontal chains can be taken to be towers of cofibrations whose first object is cofibrant in $\mathcal{V}^O$. Therefore, the horizontal colimits are homotopical and the induced arrow on colimits is an equivalence. One concludes the claim whenever the induced arrow on colimits coincides with the natural transformation $\mathbb{L}\phi_!\overline{A}_\lambda \longrightarrow \overline{\mathbb{L}\phi_\sharp A_\lambda}$. This fact is just a commutation of sequential homotopy colimits with all the functors involved in the definition of the mate natural transformation. □

*Proof of Theorem 4.52.* Since equivalences are closed under retracts, we have to check the claim just for cellular proj-cofibrant $\mathcal{O}$-algebras. This claim is consequently reduced to show that it holds in each step in the construction of a cellular algebra. The free case is treated in Lemma 4.54, the pushout case in Lemma 4.55 and the transfinite one in Lemma 4.56. □

*Remark* 4.57. To the best of our knowledge, Theorem 4.52 is the first result of this form in the literature.

# A Miscellanea

**Semimodel categories.** We will adopt the following notation:

**Definition A.1.** A *core cofibration* (resp. *core fibration*) is a cofibration that has cofibrant source (resp. fibrations with fibrant target). We denote the class of core cofibrations (resp. core fibrations) by $\mathtt{Cof}_\circ$ (resp. $\mathtt{Fib}_\circ$). Also, $\mathtt{ACof}_\circ \equiv \mathtt{Cof}_\circ \cap \mathtt{Eq}$ and $\mathtt{AFib}_\circ \equiv \mathtt{Fib}_\circ \cap \mathtt{Eq}$.

Let us recall the notion of semimodel structure that we use in this document. See [1, 8, 10] for more.

**Definition A.2.** A *structured homotopical category* $\mathcal{M}$ is a category $\mathtt{M}$ endowed with three classes of maps $(\mathtt{Cof}, \mathtt{Eq}, \mathtt{Fib})$ which satisfy:

- $\mathtt{M}$ is bicomplete, i.e. has all limits and colimits.
- $\mathtt{Eq}$ is closed under retracts and satisfies 2-out of-3.
- $\mathtt{Cof}$ (resp. $\mathtt{Fib}$) is closed under retracts and pushouts (resp. pullbacks).

**Definition A.3.** A *(Spitzweck) left semimodel category* $\mathcal{M}$ is a structured homotopical category $(\mathtt{M}, \mathtt{Cof}, \mathtt{Eq}, \mathtt{Fib})$ which satisfies: $(a)$ $(\mathtt{Cof}, \mathtt{AFib})$ is a weak factorization system (wfs), while $(b)$ the pair $(\mathtt{ACof}, \mathtt{Fib})$ only satisfies the following relaxed conditions:

$(b1)$ (Core lifting axiom) $\mathtt{ACof}_\circ \boxslash \mathtt{Fib}$.

(b2) (Core factorization axiom) any map f with cofibrant domain can be factored as

$$f\colon \bullet \xrightarrow{\texttt{ACof}_\circ} \bullet \xrightarrow{\texttt{Fib}} \bullet.$$

*Remark* A.4. The adjective "Spitzweck" in Definition A.3 refers to a choice that has been made (see [8]). In this document, we omit this adjective since all the semimodel structures that we consider are of this form.

**Relative cell complexes.**

**Definition A.5.** Let K be a set of maps in a category C.

- We say that a morphism $X \to X'$ is a K-cell attachment if it fits into a pushout square

$$\begin{array}{ccc} \bullet & \xrightarrow{\coprod_i k_i} & \bullet \\ \downarrow & \ulcorner & \downarrow \\ X & \longrightarrow & X' \end{array}$$

  where $k_i \in K$ for all $i$.

- A map $f\colon X \to X'$ is a *relative* K-*cell complex*, or equivalently $f \in K$-cell, if it can be obtained by a (possibly transfinite) composition[i]

$$f\colon X = X_0 \to \cdots \to X_\alpha \xrightarrow{f_\alpha} X_{\alpha+1} \to \cdots \to X_\mu = X'$$

  over an ordinal $\mu$ of K-cell attachments $f_\alpha$'s.

- A K-*cell complex* is an object $X \in C$ such that $\mathbb{0} \to X$ is in K-cell.

*Remark* A.6. If one is willing to consider longer transfinite compositions, each K-cell attachment in a relative K-cell complex can be taken to be a cobase change of just one morphism in K.

We follow the convention

$$\begin{array}{ccccccccc} & & & \bullet & \xrightarrow{\text{cell}_t} & \bullet & & & \\ & & & {\scriptstyle \text{attach}_t}\downarrow & \ulcorner & \downarrow{\scriptstyle \text{char}_t} & & & \\ X = X_0 & \longrightarrow & \cdots \longrightarrow & X_{t-1} & \xrightarrow[\text{bond}_t]{} & X_t & \longrightarrow & \cdots \longrightarrow & X_\mu = X' \end{array}$$

to make reference to the components of a cellular decomposition. The upper horizontal map is the attached cell, the left vertical map in the square is called attaching map, the one on the right, characteristic map and the lower horizontal one, bonding map.

The following lemma and comment are just easy results about relative cell complexes indexed by the first infinite ordinal $\omega$ that we need in §3.

[i]As always assumed for transf.compositions, for any limit ordinal $\lambda \leqslant \mu$, one sets $X_\lambda = \operatorname{colim}_{\beta<\lambda} X_\beta$.

**Lemma A.7.** *Assume that* $\mathsf{C}$ *has a (symmetric) monoidal structure. If* $X_0 \to X_\omega$ *and* $Y_0 \to Y_\omega$ *are relative* $\mathsf{K}$*-cell* $\omega$*-complexes, then* $X_0 \otimes Y_0 \to X_\omega \otimes Y_\omega$ *is a relative cell* $\omega$*-complex with:*

- *stages* $(X \otimes Y)_t = \mathrm{colim}_{p+q \leqslant t}\, X_p \otimes Y_q$*;*
- *cells* $\mathrm{cell}_{X\otimes Y,t} = \coprod_{p+q=t} \mathrm{cell}_{X,p} \,\square\, \mathrm{cell}_{Y,q}$*;*
- *characteristic maps induced by* $\mathrm{char}_{X,p} \otimes \mathrm{char}_{Y,q}$*; and*
- *attaching maps induced by* $\mathrm{attach}_{X,p} \otimes \mathrm{char}_{Y,q}$ *and* $\mathrm{char}_{X,p} \otimes \mathrm{attach}_{Y,q}$.

*Proof.* Rutinary arguments over the commutative diagram

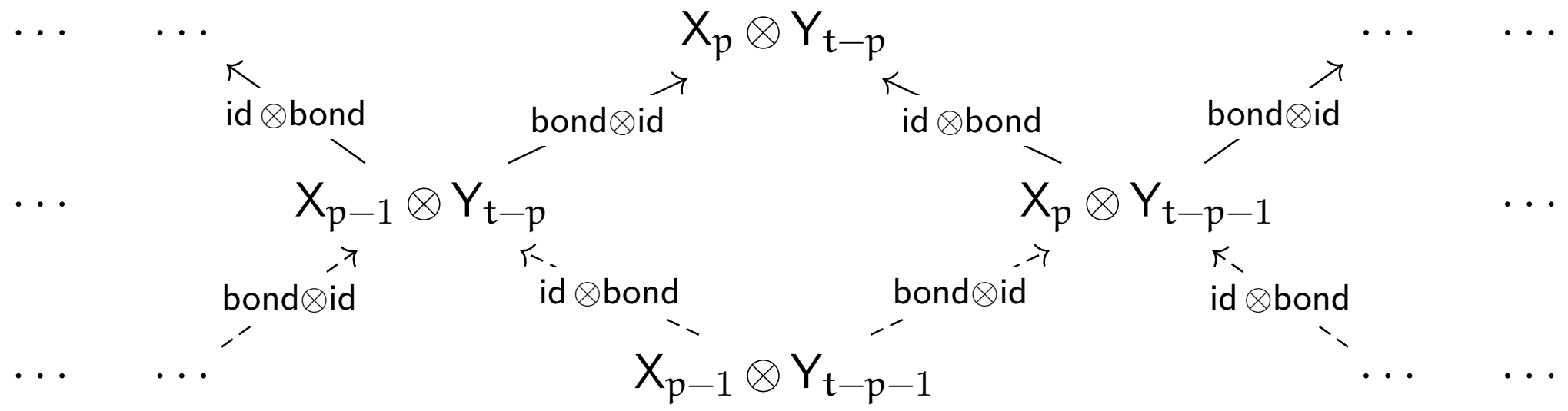

such that the colimit of the solid (resp. dashed) subdiagram is the stage $(X \otimes Y)_t$ (resp. $(X \otimes Y)_{t-1}$) above. ☐

*Remark* A.8. Any relative cell $\omega$-complex $X_0 \to X_\omega$ admits a *twisted structure*, that is, a modified relative $\omega$-cell complex structure where one modifies the cells, but keeping the same stages and bonding maps. Just consider

$$\begin{array}{ccccc} & \bullet & \xrightarrow{\mathrm{cell}_t \amalg \mathrm{id}_{X_{t-1}}} & \bullet & \\ & {\scriptstyle (\mathrm{attach}_t,\, \mathrm{id}_{X_{t-1}})}\downarrow & \ulcorner & \downarrow{\scriptstyle (\mathrm{char}_t,\, \mathrm{bond}_t)} & \\ X_0 \longrightarrow \cdots \longrightarrow & X_{t-1} & \xrightarrow[\mathrm{bond}_t]{} & X_t & \longrightarrow \cdots \longrightarrow X_\omega \end{array} .$$

**Some useful tricks.** Let us collect some direct, while important, observations that we apply in several places. The first one is discussed in the following:

*Remark* A.9. In some situations, it is helpful to add cells which are core (acyclic) cofibrations to construct cellular objects, e.g. when dealing with Bousfield localizations or to show homotopical invariance and/or cofibrancy properties of operadic constructions.

1. If the model structure is tractable, i.e. core cofibrantly generated, all the (acyclic) cells that one might want to attach are core cofibrations.

2. If the model structure is left proper and we consider a cell attachment

$$\begin{array}{ccc} \bullet & \longrightarrow & X \\ {\scriptstyle j}\downarrow & \ulcorner & \downarrow \\ \bullet & \longrightarrow & X[j] \end{array} ,$$

by factoring maps using the wfs $(\mathtt{Cof}, \mathtt{AFib})$, i.e. finding a core replacement $j_c$ of $j$, we form a cube whose black faces are pushouts. By left properness, these pushouts are homotopy pushouts and so $X[j_c] \to X[j]$ is an equivalence. Using this fact, one may replace each cell attachment by a core cell attachment with an equivalent result.

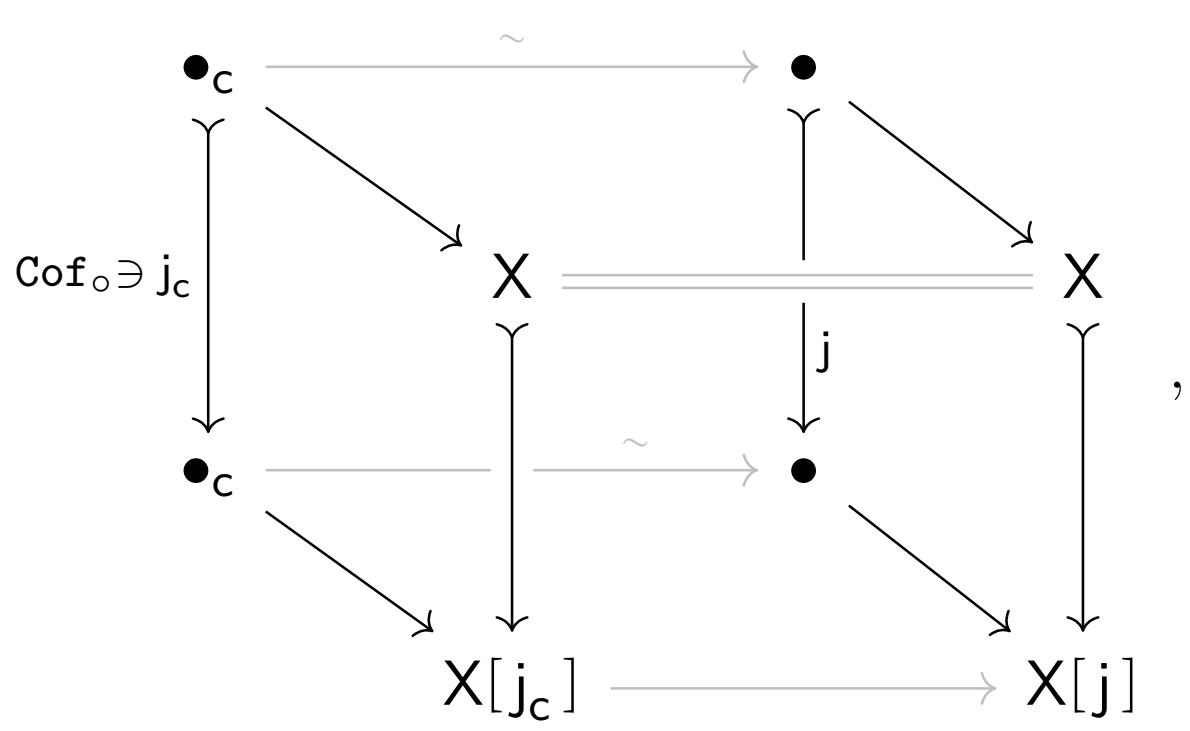

3. A trick for left-transferred (semi)model structures [26, Proof of Proposition 4.2]. Consider a Quillen adjunction $F \dashv \mathrm{fgt}$. If we attach a free cell $F(j)$ to a fgt-cofibrant object $A$, we may add a different core cell to $A$ yielding the same result: starting from

$$\begin{array}{ccc} F(X) & \longrightarrow & A \\ {\scriptstyle F(j)}\downarrow & \ulcorner & \downarrow \\ F(Y) & \longrightarrow & A' \end{array},$$

apply fgt and consider

$$\begin{array}{ccccc} X & \longrightarrow & \mathrm{fgt}(A) & = & \mathrm{fgt}(A) \\ {\scriptstyle j}\downarrow & \ulcorner & \downarrow{\scriptstyle k \in \mathtt{Cof}_\circ} & & \downarrow \\ Y & \longrightarrow & Z & \longrightarrow & \mathrm{fgt}(A') \end{array}.$$

Then, applying F, one gets

$$\begin{array}{ccccc} F(X) & \longrightarrow & F(\mathrm{fgt}(A)) & \longrightarrow & A \\ {\scriptstyle F(j)}\downarrow & \ulcorner & \downarrow{\scriptstyle F(k) \in \mathtt{Cof}_\circ} & & \downarrow \\ F(Y) & \longrightarrow & F(Z) & \longrightarrow & A' \end{array}$$

where the full rectangle is a pushout. Since the left square is also a pushout, the right square is a pushout and it is a cobase change of a core cofibration as we wanted.

The last remark concerns the perspective on unitality as structure or as a property. Looking at unitality as a property simplifies some computations and constructions in (homotopical) algebra. For instance, one has the following simple result in this direction.

**Proposition A.10.** *Let* $\mathtt{U}$ *be a monoidal category. Then, the functor that forgets the identity element* $\mathtt{Mon}(\mathtt{U}) \to \mathtt{Semigp}(\mathtt{U})$, $M \mapsto M^\dagger$ *creates pushouts and sequential colimits (when they exist). More precisely, let* $\mathscr{D}$ *be a diagram of one of these classes in* $\mathtt{Mon}(\mathtt{U})$ *and* $\mathscr{D}^\dagger \Rightarrow M$ *be a colimit cocone in* $\mathtt{Semigp}(\mathtt{U})$. *Then, the cocone lives in* $\mathtt{Mon}(\mathtt{U})$ *and is colimiting there.*

*Proof.* To endow $M$ with a (two-sided) unit, consider the composition

$$1_M \colon \mathbb{1}_U \xrightarrow{\text{unit of } D_0} D_0 \xrightarrow[\text{cocone}]{\text{component of}} M \ ,$$

where $D_0$ is the initial object in $\mathscr{D}$. All maps in the cocone $\mathscr{D}^\dagger \Rightarrow M$ preserve units (equipping $M$ with $1_M$) since $D_0$ is initial in $\mathscr{D}$. This fact together with the universal property of $M$ in $\mathtt{Semigp}(\mathtt{U})$ implies that $1_M$ is in fact a two-sided unit of $M$. Hence, we have checked that the cocone lives in $\mathtt{Mon}(\mathtt{U})$.

To show that the cocone is colimiting in this category, take another cocone $\mathscr{D} \Rightarrow N$ in monoids. Then, we get a unique map $M \to N^\dagger$ compatible with the product and it remains to check that it preserves units. This follows from the commutative diagram

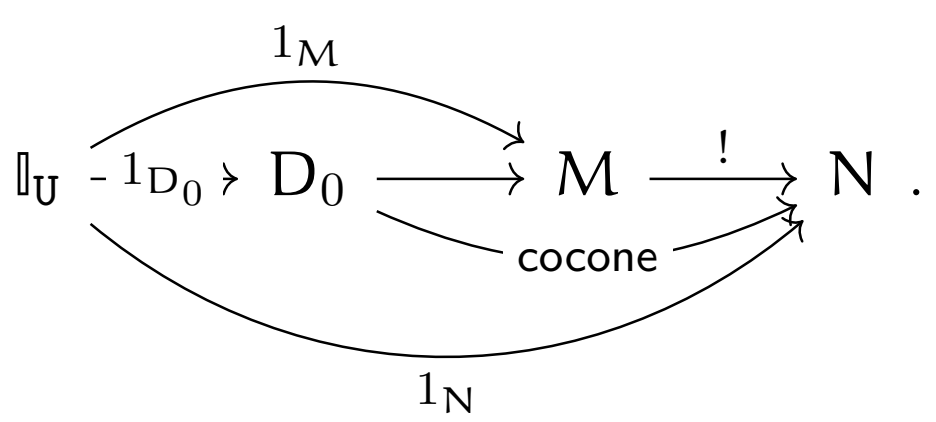

□

Since categories and operads (with a fixed set of objects/colors) are monoids for certain monoidal categories, this result applies to both of them. However, the idea is quite general and regarding operads, we prefer to consider a different category of "operads without identities". Using a colored version of Markl's pseudo-operads (see [21]) one has an analogous:

**Proposition A.11.** *The functor that forgets identities*

$$\mathtt{Opd}_O(\mathtt{V}) \to \{O\text{-colored pseudo-operads in } \mathtt{V}\}$$

*creates pushouts and sequential colimits.*

Note that even more general colimits are created by forgetting units/identities, but we prefer less generality here since this is all what we will need: pushouts and sequential colimits are the ingredients to construct cellular objects.

**$\mathbb{1}$-cofibrancy.** We finish this appendix by presenting some technical observations applied several times in the body of the text to deal with monoidal units that are not cofibrant. The following material is based on [26, Appendix A, B]. Let us fix a general (possibly non-symmetric nor closed) monoidal model category $(\mathcal{C}, \otimes, \mathbb{1})$ in this digression. Notice that we do not assume the monoid axiom on $\mathcal{C}$, in contrast to [25, 26, 29].

**Definition A.12.** [B.1 in [26]] An object $X \in \mathcal{C}$ is $\mathbb{I}$*-cofibrant* if there exists a cofibration

$$\mathbb{I} \rightarrowtail X.$$

*Example* A.13. O-colored well-pointed operads in $\mathcal{V}$ are $\mathbb{I}$-cofibrant objects in the model category $(\Sigma\mathcal{C}oll_O(\mathcal{V}), \circ, \mathcal{J}_O)$.

*Remark* A.14. As a consequence of the pushout-product axiom, if $X$ is $\mathbb{I}$-cofibrant, $X \otimes \star \colon \mathcal{C} \to \mathcal{C}$ is left-Quillen. Moreover, there are also simple closure properties of $\mathbb{I}$-cofibrant objects which follow from this axiom (see [26, A.4-A.8]).

The following technical axiom is key to obtain a nice homotopical behavior of $\mathbb{I}$-cofibrant objects. It is an obvious generalization of [26, Definition A.9].

**Definition A.15.** We say that $\mathcal{C}$ satisfies the $\mathbb{I}$*-strong unit axiom*, or that such an axiom holds, if the following is satisfied: If $X$ is an $\mathbb{I}$-cofibrant object and $q \colon Q\mathbb{I} \to \mathbb{I}$ is a cofibrant replacement of the monoidal unit, the maps $X \otimes q$ and $q \otimes X$ are equivalences.

*Remark* A.16. If this axiom holds for a certain cofibrant replacement of $\mathbb{I}$, then it holds for any cofibrant replacement of $\mathbb{I}$. In particular, it is satisfied if: $(a)$ the unit $\mathbb{I}$ is cofibrant; or $(b)$ cofibrant objects are flat.

**Lemma A.17.** *Assume $\mathcal{C}$ satisfies the $\mathbb{I}$-strong unit axiom. Then,*

- *Equivalences between cofibrant or $\mathbb{I}$-cofibrant objects are closed under tensor products $\otimes$ and arbitrary coproducts $\amalg$.*

- *Consider a commutative diagram in $\mathcal{C}$*

$$\begin{array}{ccccc} & & X & & \\ & \swarrow & \downarrow\wr & \searrow & \\ Y & & & & Z \\ \downarrow\wr & & X' & & \downarrow\wr \\ & \swarrow & & \searrow & \\ Y' & & & & Z' \end{array}.$$

  *If each of the objects is $\mathbb{I}$-cofibrant or cofibrant, the induced map between pushouts is an equivalence (we will say that these pushouts are homotopical).*

- *Consider a limit ordinal $\lambda$, two cocontinuous functors $X_\bullet, Y_\bullet \colon \lambda \to \mathcal{C}$ and a natural equivalence $\psi_\bullet \colon X_\bullet \Rightarrow Y_\bullet$, i.e. $\psi_\beta \colon X_\beta \xrightarrow{\sim} Y_\beta$ for all $\beta < \lambda$. Then,*

$$\operatorname*{colim}_{\beta<\lambda} \psi_\beta \colon X_\lambda \longrightarrow Y_\lambda$$

  *is an equivalence provided $X_0$, $Y_0$ are cofibrant or $\mathbb{I}$-cofibrant, and the morphisms $X_\alpha \to X_{\alpha+1}$, $Y_\alpha \to Y_{\alpha+1}$ are cofibrations for all $\alpha < \lambda$.*

*Proof.* Obvious adaptations of [26, Corollary A.14, Lemmas A.15-A.17]. □

## Acknowledgments

The author would like to thank Fernando Muro for his continued support and for all his explanations about the intricacies of his former work on the homotopy theory of non-symmetric operads and algebras.